\documentclass[12pt]{amsart}
\usepackage{amsmath,amsthm,amsfonts,amssymb,xcolor,enumitem,mathscinet}
\usepackage{graphicx}
\usepackage{tikz}
\usetikzlibrary{decorations.pathreplacing,patterns}
\usepackage{pgfplots}
\pgfplotsset{compat=1.16}
\usepackage{mathrsfs}
\usepackage{mathtools}
\usepackage{etoolbox}
\apptocmd{\sloppy}{\hbadness 10000\relax}{}{}
\apptocmd{\sloppy}{\vbadness 10000\relax}{}{}
\usepackage[pdfpagelabels]{hyperref}
\usepackage[letterpaper,margin=1.1in]{geometry}

\definecolor{PTblue}{HTML}{0077BB}
\definecolor{PTteal}{HTML}{009988}
\definecolor{PTmagenta}{HTML}{EE3377}

\newtheorem{theorem}{Theorem}[section]
\newtheorem{lemma}[theorem]{Lemma}
\newtheorem{corollary}[theorem]{Corollary}
\newtheorem{proposition}[theorem]{Proposition}

\theoremstyle{definition}
\newtheorem{definition}[theorem]{Definition}
\newtheorem{convention}[theorem]{Convention}

\theoremstyle{remark}
\newtheorem{remark}[theorem]{Remark}
\newtheorem{example}[theorem]{Example}

\newcommand{\RR}{\mathbb{R}}

\newcommand{\ZZ}{\mathbb{Z}}
\newcommand{\MM}{\mathbb{M}}
\newcommand{\XX}{\mathbb{X}}

\newcommand{\HH}{\mathbb{H}}
\newcommand{\EE}{\mathbb{E}}
\newcommand{\PP}{\mathbb{P}}
\newcommand{\Haus}{\mathcal{H}}
\newcommand{\Pack}{\mathcal{P}}

\newcommand{\diam}{\mathop\mathrm{diam}\nolimits}
\newcommand{\side}{\mathop\mathrm{side}\nolimits}

\newcommand{\res}{\hbox{ {\vrule height .22cm}{\leaders\hrule\hskip.2cm} }} 
\newcommand{\gap}{\mathop\mathrm{gap}\nolimits}
\newcommand{\muse}{\mu_{\mathbf{s}}}
\newcommand{\leaves}{\mathsf{Leaves}}
\newcommand{\Child}{\mathsf{Child}}
\newcommand{\Collar}{\mathsf{Collar}}
\newcommand{\Outer}{\mathsf{Outer}}
\newcommand{\Inner}{\mathsf{Inner}}
\newcommand{\Top}{\mathop\mathsf{Top}}

\numberwithin{figure}{section}
\numberwithin{table}{section}
\numberwithin{equation}{section}

\begin{document}

\title[Rectifiable measures in Ahlfors regular spaces]{Rectifiable measures in Ahlfors regular spaces: endpoint constructions}
\thanks{The author was partially supported by NSF DMS grant 2154047.}
\date{November 10, 2026}
\author{Matthew Badger}
\address{Department of Mathematics\\ University of Connecticut\\ Storrs, CT 06269-1009}
\email{matthew.badger@uconn.edu}

\begin{abstract}
A theorem of Balka and Keleti implies that on any complete separable metric space $\XX$, every locally finite Borel measure $\mu$ of packing dimension strictly less than $m$ is ``Euclidean-like'' insofar as there exist Lipschitz maps $f_i:E_i\subset\RR^m\rightarrow\XX$ such that $\mu(\XX\setminus\bigcup_1^\infty f_i(E_i))=0$. An ongoing challenge in geometric measure theory is to understand when this property, called $m$-rectifiability, holds for measures of packing dimension $m$. Towards this goal, we prove that every complete Ahlfors $q$-regular metric space with $q\geq m$ supports $m$-rectifiable doubling measures of packing dimension $m$. The proof of existence of such measures is new even when $\XX=\RR^n$ and $n>m$. More generally, on every complete Ahlfors $q$-regular metric space with $q>m-1$, we construct $m$-rectifiable, purely $(m-1)$-unrectifiable measures with any prescribed exact Hausdorff and packing dimensions $0\leq s_H\leq s_P\leq q$ with $m-1\leq s_P\leq m$, which are doubling when $s_H>0$. This verifies a conjecture of Badger and Schul. The design of the measures relies on metric cubes with variable scaling, in which the ratio of consecutive side lengths tends to infinity but jumps only on a sparse set of generations, together with Badger and Schul's square packing construction of Lipschitz maps.
\end{abstract}

\maketitle

\setcounter{tocdepth}{1}
\tableofcontents

\section{Introduction}\label{s:intro}

A recent advance in metric geometry is the following theorem of Balka and Keleti, which we might think of as a ``universal'' Lipschitz parameterization theorem insofar as \emph{any} compact metric space $\XX$ can be modeled by sets in \emph{any} compact metric space $\MM$ of strictly ``larger'' dimension.

\begin{theorem}[{\cite[Theorem 1.5]{balka-keleti}}] \label{t:bk} Suppose that $\XX$ (``the target'') and $\MM$ (``the model'') are compact metric spaces. If $\dim_H\MM>\overline{\dim}_M\, \XX$, then there exists a compact set $E\subset\MM$ and a Lipschitz map $f:E\rightarrow \XX$ with $f(E)=\XX$.\end{theorem}

In the statement above, $\dim_H$ denotes Hausdorff dimension and $\overline{\dim}_M$ denotes upper Minkowski dimension; we briefly review these and related notions of dimension in \S\ref{ss:dim}.
By invoking the theorem with $\MM=[0,1]^m$ countably many times, one easily obtains the following application in geometric measure theory, which has not previously appeared in the literature. See \S\ref{ss:bk} for a detailed derivation.

\begin{corollary}\label{c:bk} Let $\XX$ be a complete, separable metric space and let $\mu$ be a Borel measure on $\XX$ that is locally finite in the sense that for every $x\in \XX$, there exists $r>0$ such that $\mu(B(x,r))<\infty$. For every integer $m\geq 1$, the measure $\mu\res A$,  \begin{equation}\label{e:small-upper-dimension}A:=\left\{x\in \XX:\limsup_{r\downarrow 0}\frac{\log \mu(B(x,r))}{\log r}<m\right\},\end{equation} is $m$-rectifiable in the sense of Federer, i.e.~there exist countably many Lipschitz maps $f_i:E_i\subset[0,1]^m\rightarrow\XX$ with $\mu(A\setminus\bigcup_i f_i(E_i))=0$. In fact, $A$ is $\sigma$-compact.\end{corollary}

The notion of rectifiability of measures appeared early in the development of geometric measure theory \cite{MR44,Fed47}, but was primarily studied (see \cite{rectifiability-survey} for a comprehensive survey) within the class of locally finite measures $\mu$ on $\XX=\RR^n$ such that \begin{equation}\label{e:good-upper-density} \limsup_{r\downarrow 0} \frac{\mu(B(x,r))}{r^m}<\infty\quad\text{at $\mu$-a.e.~$x\in\XX$,}\end{equation} i.e.~measures that vanish on sets of zero $m$-dimensional Hausdorff measure: $\mu\ll\Haus^m$. For details on the equivalence between \eqref{e:good-upper-density} and the absolute continuity property, see e.g.~\cite[\S2.2]{BLZ}. In particular, classical theorems about rectifiable measures, including those in Federer's monograph \cite{Federer} and in Preiss \cite{Preiss}, assume \eqref{e:good-upper-density} in addition to rectifiability. The first positive results on rectifiability in a class of measures with \begin{equation}\label{e:infinite-upper-density} \limsup_{r\downarrow 0} \frac{\mu(B(x,r))}{r^m}=\infty\quad\text{at $\mu$-a.e.~$x\in\XX$},\end{equation} i.e.~measures that are carried by $\Haus^m$ null sets (equivalently: $\mu\perp \Haus^m$), were obtained by Mart\'in and Mattila \cite{MM1988}. Corollary \ref{c:bk} to Balka and Keleti's theorem is a substantial extension of \cite[Theorem 4.1(1)]{MM1988} as well as \cite[Theorem C]{BV}. For further results on rectifiability of general measures, see \cite{ident,Badger-Naples}.

It is an open problem to give a pointwise characterization of the $m$-rectifiable Radon measures on $\RR^n$ when $2\leq m\leq n-1$. The case $m=1$ was solved in \cite{BS3} in terms of anisotropic $L^2$ Jones $\beta$ numbers (also see \cite{AM-rectifiable}). Corollary \ref{c:bk} shows that the search for such a characterization may focus on the part of the measure where the upper local dimension is equal to $m$. In other words, a characterization will not depend on any exotic geometry of sets of packing dimension less than $m$.

In this paper, we continue a study begun in \cite{BS4} and complete the proof of the following theorem, which asserts the existence of rectifiable measures with prescribed lower and upper local dimensions. We are particularly interested in constructions of measures whose rectifiability cannot be verified using Corollary \ref{c:bk}. Instead, we will demonstrate rectifiability using Badger and Schul's ``square packing construction'' of Lipschitz maps with Euclidean domains and metric targets, which we recall in \S\ref{ss:spc}. A metric space $\XX$ is \emph{Ahlfors $q$-regular} with \emph{distortion} $D_\nu\geq 1$ if $\diam\XX>0$ and there exist a Borel measure $\nu$ on $\XX$ and a constant $C_\nu>0$ such that $C_\nu r^q\leq \nu(B(x,r))\leq C_\nu D_\nu r^{q}$ for all $x\in\XX$ and $0<r\leq\diam \XX$. In particular, $\nu$ is finite on bounded sets.

\begin{theorem}\label{t:main} Let $\XX$ be a complete Ahlfors $q$-regular metric space. For all integers $m\geq 1$ with $q>m-1$ and for all real-valued dimensions $0\leq s_H\leq s_P\leq q$ with $m-1\leq s_P\leq m$, we can find a Radon measure $\mu$ on $\XX$ such that
\begin{equation}\label{full-support} 0<\mu(B(x,r))<\infty\quad\text{for all $x\in\XX$ and $r>0$};\end{equation}
\begin{equation}\label{haus-dim-s}\liminf_{r\downarrow 0} \frac{\log\mu(B(x,r))}{\log r}=s_H\quad\text{at $\mu$-a.e.~$x\in\XX$;}\end{equation}
\begin{equation}\label{pack-dim-s}\limsup_{r\downarrow 0} \frac{\log\mu(B(x,r))}{\log r}=s_P\quad\text{at $\mu$-a.e.~$x\in\XX$;}\end{equation}
\begin{equation}\label{unrectifiable}\text{$\mu(f(E))=0$ whenever $f:E\subset\RR^{m-1}\rightarrow\XX$ is Lipschitz;}\end{equation}
\begin{equation}\label{m-rectifiable}\text{there exist Lipschitz maps $f_i:E_i\subset\RR^m\rightarrow\XX$ such that }\mu\left(\XX\setminus \textstyle\bigcup_1^\infty f_i(E_i)\right)=0.\end{equation} When $s_H>0$, we may further arrange that $\mu$ is doubling in the sense that
\begin{equation}\label{doubling} \mu(B(x,2r))\leq D_{\mathsf{db}}\,\mu(B(x,r))\quad\text{for all $x\in\XX$ and $r>0$}\end{equation} for some constant $1\leq D_{\mathsf{db}}<\infty$ depending on at most $s_H$ and $\XX$.
\end{theorem}

\begin{remark}Properties \eqref{haus-dim-s} and \eqref{pack-dim-s} are equivalent, respectively, to the assertions that the exact Hausdorff dimension $\dim_H\mu$ of $\mu$ is $s_H$ and the exact packing dimension $\dim_P\mu$ of $\mu$ is $s_P$; see e.g.~\cite{tamashiro} or \cite{Falconer-techniques}. In the language of geometric measure theory, equation \eqref{unrectifiable} says that $\mu$ is purely $(m-1)$-unrectifiable; when $m-1=0$, one should interpret this to mean that $\mu(\{x\})=0$ for every $x\in\XX$. Equation \eqref{m-rectifiable} says that $\mu$ is $m$-rectifiable or carried by Lipschitz images of subsets of $[0,1]^m$.\end{remark}

The special case of Theorem \ref{t:main} with parameters $s_H>0$ and $m-1<s_P<\min\{m,q\}$ was established in \cite{BS4}. Thus, the goal is to prove that the theorem remains valid at the endpoints, i.e.~we verify \cite[Conjecture 1.3]{BS4} in the affirmative. The following examples illustrate the upper endpoint $s_P=m$.

\begin{example}
Let $n\geq 2$. Every doubling measure on $\RR^n$ vanishes on porous sets, and in particular on bi-Lipschitz curves (see \cite[\S6.2]{BS4}). Nevertheless, Garnett, Killip, and Schul \cite{GKS} constructed a doubling measure $\nu$ on $\RR^n$ that charges a rectifiable curve, which is necessarily not bi-Lipschitz. This surprising result raised the question of which Radon measures on $\RR^n$ are carried by rectifiable curves, which was later answered in \cite{BS3}. Since $\dim_P\nu\ll 1$, the 1-rectifiability of $\nu$ now also follows from Corollary \ref{c:bk}. Applying Theorem \ref{t:main} with $\XX=\RR^n$ and $s_H=s_P=m=1$, we obtain a 1-rectifiable doubling measure $\mu$ on $\RR^n$ with $\dim_H\mu=\dim_P\mu=1$, which is the largest possible dimension of a 1-rectifiable measure. The measure $\mu$ may be viewed as completing the inquiry of Garnett, Killip, and Schul.
\end{example}

\begin{example}
Let $\XX=\HH^1$ be the Heisenberg group, modeled on $\RR^3$ with a nonabelian group law and a left-invariant metric under which $\HH^1$ is complete and Ahlfors $4$-regular. Ambrosio and Kirchheim \cite{AK} proved that $\HH^1$ is purely $4$-unrectifiable, i.e.~$\Haus^4(f(E))=0$ for every Lipschitz map $f:E\subset\RR^4\rightarrow\HH^1$. By contrast, applying Theorem \ref{t:main} with $\XX=\HH^1$ and $s_H=s_P=m=q=4$, we obtain a $4$-rectifiable doubling measure $\mu$ on $\HH^1$ with $\dim_H\mu=\dim_P\mu=4$. Thus, doubling measures on $\HH^1$ can charge sets of Hausdorff dimension 4 that are Lipschitz images of subsets of $\RR^4$, albeit with $\Haus^4$ measure zero.
\end{example}

We view Corollary \ref{c:bk} and Theorem \ref{t:main} as evidence that Federer was right to omit \eqref{e:good-upper-density} from the definition of rectifiability. These results show that rectifiability is more closely tied to packing dimension than to Hausdorff dimension. The reliance of the classical theorems on \eqref{e:good-upper-density} may simply reflect the fact that Tricot introduced packing measure and packing dimension \cite{tricot-packing} in 1982, more than sixty years after Hausdorff measures and Hausdorff dimension \cite{Hausdorff-measure-paper} appeared in 1918.

When $s_H>0$, Theorem \ref{t:main} also contributes to the study of doubling measures with prescribed dimensions. Vol'berg and Konyagin \cite{VK-doubling} constructed doubling measures on compact doubling metric spaces, and Luukkainen and Saksman \cite{doubling-exists} extended the construction to complete doubling metric spaces. Wu \cite{Wu-dimension} and K\"aenm\"aki, Rajala, and Suomala \cite{KRS-cubes} showed that doubling measures can be concentrated on sets of arbitrarily small Hausdorff dimension and packing dimension, respectively. Peng and Wen \cite{Peng-Wen} constructed doubling measures on $[0,1]^n$ of every exact dimension in $(0,n]$, and more generally doubling measures whose lower local dimension takes a prescribed range of values. By Theorem \ref{t:main}, for all $0<s_H\leq s_P\leq q$, every complete Ahlfors $q$-regular metric space carries a doubling measure $\mu$ with $\dim_H\mu=s_H$ and $\dim_P\mu=s_P$.

As in \cite{BS4}, the measures that we construct are defined by distributing mass down systems of metric cubes. More specifically, we let $\nu$ be an Ahlfors $q$-regular background measure on $\XX$, let $(\Delta_k)_{k\in\ZZ}$ be an appropriate system of metric cubes on $\XX$, and let $\mathbf{s}=(s_k)_{k=1}^\infty$ be a sequence of ``target entropies'' whose Ces\`aro averages have lower limit $s_H$ and upper limit $s_P$. We define a Radon measure $\mu=\muse$ by specifying the value of $\muse(Q)$ on cubes $Q\in\bigcup_{k=0}^\infty \Delta_k$ so that $\muse(Q)=\nu(Q)$ for all $Q\in\Delta_0$ and the entropy associated to cubes in level $k$ is $s_k$; that these prescribed values are attained by a Radon measure follows from \cite{Badger-MethodI}. For systems of cubes with a fixed ratio of consecutive side lengths, the strong law of large numbers then gives, at $\muse$-a.e.~$x\in\XX$, \begin{align}\liminf_{r\downarrow 0} \frac{\log\muse(B(x,r))}{\log r}&=\liminf_{k\rightarrow\infty} \frac{s_1+\dots+s_k}{k},\\
\limsup_{r\downarrow 0} \frac{\log\muse(B(x,r))}{\log r}&=\limsup_{k\rightarrow\infty} \frac{s_1+\dots+s_k}{k}.\end{align} Thus, in principle, constructing $\muse$ with prescribed Hausdorff and packing dimensions is as simple as choosing an appropriate sequence of target entropies. However, to prove that in any Ahlfors regular metric space, we can specify the entropy associated to cubes in level $k$ to be $s_k$ and verify the required properties in Theorem \ref{t:main}, we need to solve three technical challenges.

\emph{Challenge 1.} To design a measure with $s_H<q$ and $s_P=q$, we cannot use a standard system of $b$-adic metric cubes, because the primary construction in \cite{BS4} demands that the scaling factor $b\rightarrow\infty$ as $s_P\rightarrow q$. (This issue does not occur when $\XX=\RR^n$ and $s_P=n$.) The obstruction is a matter of counting: at a fixed ratio of consecutive side lengths, the estimates on the number of children of a metric cube in an Ahlfors regular space do not guarantee enough children to prescribe entropy arbitrarily close to the ambient dimension. After a short time, one may arrive at the idea of using a system of metric cubes whose ratios of consecutive side lengths tend to infinity \emph{slowly}. This is a natural way to reach the ambient entropy, but it introduces error terms into the computation of local dimensions, and these errors must be controlled in order to use the strong law of large numbers. The approach that we adopt is to let the ratios grow slowly \emph{and} to keep them constant through long stretches of generations, so that the ratio changes only on a \emph{sparse} set of generations. This provides the full range of entropies while keeping the errors under control. As one might imagine, using systems of metric cubes with growing ratios also complicates the proof that the measure is doubling when $s_H>0$. For that part of the construction, we introduce intermediate levels that divide a large change of scale into bounded steps.

\emph{Challenge 2.} When $s_P=m$, we need to design measures which are $m$-rectifiable, but cannot check this property using Corollary \ref{c:bk}. Instead we show how to use Badger and Schul's square packing construction to build Lipschitz maps $f:E\subset[0,1]^m\rightarrow\XX$ whose images have positive measure. The underlying idea is basic enough: try using target entropies $s_k<m$ such that $\limsup_{k\rightarrow\infty} k^{-1}(s_1+\dots+s_k)=m$ slowly in order to check the hypothesis of Theorem \ref{t:pack} for trees of cubes whose leaves have significant measure. But at the time that we wrote \cite{BS4}, we struggled to overcome issues related to rational approximation that arise when building appropriate trees of cubes when $s_P=m$ or $s_P=q$; see the argument starting near \cite[(6.4)]{BS4}. In this paper, we bypass the rational approximation issue by reversing the order of the argument. In \cite{BS4}, we first built trees with a controlled number of cubes in each level, by requiring that long chains of cubes contain a minimum number of central children, and then estimated the measure of the leaves. Here, we first identify the descendants whose measure is significant relative to the entropy, and then bound the number of such cubes from above.

\emph{Challenge 3.} When $s_P=m-1$, the measures that we design are $m$-rectifiable by Corollary \ref{c:bk}. However, we must ensure that the measures are purely $(m-1)$-unrectifiable. To do so, we arrange for our measures to have zero lower $(m-1)$-density, meaning that \begin{equation}\label{zero-m-1-density} \liminf_{r\downarrow 0} \frac{\muse(B(x,r))}{r^{m-1}}=0\quad\text{at $\muse$-a.e. $x\in\XX$,}\end{equation} because this implies that $\muse$ vanishes on every set of $\sigma$-finite $(m-1)$-dimensional packing measure $\Pack^{m-1}$, including Lipschitz images of $E\subset [0,1]^{m-1}$. In Euclidean spaces $\XX=\RR^n$, for a self-similar Bernoulli measure, we could check \eqref{zero-m-1-density} using the law of the iterated logarithm for i.i.d.~random variables. Unfortunately, in a general Ahlfors regular space, the number of children of a metric cube fluctuates cube-by-cube and the natural sequence of random variables associated to $\muse$ is no longer identically distributed. Thus, this simple argument is no longer available. Instead, we choose target entropies whose averages have upper limit $m-1$, but take excursions above that value with enough accumulated excess to outweigh fluctuations of the measure and achieve the required vanishing density.

The rest of the paper is organized as follows. In a preliminary section, \S\ref{s:gmt}, we recall basic notions of dimensions of sets and measures (\S\ref{ss:dim}), discuss how to detect dimension using very general partitions (see Lemma \ref{t:mf}, borrowed from \cite{KRS-multifractal}), record the proof of Corollary \ref{c:bk} (\S\ref{ss:bk}), and recall the statement of the square packing construction of Lipschitz maps from \cite{BS4} (\S\ref{ss:spc}). In \S\ref{s:cubes}, we develop metric cubes with variable scaling, specialize to $b$-adic cubes with sparse jumps, and state a criterion from \cite{Badger-MethodI} for prescribed masses of cubes to be the values of a Radon measure. In \S\ref{s:entropy}, we construct measures with prescribed entropies and compute their local dimensions. Then, in \S\ref{s:doubling}, we introduce intermediate levels between consecutive generations, construct doubling measures with prescribed entropy, and prove the existence of doubling measures with prescribed lower and upper Hausdorff and packing dimensions. Section \ref{s:upper} records a condition on the target entropies that guarantees rectifiability and treats the upper endpoint $s_P=m$, including the critical case $m=q$. Section \ref{s:lower} treats the lower endpoint $s_P=m-1$ by establishing the vanishing density property \eqref{zero-m-1-density}. Finally, \S\ref{s:zero} treats the cases with $s_H=0$, using the construction of \S\ref{s:entropy}.

Table \ref{tab:cases} records where each case of Theorem \ref{t:main} is proved, which system of cubes is used, and which results give $m$-rectifiability and pure $(m-1)$-unrectifiability.

\begin{table}[htb]
\caption{When $s_H>0$, the measures are given by Definition \ref{d:endpoint-measure} and are doubling. When $s_H=0$, they are given by Definition \ref{d:central}.}\label{tab:cases}
\centering
\small
\renewcommand{\arraystretch}{1.25}
\setlength{\tabcolsep}{3.5pt}
\begin{tabular}{llllll}
\hline
 & case & proof & cubes & rectifiability & pure unrectifiability\\
\hline
$s_H>0$ & $m-1<s_P<m$ & \S\ref{s:upper} & sparse jumps & Corollary \ref{c:bk} & Remark \ref{r:pure-unrectifiable}\\
 & $s_P=m<q$ & \S\ref{ss:upper-fixed} & no jumps & Proposition \ref{p:square} & Remark \ref{r:pure-unrectifiable}\\
 & $s_P=m=q$ & \S\ref{ss:upper-q} & sparse jumps & Proposition \ref{p:square} & Remark \ref{r:pure-unrectifiable}\\
 & $s_P=m-1$ & \S\ref{s:lower} & no jumps & Corollary \ref{c:bk} & Lemma \ref{l:lower-fluctuation}\\
\hline
$s_H=0$ & $m-1<s_P<\min\{m,q\}$ & \S\ref{ss:zero-fixed} & no jumps & Corollary \ref{c:bk} & Remark \ref{r:pure-unrectifiable}\\
 & $s_P=m<q$ & \S\ref{ss:zero-fixed} & no jumps & Proposition \ref{p:square} & Remark \ref{r:pure-unrectifiable}\\
 & $s_P=m-1>0$ & \S\ref{ss:zero-fixed} & no jumps & Corollary \ref{c:bk} & Lemma \ref{l:zero-fluctuation}\\
 & $s_P=m-1=0$ & \S\ref{ss:zero-fixed} & no jumps & Corollary \ref{c:bk} & no atoms\\
 & $s_P=q<m$ & \S\ref{ss:zero-q} & sparse jumps & Corollary \ref{c:bk} & Remark \ref{r:pure-unrectifiable}\\
 & $s_P=q=m$ & \S\ref{ss:zero-q} & sparse jumps & Proposition \ref{p:square} & Remark \ref{r:pure-unrectifiable}\\
\hline
\end{tabular}
\end{table}

\subsection*{Acknowledgements} Several of the ideas brought to fruition in this paper stemmed from conversations with Raanan Schul in the development of our joint paper \cite{BS4}.

\subsection*{Development of the results and AI assistance} Initial drafts of this paper, defining metric cubes with variable scaling, were written in January 2024. The author presented the main idea in the talk ``Metric cubes with variable scaling and applications'' in the AMS Special Session on Non-smooth Analysis and Geometry at the 2025 Joint Mathematics Meetings in Seattle (January 2025), and in the talk ``Doubling measures that charge rectifiable surfaces but not rectifiable curves'' at the conference Quasiweekend III in Helsinki (June 2025). These talks covered the scaling factors that tend to infinity, the new error term that arises when the strong law of large numbers is used to compute local dimensions, and the requirement that the scaling factors change only on a sparse set of generations. At that time, the author had a working construction of rectifiable measures with prescribed Hausdorff and packing dimensions with $s_P=m$ and $s_P=q$. The outstanding issues were the computations needed to make the measures doubling when $s_H>0$ and the estimates needed to establish the vanishing lower density \eqref{zero-m-1-density} at the lower endpoint.

These final computations were carried out in September 2026 with the assistance of Claude Opus 5.5 (Anthropic). In particular, exchanges with the LLM led to the construction in \S\ref{s:doubling}, in which each generation is divided into intermediate levels that carry Type I weights, followed by a single Type II transition, and to Lemma \ref{l:lower-fluctuation} establishing vanishing lower density. The manuscript was completed and revised in October 2026 and November 2026 using a combination of Claude Opus 5.5, Claude Fable 5.2, and GPT-6 Astra Pro (OpenAI). These models were used to check proofs, to propose revisions, and to draft and edit the exposition and figures. The author directed this work and verified all of the mathematics, and takes full responsibility for the contents of the paper.

\section{Packing dimension and rectifiability}\label{s:gmt}

\subsection{Dimensions of sets and measures}\label{ss:dim} In a metric space $\XX$, we denote the distance between points by $|x-y|$ and the diameter of a set by $\diam E:=\sup\{|x-y|:x,y\in E\}$. We let $U(x,r):=\{y\in\XX:|x-y|<r\}$ and $B(x,r):=\{y\in\XX:|x-y|\leq r\}$ denote the \emph{open ball} and \emph{closed ball} with \emph{center} $x\in\XX$ and \emph{radius} $r>0$, respectively. If $B$ is a ball and $\lambda>0$, we let $\lambda B$ denote the ball with the same center and topological type as $B$, but with radius multiplied by $\lambda$.

For $s\geq 0$, we let $\Haus^s$ denote the $s$-dimensional Hausdorff measure, in which a set $E_i$ in a cover has cost $(\diam E_i)^s$, and we let
\[
 \Pack^s(E):=\inf\left\{\sum_{i=1}^\infty P^s(E_i):E\subset\bigcup_{i=1}^\infty E_i\right\}\quad\text{for all }E\subset\XX
\]
denote the $s$-dimensional packing measure, where $P^s(\emptyset):=0$ and
\[
 P^s(E):=\lim_{\delta\downarrow0}\sup\left\{\sum_{i=1}^\infty(2r_i)^s:B(x_i,r_i)\text{ pack in }E,\ 0<r_i\leq\delta/2\right\},
\] where we say that balls $B(x_i,r_i)$ \emph{pack in} $E$ if $x_i\in E$ and the balls are pairwise disjoint.
Since a general metric space need not be uniformly perfect, we use this radial definition, following Cutler \cite{Cutler-density-theorem}. The \emph{Hausdorff dimension} $\dim_H E$ and \emph{packing dimension} $\dim_P E$ are the infima of the exponents $s\geq 0$ with $\Haus^s(E)=0$ and $\Pack^s(E)=0$, respectively. Letting $N(E,r)$ denote the least number of closed balls of radius $r$ needed to cover $E$, the \emph{upper Minkowski dimension} of $E$ is
\[
 \overline{\dim}_M\,E:=\limsup_{r\downarrow0}\frac{\log N(E,r)}{\log(1/r)}.
\]
Let $f:E\to\XX$ be an $L$-Lipschitz map with $L>0$ defined on a subset $E$ of a metric space. If $x_i,x_j\in E$ and the closed balls $B(f(x_i),r_i)$ and $B(f(x_j),r_j)$ are disjoint, then $|f(x_i)-f(x_j)|>\max\{r_i,r_j\}$. It follows that $|x_i-x_j|>(r_i+r_j)/(2L)$ and the balls $B(x_i,r_i/(2L))$ and $B(x_j,r_j/(2L))$ are disjoint. Therefore,
\begin{equation}\label{lip-packing}
 \Pack^s(f(E))\leq(2L)^s\Pack^s(E)\quad\text{for all }s\geq0.
\end{equation}
Since constant maps have one-point images, Lipschitz maps do not increase packing dimension.

Let $\mu$ be a Borel measure on a separable metric space $\XX$. The \emph{lower} and \emph{upper Hausdorff dimensions} of $\mu$ are
\[
 \underline{\dim}_H\,\mu:=\inf\{\dim_H E:\mu(E)>0\}\quad\text{and}\quad\overline{\dim}_H\,\mu:=\inf\{\dim_H E:\mu(\XX\setminus E)=0\},
\]
where $E$ ranges over Borel sets. The \emph{lower} and \emph{upper packing dimensions} $\underline{\dim}_P\,\mu$ and $\overline{\dim}_P\,\mu$ are defined in the same way with $\dim_P$ in place of $\dim_H$. When the lower and upper Hausdorff dimensions of $\mu$ agree, we write $\dim_H\mu$ for the common value and say that $\mu$ has \emph{exact} Hausdorff dimension. Exact packing dimension $\dim_P\mu$ is defined similarly. The \emph{lower} and \emph{upper local dimensions} of $\mu$ at $x\in\XX$ are
\[
 \liminf_{r\downarrow0}\frac{\log\mu(B(x,r))}{\log r}\quad\text{and}\quad\limsup_{r\downarrow0}\frac{\log\mu(B(x,r))}{\log r}.
\]
By a theorem of Tamashiro \cite[Theorem 1.8]{tamashiro}, if $\mu$ is finite, then $\underline{\dim}_H\,\mu$ and $\overline{\dim}_H\,\mu$ are the $\mu$-essential infimum and supremum of the lower local dimension, and $\underline{\dim}_P\,\mu$ and $\overline{\dim}_P\,\mu$ are the $\mu$-essential infimum and supremum of the upper local dimension. We apply this to a locally finite Borel measure $\mu$ by restricting $\mu$ to the sets in a countable cover of $\XX$ by open sets of finite measure, which exists because $\XX$ is separable. If $U$ is open and $x\in U$, then $\mu(B(x,r))=\mu\res U(B(x,r))$ for all sufficiently small $r>0$.

\begin{remark}[packing-density null criterion]\label{r:packing-density-null}
Let $s>0$. By \cite[Corollary~2.4]{BLZ}, the part of a locally finite Borel measure $\mu$ on which
\[
 \liminf_{r\downarrow0}\frac{\mu(B(x,r))}{r^s}=0
\]
is singular to sets of finite $\Pack^s$ measure. In particular, if this lower $s$-density vanishes at $\mu$-a.e.~$x$, then $\mu(E)=0$ for every set $E$ of $\sigma$-finite $\Pack^s$ measure.
\end{remark}

\begin{remark}[pure unrectifiability]\label{r:pure-unrectifiable}
Let $m\geq1$ be an integer and let $\mu$ be a locally finite Borel measure on a metric space $\XX$. Then $\mu$ is purely $(m-1)$-unrectifiable in the sense of \eqref{unrectifiable} whenever
\[
 \liminf_{r\downarrow0}\frac{\mu(B(x,r))}{r^{m-1}}=0\quad\text{at $\mu$-a.e.~}x\in\XX.
\]
When $m\geq2$, this follows from Remark \ref{r:packing-density-null}, because every Lipschitz image of a subset of $\RR^{m-1}$ has $\sigma$-finite $\Pack^{m-1}$ measure; see e.g.~\cite[Lemma 2.8]{BS1}. When $m=1$, the condition says that $\mu$ has no atoms. The condition holds, in particular, when the upper local dimension of $\mu$ is greater than $m-1$ at $\mu$-a.e.~point.
\end{remark}

Recall that a metric space $\XX$ is \emph{doubling} if $\sup_{x\in\XX}\sup_{r>0} N(B(x,r),r/2)<\infty$, i.e.~every closed ball can be covered by a uniform number of balls of half the radius. On any doubling space, one may compute local dimensions using sets drawn from a sequence of shrinking partitions of $\XX$ in place of balls (e.g.~half-open dyadic cubes in $\RR^n$).

\begin{definition}[{\cite[\S3]{KRS-multifractal}}] \label{KRS-partition} Let $\XX$ be a metric space, let $\delta>0$, and let $\Lambda\geq 1$. We say that $\mathcal{P}$ is a \emph{$(\delta,\Lambda)$-partition} if $\mathcal{P}$ is a collection of Borel sets that partition $\XX$ (i.e.~$\XX=\bigcup_{P\in\mathcal{P}}P$ and $P\cap P'=\emptyset$ for all distinct $P,P'\in\mathcal{P}$) and there exist pairwise disjoint closed balls $\{B_P:P\in\mathcal{P}\}$ of radius $\delta$ such that $P\subset \Lambda B_P$ for all $P\in\mathcal{P}$.
\end{definition}

\begin{definition}[{\cite[\S3]{KRS-multifractal}}] \label{KRS-partition-2} Let $\XX$ be a metric space and let $(\delta_k)_{k\geq k_0}$ be a sequence of positive numbers such that
\begin{equation}\label{geometric-decay}
 \delta_{k+1}\leq c\delta_k\quad\text{for all $k\geq k_0$},\quad\text{for some $0<c<1$, and}
\end{equation}
\begin{equation}\label{asymptotic-log-stability}
 \lim_{k\rightarrow\infty}\frac{\log\delta_k}{\log\delta_{k+1}}=1.
\end{equation}
We say that $(\mathcal{P}_k)_{k\geq k_0}$ is an \emph{admissible family of $\delta_k$-partitions} if there is $\Lambda\geq1$ such that $\mathcal{P}_k$ is a $(\delta_k,\Lambda)$-partition for all $k\geq k_0$.
\end{definition}

\begin{lemma}[{\cite[Proposition 3.1]{KRS-multifractal}}] \label{t:mf} Let $\XX$ be a doubling metric space and let $(\mathcal{P}_k)_{k\geq k_0}$ be an admissible family of $\delta_k$-partitions. Let $P_k(x)$ denote the unique set in $\mathcal{P}_k$ containing $x\in\XX$. If $\mu$ is a locally finite Borel measure on $\XX$, then at $\mu$-a.e.~$x\in\XX$,
\begin{align}
 \label{dim5} \liminf_{r\downarrow 0}\frac{\log\mu(B(x,r))}{\log r}
 &=\liminf_{k\rightarrow\infty}\frac{\log\mu(P_k(x))}{\log\delta_k},\\
 \label{dim6} \limsup_{r\downarrow 0}\frac{\log\mu(B(x,r))}{\log r}
 &=\limsup_{k\rightarrow\infty}\frac{\log\mu(P_k(x))}{\log\delta_k}.
\end{align}
\end{lemma}

\subsection{Proof of Corollary \ref{c:bk}}\label{ss:bk}

Let $\XX$ be a complete, separable metric space, let $Z\subset\XX$ be a countable set whose closure $\overline{Z}=\XX$, let $\mu$ be a locally finite Borel measure on $\XX$, and let $m\geq 1$ be an integer. The set of points where the upper local dimension of $\mu$ is less than $m$ satisfies $ A=\bigcup_{z\in Z}\bigcup_{n=3}^\infty A_{z,n}$,
where $$A_{z,n}=\left\{x\in B\Big(z,\frac{1}{n}\Big):r^{m-1/n}\leq\mu(B(x,r))\leq n\quad\text{for all }0<r<\frac{3}{n}\right\},
$$ where $z\in Z$ and $n\geq 3$ is an integer. To proceed, suppose that $F:=A_{z,n}$ is nonempty for some $z$ and $n$. It suffices to prove that $F$ is the image of a Lipschitz map defined on a compact subset of $[0,1]^m$. Steps 1 and 2 below adapt the proof of \cite[Theorem 3.1]{BS4}. Step 3 replaces the square packing construction used there with Theorem \ref{t:bk}.

\emph{1. $F$ is closed.} Suppose that $x_i\in F$ for all $i\geq 1$ and $\lim_{i\rightarrow\infty} x_i=x$. Then $x\in B(z,1/n)$, because $x_i$ belongs to $B(z,1/n)$ for all $i\geq 1$ and the ball is closed. Further, since $B(x_i,r-|x_i-x|)\subset B(x,r)\subset B(x_i,r+|x_i-x|)$ for all $i\geq 1$, it follows that $$r^{m-1/n}\leq \limsup_{i\rightarrow\infty} \mu(B(x_i,r-|x_i-x|))\leq \mu(B(x,r))\leq \liminf_{i\rightarrow\infty} \mu(B(x_i,r+|x_i-x|))\leq n$$ for all $0<r<3/n$. Hence $x\in F$ and this shows that $F$ is closed.

\emph{2. $F$ has upper Minkowski dimension at most $m-1/n$.} Let $0<r<1/n$. If $x_1,\dots,x_j\in F$ and $B(x_1,r),\dots, B(x_j,r)$ are pairwise disjoint, then $$j r^{m-1/n} \leq \sum_{i=1}^j \mu(B(x_i,r))=\mu\bigg(\bigcup_{i=1}^j B(x_i,r)\bigg)\leq \mu\left(B\Big(x_1,r+\frac{2}{n}\Big)\right)\leq n.$$ Hence $j$ is finite. Suppose that $x_1,\dots,x_j$ are maximal in the sense that for every $x\in F$, the ball $B(x,r)$ intersects $B(x_i,r)$ for some $1\leq i\leq j$. Then $F\subset B(x_1,2r)\cup\dots\cup B(x_j,2r)$ and $j\leq n/r^{m-1/n}$.  To summarize, for any $0<r<2/n$, the set $F$ can be covered by $n/(r/2)^{m-1/n}=C(n,m)r^{-(m-1/n)}$ or fewer balls of radius $r$. It follows immediately that $F$ is totally bounded and the upper Minkowski dimension of $F$ is at most $m-1/n$.

\emph{3. $F$ is a Lipschitz image of a compact subset of $[0,1]^m$.} Since $\XX$ is complete and $F$ is closed and totally bounded, $F$ is compact. By Theorem \ref{t:bk}, since $F$ is compact and $\overline{\dim}_M\, F\leq m-1/n<m=\dim_H [0,1]^m$, there exists a compact set $E\subset[0,1]^m$ and a Lipschitz surjection $f:E\rightarrow F$. This completes the proof of the corollary.

\subsection{Square packing construction of Lipschitz maps} \label{ss:spc}

\begin{definition}[{\cite[Definition 2.1]{BS4}}]\label{def:tree} Let $\XX$ be a metric space. We call $\mathcal{T}=\bigsqcup_{j=0}^\infty\mathcal{T}_j$ a \emph{tree of sets} in $\XX$ if \begin{itemize}
\item $\mathcal{T}_j$ is a finite multiset\footnote{A finite multiset is a finite unordered list with repetition allowed. For example, $\{2,0,2,6\}$ is a multiset with 4 elements. A function between multisets is a function between sets formed by assigning any repeated elements a different color. For example, we could define a surjective function $f:\{2,0,2,6\}\rightarrow\{0,0,0\}$ by forming sets $\{2,0,\color{PTblue}2\color{black},6\}$ and $\{0,\color{PTblue}0\color{black},\color{PTmagenta}0\color{black}\}$ and defining $f(2)=0$, $f(0)=\color{PTblue}0\color{black}$, and $f(\color{PTblue}2\color{black})=f(6)=\color{PTmagenta}0\color{black}$.} of nonempty subsets of $\XX$ for all $j\geq 0$;
\item $\#\mathcal{T}_0=1$ (counting multiplicity); and,
\item there is a function $\uparrow:\bigsqcup_{j=1}^\infty \mathcal{T}_j\rightarrow \mathcal{T}$ such that $Q\subset Q^\uparrow\in\mathcal{T}_{j-1}$ for all $j\geq 1$ and $Q\in\mathcal{T}_j$.
\end{itemize} We call $Q^\uparrow$ the \emph{parent} of $Q$ and call $Q$ a \emph{child} of $Q^\uparrow$. For all $j\geq 0$ and $Q\in\mathcal{T}_j$, we let $\Child(Q):=\{R\in\mathcal{T}_{j+1}:Q=R^\uparrow\}$ denote the \emph{multiset of children} of $Q$. We denote the unique set in $\mathcal{T}_0$ by $\Top(\mathcal{T})$.
Finally, the \emph{set of leaves} of $\mathcal{T}$ is defined by \begin{equation}\leaves(\mathcal{T}):=\bigcap_{j=0}^\infty \bigcup_{Q\in\mathcal{T}_j} Q.\end{equation}
\end{definition}

\begin{theorem}[square packing construction {\cite[Theorem 2.5]{BS4}}] \label{t:pack} Let $m\geq 1$ be an integer. Let $\mathcal{T}=\bigsqcup_{j=0}^\infty\mathcal{T}_j$ be a tree of sets in a metric space $\XX$. For each $j\geq 0$, assign \begin{equation}\label{N-def} N_j:=\max_{Q\in\mathcal{T}_j}\#\Child(Q)\quad \text{and}\quad D_j:=\max_{Q\in\mathcal{T}_j}\,\diam Q.\end{equation} If the closure of $\Top(\mathcal{T})$ is compact and \begin{equation}\label{root-condition} S:=\sum_{j=0}^\infty \left(\prod_{i=0}^j \lceil N_{i}^{1/m}\rceil \right) D_j<\infty,\end{equation} then there exists a compact set $E\subset [0,1]^m$ and an $S$-Lipschitz map $f:E\rightarrow\XX$ such that $f(E)$ contains $\leaves(\mathcal{T})$.\end{theorem}

What is the connection between Theorem \ref{t:pack} and square packings? The proof is based on the following elementary observation: any collection of $N$ squares in the plane of side length $s$ can be packed inside of a square of side length $\lceil N^{1/2}\rceil s$, where the \emph{ceiling} of a real number $x$ is $\lceil x\rceil:=\inf\{n\in\ZZ:x\leq n\}$. More generally, any collection of $N$ cubes in $\RR^m$ of side length $s$ can be packed inside of a cube of side length $\lceil N^{1/m}\rceil s$.

\section{Metric cubes with variable scaling}\label{s:cubes}

\subsection{KRS cubes with variable scaling}\label{ss:cubes}

A metric space $\XX$ has \emph{totally bounded balls} if every closed ball in $\XX$ can be covered by finitely many closed balls of any given radius. On such spaces, K\"aenm\"aki, Rajala, and Suomala \cite{KRS-cubes} gave a short construction of nested partitions into sets that behave like $b$-adic cubes. Our construction at the endpoint $s_P=q$ uses ratios of consecutive side lengths tending to infinity. The following version of their construction only requires the side lengths to decay at least geometrically.

\begin{theorem}[KRS cubes with variable scaling]\label{vKRS-cubes}
Let $\XX$ be a metric space with totally bounded balls. Choose an \emph{origin} $o\in\XX$, a number $b>3$, and a sequence $(r_k)_{k\in\ZZ}$ of positive numbers, called \emph{side lengths}, such that $r_k=b^{-k}$ for all $k\leq 0$ and
\begin{equation}\label{geometric-sides}
 r_k\leq r_{k-1}/b\quad\text{for all }k\geq 1.
\end{equation}
Put
\[
 c_b:=\frac12-\frac1{b-1}\quad\text{and}\quad C_b:=\frac b{b-1}.
\]
For every $k\in\ZZ$, there is a partition $\Delta_k$ of $\XX$ into nonempty bounded Borel sets, called \emph{cubes}, together with points $x_Q$ for $Q\in\Delta_k$, called \emph{centers}, with the following properties.
\begin{enumerate}
\item nesting: if $Q\in\Delta_k$, $R\in\Delta_l$, and $l\geq k$, then $Q\cap R=\emptyset$ or $R\subset Q$;
\item roundness: $U(x_Q,c_b\,r_k)\subset Q\subset B(x_Q,C_b\,r_k)$ for all $Q\in\Delta_k$;
\item origin: for every $k\in\ZZ$, there is $Q\in\Delta_k$ with $x_Q=o$;
\item inheritance: for every $Q\in\Delta_k$, there is $R\in\Delta_{k+1}$ with $R\subset Q$ and $x_R=x_Q$.
\end{enumerate}
\end{theorem}

\begin{proof}
When $r_k=b^{-k}$ for all $k$, this is \cite[Theorem 2.1]{KRS-cubes} with $r=1/b$, for which the constants $\frac12-\frac r{1-r}$ and $\frac1{1-r}$ in that theorem equal $c_b$ and $C_b$. In general, we run the construction of \cite{KRS-cubes} with $r_k$ in place of $r^k$. For $k\leq 0$, the level-$k$ centers are chosen exactly as there. For $k\geq 1$, the level-$k$ centers form a maximal $r_k$-separated set that contains the level-$(k-1)$ centers, which is possible because $r_k\leq r_{k-1}$. The partial order on the indices and the definition of the cubes are unchanged. In the proof of \cite[Theorem 2.1]{KRS-cubes}, the scales enter only through the separation of the level-$k$ centers, the fact that every level-$(k+1)$ center lies within distance $r_k$ of the level-$k$ center to which it is assigned, and the tail sum $\sum_{m>k}r^m$. With the tail bound
\[
 \sum_{i=1}^\infty r_{k+i}\leq r_k\sum_{i=1}^\infty b^{-i}=\frac{r_k}{b-1},
\]
the argument there gives $Q\subset B(x_Q,(1+\frac1{b-1})r_k)$ and $U(x_Q,(\frac12-\frac1{b-1})r_k)\subset Q$ for all $Q\in\Delta_k$. The remaining properties follow as in \cite{KRS-cubes}.
\end{proof}

\begin{remark}[notation]\label{r:KRS}
Let $b$ and $(r_k)_{k\in\ZZ}$ be as in Theorem \ref{vKRS-cubes}. We call a family $(\Delta_k)_{k\in\ZZ}$ of partitions of $\XX$ into nonempty bounded Borel sets, together with centers $x_Q$, a \emph{system of cubes} with side lengths $(r_k)_{k\in\ZZ}$ if properties (1)--(4) of Theorem \ref{vKRS-cubes} hold. Thus, Theorem \ref{vKRS-cubes} provides a system of cubes with any such side lengths. When $r_k=b^{-k}$ for all $k$, we call it a system of \emph{$b$-adic cubes}. We call $\Delta_k$ the $k$-th \emph{level} or \emph{generation} of $\Delta:=\bigsqcup_{k\in\ZZ}\Delta_k$. The same set may belong to $\Delta_k$ for several $k$, and a cube is always considered together with its level. For $Q\in\Delta_k$, we call
\[
 \side Q:=r_k,\qquad U_Q:=U(x_Q,c_b\side Q),\qquad B_Q:=B(x_Q,C_b\side Q)
\]
the \emph{side length}, \emph{inner ball}, and \emph{outer ball} of $Q$. Thus, $\diam Q\leq 2C_b\side Q$. For $j\geq 0$, let $\Child^j(Q)$ be the set of cubes in $\Delta_{k+j}$ contained in $Q$, and put $\Child(Q):=\Child^1(Q)$. We call the cubes in $\Child(Q)$ the \emph{children} of $Q$. Since $x_Q\in U_Q$, exactly one child of $Q$ contains $x_Q$. We call it the \emph{central child} of $Q$ and denote it by $Q^\downarrow$. By inheritance, the center of $Q^\downarrow$ is $x_Q$. For every $k\in\ZZ$, put
\begin{equation}\label{xi-def}
 \xi_k:=\frac{r_{k-1}}{r_k}.
\end{equation}
Thus, $\xi_k\geq b$ is the ratio of the side lengths of cubes in the generations $k-1$ and $k$, and $\xi_k=b$ for all $k\leq 0$. The side lengths are determined by the ratios, since $r_k=(\xi_1\cdots\xi_k)^{-1}$ for all $k\geq 1$. For a system of $b$-adic cubes, $\xi_k=b$ for all $k$.
\end{remark}

\begin{remark}[exhaustion]\label{r:exhaustive}
For $k\in\ZZ$, let $Q^o_k$ be the cube in $\Delta_k$ with center $o$, the origin of system $\Delta$. By inheritance, $Q^o_k\subset Q^o_{k-1}$ for all $k$, and $Q^o_k\supset U(o,c_bb^{-k})$ for all $k\leq 0$. Hence $\XX=\bigcup_{k\leq 0}Q^o_k$. It follows that every cube is contained in $Q^o_k$ for some $k$, and any two cubes have a common ancestor.
\end{remark}

For nonempty sets $A,B\subset\XX$, the \emph{gap} between $A$ and $B$ is
\begin{equation}\label{gap-def}
 \gap(A,B):=\inf\{|x-y|:x\in A,\ y\in B\}.
\end{equation}
If $A\subset A'$ and $B\subset B'$, then $\gap(A,B)\geq\gap(A',B')$.

\begin{lemma}[counting cubes]\label{l:vAR-count}
Let $\XX$ be an Ahlfors $q$-regular metric space, let $\nu$ be a Borel measure on $\XX$, and let $C_\nu>0$ and $D_\nu\geq 1$ be constants such that
\begin{equation}\label{AR-background}
 C_\nu r^q\leq \nu(B(x,r))\leq C_\nu D_\nu r^q\quad\text{for all }x\in\XX\text{ and }0<r\leq\diam\XX.
\end{equation}
Let $(\Delta_k)_{k\in\ZZ}$ be a system of cubes with $b\geq 47$. If $j<k$, $P\in\Delta_j$, $P\neq\XX$, and $L:=r_j/r_k$, then
\begin{equation}\label{vAR-one-generation}
 D_\nu^{-1}3^{-q}L^q
 \leq \#\Child^{k-j}(P)\leq D_\nu3^qL^q,\text{ and}
\end{equation}
\begin{equation}\label{vAR-nu-ratio}
 D_\nu^{-1}3^{-q}L^{-q}\leq\frac{\nu(R)}{\nu(P)}\leq D_\nu3^qL^{-q}\quad\text{for all }R\in\Child^{k-j}(P).
\end{equation}
Moreover, put $C_\XX:=D_\nu6^q$. For all $k\in\ZZ$ and $P\in\Delta_{k-1}$ with $P\neq\XX$, the sets
\begin{align}
 \Inner(P)&:=\{R\in\Child(P):R\cap U(x_P,\tfrac12c_b\side P)\neq\emptyset\},\label{inner-def}\\
 \Outer(P)&:=\Child(P)\setminus\Inner(P)\label{outer-def}
\end{align}
satisfy
\begin{equation}\label{vin-out-count}
 C_\XX^{-1}\xi_k^q\leq\#\Inner(P)\quad\text{and}\quad\#\Child(P)\leq C_\XX\xi_k^q,
\end{equation}
and every $R\in\Child(P)$ with $\gap(R,\XX\setminus P)\leq 9\side R$ belongs to $\Outer(P)$. In particular, $P^\downarrow\in\Inner(P)$ and $\gap(P^\downarrow,\XX\setminus P)>9\side P^\downarrow$.
\end{lemma}

\begin{proof}
Since $b\geq 47$, we have $c_b\geq 22/46\geq 2/5$ and $C_b\leq 47/46\leq 6/5$. Since $U_P\subset P\neq\XX$, we have $c_b\side P\leq\diam\XX$. Hence the lower bound in \eqref{AR-background} applies to balls of radius at most $c_b\side P$. Write $\sigma:=\side P$ and $\tau:=r_k$. The inner balls $U_R$ with $R\in\Child^{k-j}(P)$ are pairwise disjoint and lie in $B_P$. Hence
\[
 \#\Child^{k-j}(P)\,C_\nu(2\tau/5)^q\leq C_\nu D_\nu(6\sigma/5)^q.
\]
Conversely, $U_P\subset P\subset\bigcup_{R\in\Child^{k-j}(P)}B_R$, and therefore
\[
 C_\nu(2\sigma/5)^q\leq \#\Child^{k-j}(P)\,C_\nu D_\nu(6\tau/5)^q.
\]
These are \eqref{vAR-one-generation}. Similarly, $C_\nu(2\tau/5)^q\leq\nu(R)\leq C_\nu D_\nu(6\tau/5)^q$ and $C_\nu(2\sigma/5)^q\leq\nu(P)\leq C_\nu D_\nu(6\sigma/5)^q$, which give \eqref{vAR-nu-ratio}.

Now let $P\in\Delta_{k-1}$ with $P\neq\XX$, and put $\sigma:=\side P$. By \eqref{vAR-one-generation}, $\#\Child(P)\leq D_\nu3^q\xi_k^q\leq C_\XX\xi_k^q$. Since $\tfrac12c_b\geq 11/46>1/5$, the ball $B(x_P,\sigma/5)$ lies in $U(x_P,\tfrac12c_b\sigma)\subset P$, and it is covered by the cubes in $\Inner(P)$. Every child $R$ of $P$ satisfies $\nu(R)\leq\nu(B_R)\leq C_\nu D_\nu(6\sigma/(5\xi_k))^q$, whereas $\nu(B(x_P,\sigma/5))\geq C_\nu(\sigma/5)^q$. Hence $\#\Inner(P)\geq D_\nu^{-1}6^{-q}\xi_k^q$. Finally, let $R\in\Inner(P)$. Then $R\subset B(x_P,\tfrac12c_b\sigma+\diam R)$ and $\XX\setminus P\subset\XX\setminus U(x_P,c_b\sigma)$. Since $\diam R\leq 2C_b\side R$, $\sigma=\xi_k\side R$, and $\xi_k\geq 47$,
\[
 \gap(R,\XX\setminus P)\geq\left(\tfrac12c_b\xi_k-2C_b\right)\side R\geq (423/46)\side R>9\side R.
\]
Thus, every $R\in\Child(P)$ with $\gap(R,\XX\setminus P)\leq 9\side R$ belongs to $\Outer(P)$. Finally, $P^\downarrow$ contains $x_P$, and hence $P^\downarrow\in\Inner(P)$.
\end{proof}

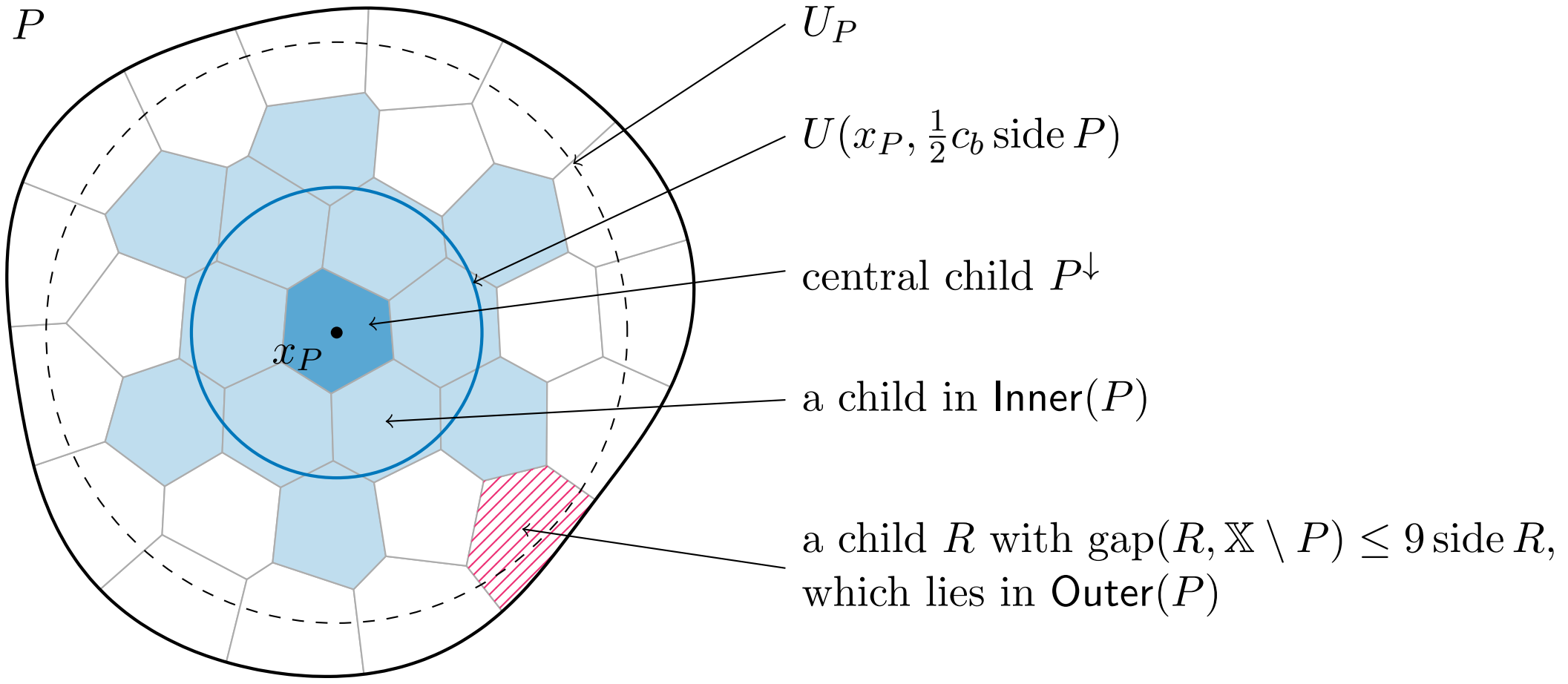
\begin{figure}[htb]
\centering
\begin{tikzpicture}[font=\small]
\path[fill=white,draw=gray!65,line width=0.4pt] (-1.556,-1.802) -- (-0.801,-2.200) -- (-0.984,-2.938) -- (-1.009,-2.931) -- (-1.061,-2.915) -- (-1.113,-2.898) -- (-1.164,-2.881) -- (-1.215,-2.862) -- (-1.266,-2.843) -- (-1.316,-2.822) -- (-1.366,-2.801) -- (-1.415,-2.778) -- (-1.465,-2.754) -- (-1.513,-2.730) -- (-1.561,-2.704) -- (-1.609,-2.678) -- (-1.656,-2.650) -- (-1.702,-2.621) -- (-1.748,-2.592) -- (-1.793,-2.561) -- (-1.837,-2.529) -- (-1.881,-2.496) -- (-1.923,-2.462) -- (-1.965,-2.427) -- (-2.006,-2.391) -- (-2.046,-2.354) -- (-2.085,-2.315) -- (-2.110,-2.290) -- cycle;
\path[fill=white,draw=gray!65,line width=0.4pt] (-0.801,-2.200) -- (-0.571,-2.050) -- (0.151,-2.289) -- (0.243,-3.046) -- (0.213,-3.049) -- (0.160,-3.055) -- (0.107,-3.059) -- (0.053,-3.062) -- (-0.000,-3.065) -- (-0.054,-3.066) -- (-0.107,-3.066) -- (-0.161,-3.065) -- (-0.214,-3.064) -- (-0.268,-3.061) -- (-0.321,-3.057) -- (-0.375,-3.053) -- (-0.428,-3.048) -- (-0.482,-3.041) -- (-0.535,-3.034) -- (-0.588,-3.026) -- (-0.641,-3.017) -- (-0.694,-3.008) -- (-0.747,-2.997) -- (-0.800,-2.985) -- (-0.853,-2.973) -- (-0.905,-2.960) -- (-0.957,-2.946) -- (-0.984,-2.938) -- cycle;
\path[fill=white,draw=gray!65,line width=0.4pt] (0.151,-2.289) -- (0.438,-1.994) -- (1.157,-2.080) -- (1.476,-2.485) -- (1.449,-2.510) -- (1.410,-2.543) -- (1.369,-2.576) -- (1.329,-2.608) -- (1.287,-2.639) -- (1.244,-2.669) -- (1.201,-2.698) -- (1.157,-2.726) -- (1.112,-2.754) -- (1.067,-2.780) -- (1.021,-2.805) -- (0.974,-2.829) -- (0.927,-2.852) -- (0.879,-2.874) -- (0.830,-2.894) -- (0.781,-2.914) -- (0.731,-2.932) -- (0.681,-2.949) -- (0.630,-2.965) -- (0.579,-2.980) -- (0.528,-2.993) -- (0.476,-3.005) -- (0.424,-3.017) -- (0.372,-3.026) -- (0.319,-3.035) -- (0.266,-3.043) -- (0.243,-3.046) -- cycle;
\path[fill=white,draw=gray!65,line width=0.4pt] (-2.064,-0.972) -- (-1.533,-1.370) -- (-1.556,-1.802) -- (-2.110,-2.290) -- (-2.123,-2.276) -- (-2.160,-2.236) -- (-2.196,-2.196) -- (-2.230,-2.154) -- (-2.264,-2.111) -- (-2.297,-2.068) -- (-2.328,-2.024) -- (-2.358,-1.979) -- (-2.388,-1.933) -- (-2.416,-1.887) -- (-2.442,-1.840) -- (-2.468,-1.793) -- (-2.493,-1.745) -- (-2.516,-1.697) -- (-2.539,-1.649) -- (-2.560,-1.600) -- (-2.580,-1.550) -- (-2.600,-1.501) -- (-2.618,-1.451) -- (-2.636,-1.401) -- (-2.652,-1.351) -- (-2.668,-1.301) -- (-2.683,-1.251) -- (-2.697,-1.201) -- (-2.700,-1.188) -- cycle;
\path[fill=white,draw=gray!65,line width=0.4pt] (-1.533,-1.370) -- (-1.020,-1.068) -- (-0.500,-1.380) -- (-0.571,-2.050) -- (-0.801,-2.200) -- (-1.556,-1.802) -- cycle;
\path[fill=PTblue!25,draw=gray!65,line width=0.4pt] (-0.500,-1.380) -- (-0.031,-1.142) -- (0.333,-1.336) -- (0.438,-1.994) -- (0.151,-2.289) -- (-0.571,-2.050) -- cycle;
\path[fill=white,draw=gray!65,line width=0.4pt] (0.929,-1.041) -- (1.311,-1.321) -- (1.157,-2.080) -- (0.438,-1.994) -- (0.333,-1.336) -- cycle;
\path[fill=white,draw=gray!65,line width=0.4pt] (1.157,-2.080) -- (1.311,-1.321) -- (1.871,-1.184) -- (2.307,-1.491) -- (2.303,-1.496) -- (2.274,-1.534) -- (2.246,-1.572) -- (2.217,-1.611) -- (2.188,-1.649) -- (2.160,-1.687) -- (2.131,-1.726) -- (2.102,-1.764) -- (2.073,-1.802) -- (2.044,-1.841) -- (2.015,-1.879) -- (1.986,-1.917) -- (1.956,-1.956) -- (1.926,-1.994) -- (1.895,-2.033) -- (1.865,-2.071) -- (1.833,-2.109) -- (1.802,-2.147) -- (1.769,-2.185) -- (1.736,-2.222) -- (1.703,-2.260) -- (1.669,-2.297) -- (1.634,-2.334) -- (1.598,-2.370) -- (1.562,-2.406) -- (1.525,-2.441) -- (1.487,-2.476) -- (1.476,-2.485) -- cycle;
\path[fill=white,draw=gray!65,line width=0.4pt] (-2.412,0.084) -- (-1.903,-0.399) -- (-2.064,-0.972) -- (-2.700,-1.188) -- (-2.710,-1.150) -- (-2.723,-1.100) -- (-2.735,-1.050) -- (-2.747,-1.000) -- (-2.758,-0.950) -- (-2.768,-0.900) -- (-2.779,-0.850) -- (-2.788,-0.800) -- (-2.798,-0.750) -- (-2.807,-0.700) -- (-2.815,-0.650) -- (-2.824,-0.600) -- (-2.832,-0.551) -- (-2.840,-0.501) -- (-2.848,-0.451) -- (-2.855,-0.401) -- (-2.863,-0.352) -- (-2.870,-0.302) -- (-2.877,-0.252) -- (-2.884,-0.202) -- (-2.890,-0.151) -- (-2.897,-0.101) -- (-2.903,-0.051) -- (-2.908,0.000) -- (-2.914,0.051) -- (-2.914,0.052) -- cycle;
\path[fill=PTblue!25,draw=gray!65,line width=0.4pt] (-1.903,-0.399) -- (-1.403,-0.230) -- (-0.998,-0.497) -- (-1.020,-1.068) -- (-1.533,-1.370) -- (-2.064,-0.972) -- cycle;
\path[fill=PTblue!25,draw=gray!65,line width=0.4pt] (-0.492,-0.263) -- (-0.049,-0.543) -- (-0.031,-1.142) -- (-0.500,-1.380) -- (-1.020,-1.068) -- (-0.998,-0.497) -- cycle;
\path[fill=PTblue!25,draw=gray!65,line width=0.4pt] (0.508,-0.236) -- (0.923,-0.492) -- (0.929,-1.041) -- (0.333,-1.336) -- (-0.031,-1.142) -- (-0.049,-0.543) -- cycle;
\path[fill=PTblue!25,draw=gray!65,line width=0.4pt] (0.923,-0.492) -- (1.459,-0.215) -- (1.876,-0.443) -- (1.871,-1.184) -- (1.311,-1.321) -- (0.929,-1.041) -- cycle;
\path[fill=white,draw=gray!65,line width=0.4pt] (1.871,-1.184) -- (1.876,-0.443) -- (2.374,-0.215) -- (2.978,-0.485) -- (2.960,-0.522) -- (2.936,-0.571) -- (2.910,-0.619) -- (2.884,-0.666) -- (2.857,-0.712) -- (2.830,-0.758) -- (2.802,-0.803) -- (2.773,-0.848) -- (2.745,-0.892) -- (2.716,-0.935) -- (2.686,-0.978) -- (2.657,-1.020) -- (2.627,-1.061) -- (2.598,-1.103) -- (2.568,-1.143) -- (2.538,-1.184) -- (2.509,-1.223) -- (2.479,-1.263) -- (2.449,-1.302) -- (2.420,-1.341) -- (2.391,-1.380) -- (2.361,-1.419) -- (2.332,-1.457) -- (2.307,-1.491) -- cycle;
\path[fill=white,draw=gray!65,line width=0.4pt] (-1.943,0.715) -- (-1.344,0.486) -- (-1.403,-0.230) -- (-1.903,-0.399) -- (-2.412,0.084) -- cycle;
\path[fill=PTblue!25,draw=gray!65,line width=0.4pt] (-1.344,0.486) -- (-1.068,0.638) -- (-0.438,0.387) -- (-0.492,-0.263) -- (-0.998,-0.497) -- (-1.403,-0.230) -- cycle;
\path[fill=PTblue!65,draw=gray!65,line width=0.4pt] (-0.128,0.578) -- (0.466,0.287) -- (0.508,-0.236) -- (-0.049,-0.543) -- (-0.492,-0.263) -- (-0.438,0.387) -- cycle;
\path[fill=PTblue!25,draw=gray!65,line width=0.4pt] (0.979,0.671) -- (1.426,0.401) -- (1.459,-0.215) -- (0.923,-0.492) -- (0.508,-0.236) -- (0.466,0.287) -- cycle;
\path[fill=white,draw=gray!65,line width=0.4pt] (1.459,-0.215) -- (1.426,0.401) -- (2.065,0.718) -- (2.308,0.584) -- (2.374,-0.215) -- (1.876,-0.443) -- cycle;
\path[fill=white,draw=gray!65,line width=0.4pt] (2.374,-0.215) -- (2.308,0.584) -- (3.128,0.829) -- (3.139,0.783) -- (3.151,0.727) -- (3.161,0.672) -- (3.169,0.616) -- (3.176,0.560) -- (3.181,0.504) -- (3.184,0.447) -- (3.185,0.391) -- (3.184,0.335) -- (3.182,0.278) -- (3.178,0.222) -- (3.172,0.166) -- (3.165,0.111) -- (3.156,0.055) -- (3.145,0.000) -- (3.132,-0.055) -- (3.119,-0.109) -- (3.103,-0.163) -- (3.086,-0.216) -- (3.068,-0.268) -- (3.049,-0.320) -- (3.028,-0.372) -- (3.007,-0.423) -- (2.984,-0.473) -- (2.978,-0.485) -- cycle;
\path[fill=white,draw=gray!65,line width=0.4pt] (-2.064,1.054) -- (-1.943,0.715) -- (-2.412,0.084) -- (-2.914,0.052) -- (-2.919,0.102) -- (-2.923,0.153) -- (-2.928,0.205) -- (-2.931,0.256) -- (-2.934,0.308) -- (-2.937,0.361) -- (-2.939,0.413) -- (-2.940,0.466) -- (-2.940,0.518) -- (-2.939,0.571) -- (-2.938,0.624) -- (-2.935,0.678) -- (-2.932,0.731) -- (-2.927,0.784) -- (-2.922,0.838) -- (-2.915,0.891) -- (-2.906,0.944) -- (-2.897,0.998) -- (-2.886,1.050) -- (-2.874,1.103) -- (-2.861,1.156) -- (-2.846,1.208) -- (-2.829,1.260) -- (-2.811,1.311) -- (-2.799,1.343) -- cycle;
\path[fill=PTblue!25,draw=gray!65,line width=0.4pt] (-1.943,0.715) -- (-2.064,1.054) -- (-1.562,1.624) -- (-0.976,1.473) -- (-1.068,0.638) -- (-1.344,0.486) -- cycle;
\path[fill=PTblue!25,draw=gray!65,line width=0.4pt] (-0.976,1.473) -- (-0.788,1.576) -- (-0.060,1.132) -- (-0.128,0.578) -- (-0.438,0.387) -- (-1.068,0.638) -- cycle;
\path[fill=PTblue!25,draw=gray!65,line width=0.4pt] (-0.060,1.132) -- (0.330,1.384) -- (0.939,1.035) -- (0.979,0.671) -- (0.466,0.287) -- (-0.128,0.578) -- cycle;
\path[fill=PTblue!25,draw=gray!65,line width=0.4pt] (0.979,0.671) -- (0.939,1.035) -- (1.404,1.510) -- (1.928,1.368) -- (2.065,0.718) -- (1.426,0.401) -- cycle;
\path[fill=white,draw=gray!65,line width=0.4pt] (2.308,0.584) -- (2.065,0.718) -- (1.928,1.368) -- (2.489,1.883) -- (2.493,1.879) -- (2.532,1.839) -- (2.569,1.799) -- (2.607,1.758) -- (2.643,1.717) -- (2.679,1.674) -- (2.714,1.631) -- (2.748,1.587) -- (2.782,1.542) -- (2.814,1.496) -- (2.846,1.450) -- (2.876,1.403) -- (2.905,1.355) -- (2.934,1.306) -- (2.961,1.257) -- (2.986,1.206) -- (3.011,1.156) -- (3.033,1.104) -- (3.055,1.052) -- (3.075,0.999) -- (3.093,0.946) -- (3.110,0.892) -- (3.125,0.837) -- (3.128,0.829) -- cycle;
\path[fill=white,draw=gray!65,line width=0.4pt] (-1.562,1.624) -- (-2.064,1.054) -- (-2.799,1.343) -- (-2.792,1.362) -- (-2.771,1.412) -- (-2.749,1.462) -- (-2.725,1.510) -- (-2.700,1.559) -- (-2.673,1.606) -- (-2.645,1.653) -- (-2.615,1.698) -- (-2.584,1.743) -- (-2.552,1.787) -- (-2.518,1.830) -- (-2.483,1.871) -- (-2.447,1.912) -- (-2.410,1.952) -- (-2.372,1.990) -- (-2.332,2.027) -- (-2.292,2.064) -- (-2.250,2.099) -- (-2.208,2.132) -- (-2.165,2.165) -- (-2.122,2.197) -- (-2.077,2.227) -- (-2.032,2.257) -- (-1.986,2.285) -- (-1.940,2.312) -- (-1.896,2.338) -- cycle;
\path[fill=white,draw=gray!65,line width=0.4pt] (-0.621,2.019) -- (-0.788,1.576) -- (-0.976,1.473) -- (-1.562,1.624) -- (-1.896,2.338) -- (-1.894,2.339) -- (-1.847,2.364) -- (-1.800,2.388) -- (-1.752,2.412) -- (-1.704,2.434) -- (-1.656,2.456) -- (-1.608,2.477) -- (-1.560,2.497) -- (-1.512,2.516) -- (-1.463,2.535) -- (-1.415,2.553) -- (-1.366,2.570) -- (-1.318,2.587) -- (-1.270,2.603) -- (-1.221,2.619) -- (-1.173,2.634) -- (-1.124,2.649) -- (-1.076,2.663) -- (-1.028,2.677) -- (-0.979,2.691) -- (-0.931,2.704) -- (-0.906,2.711) -- cycle;
\path[fill=PTblue!25,draw=gray!65,line width=0.4pt] (-0.788,1.576) -- (-0.621,2.019) -- (0.253,2.140) -- (0.382,1.981) -- (0.330,1.384) -- (-0.060,1.132) -- cycle;
\path[fill=white,draw=gray!65,line width=0.4pt] (0.382,1.981) -- (1.205,2.043) -- (1.404,1.510) -- (0.939,1.035) -- (0.330,1.384) -- cycle;
\path[fill=white,draw=gray!65,line width=0.4pt] (1.928,1.368) -- (1.404,1.510) -- (1.205,2.043) -- (1.485,2.605) -- (1.500,2.598) -- (1.547,2.574) -- (1.593,2.550) -- (1.639,2.524) -- (1.685,2.498) -- (1.731,2.472) -- (1.776,2.444) -- (1.821,2.416) -- (1.865,2.388) -- (1.910,2.358) -- (1.954,2.328) -- (1.997,2.298) -- (2.041,2.266) -- (2.084,2.234) -- (2.126,2.202) -- (2.169,2.169) -- (2.211,2.135) -- (2.252,2.100) -- (2.293,2.065) -- (2.334,2.029) -- (2.375,1.993) -- (2.415,1.955) -- (2.454,1.917) -- (2.489,1.883) -- cycle;
\path[fill=white,draw=gray!65,line width=0.4pt] (0.253,2.140) -- (-0.621,2.019) -- (-0.906,2.711) -- (-0.883,2.717) -- (-0.835,2.730) -- (-0.786,2.742) -- (-0.738,2.754) -- (-0.690,2.766) -- (-0.641,2.777) -- (-0.593,2.788) -- (-0.544,2.798) -- (-0.495,2.809) -- (-0.446,2.818) -- (-0.397,2.828) -- (-0.348,2.836) -- (-0.299,2.845) -- (-0.250,2.852) -- (-0.200,2.860) -- (-0.150,2.866) -- (-0.100,2.872) -- (-0.050,2.878) -- (0.000,2.882) -- (0.050,2.886) -- (0.101,2.889) -- (0.152,2.892) -- (0.202,2.893) -- (0.253,2.894) -- (0.285,2.894) -- cycle;
\path[fill=white,draw=gray!65,line width=0.4pt] (1.205,2.043) -- (0.382,1.981) -- (0.253,2.140) -- (0.285,2.894) -- (0.304,2.894) -- (0.355,2.892) -- (0.406,2.890) -- (0.457,2.887) -- (0.508,2.883) -- (0.560,2.878) -- (0.611,2.873) -- (0.662,2.866) -- (0.712,2.858) -- (0.763,2.849) -- (0.814,2.839) -- (0.864,2.828) -- (0.915,2.816) -- (0.965,2.803) -- (1.015,2.789) -- (1.065,2.774) -- (1.114,2.758) -- (1.163,2.741) -- (1.212,2.723) -- (1.261,2.704) -- (1.309,2.685) -- (1.357,2.664) -- (1.405,2.643) -- (1.453,2.621) -- (1.485,2.605) -- cycle;
\path[pattern=north east lines,pattern color=PTmagenta] (1.157,-2.080) -- (1.311,-1.321) -- (1.871,-1.184) -- (2.307,-1.491) -- (2.303,-1.496) -- (2.274,-1.534) -- (2.246,-1.572) -- (2.217,-1.611) -- (2.188,-1.649) -- (2.160,-1.687) -- (2.131,-1.726) -- (2.102,-1.764) -- (2.073,-1.802) -- (2.044,-1.841) -- (2.015,-1.879) -- (1.986,-1.917) -- (1.956,-1.956) -- (1.926,-1.994) -- (1.895,-2.033) -- (1.865,-2.071) -- (1.833,-2.109) -- (1.802,-2.147) -- (1.769,-2.185) -- (1.736,-2.222) -- (1.703,-2.260) -- (1.669,-2.297) -- (1.634,-2.334) -- (1.598,-2.370) -- (1.562,-2.406) -- (1.525,-2.441) -- (1.487,-2.476) -- (1.476,-2.485) -- cycle;
\path[draw,line width=0.8pt] (3.145,0.000) -- (3.156,0.055) -- (3.165,0.111) -- (3.172,0.166) -- (3.178,0.222) -- (3.182,0.278) -- (3.184,0.335) -- (3.185,0.391) -- (3.184,0.447) -- (3.181,0.504) -- (3.176,0.560) -- (3.169,0.616) -- (3.161,0.672) -- (3.151,0.727) -- (3.139,0.783) -- (3.125,0.837) -- (3.110,0.892) -- (3.093,0.946) -- (3.075,0.999) -- (3.055,1.052) -- (3.033,1.104) -- (3.011,1.156) -- (2.986,1.206) -- (2.961,1.257) -- (2.934,1.306) -- (2.905,1.355) -- (2.876,1.403) -- (2.846,1.450) -- (2.814,1.496) -- (2.782,1.542) -- (2.748,1.587) -- (2.714,1.631) -- (2.679,1.674) -- (2.643,1.717) -- (2.607,1.758) -- (2.569,1.799) -- (2.532,1.839) -- (2.493,1.879) -- (2.454,1.917) -- (2.415,1.955) -- (2.375,1.993) -- (2.334,2.029) -- (2.293,2.065) -- (2.252,2.100) -- (2.211,2.135) -- (2.169,2.169) -- (2.126,2.202) -- (2.084,2.234) -- (2.041,2.266) -- (1.997,2.298) -- (1.954,2.328) -- (1.910,2.358) -- (1.865,2.388) -- (1.821,2.416) -- (1.776,2.444) -- (1.731,2.472) -- (1.685,2.498) -- (1.639,2.524) -- (1.593,2.550) -- (1.547,2.574) -- (1.500,2.598) -- (1.453,2.621) -- (1.405,2.643) -- (1.357,2.664) -- (1.309,2.685) -- (1.261,2.704) -- (1.212,2.723) -- (1.163,2.741) -- (1.114,2.758) -- (1.065,2.774) -- (1.015,2.789) -- (0.965,2.803) -- (0.915,2.816) -- (0.864,2.828) -- (0.814,2.839) -- (0.763,2.849) -- (0.712,2.858) -- (0.662,2.866) -- (0.611,2.873) -- (0.560,2.878) -- (0.508,2.883) -- (0.457,2.887) -- (0.406,2.890) -- (0.355,2.892) -- (0.304,2.894) -- (0.253,2.894) -- (0.202,2.893) -- (0.152,2.892) -- (0.101,2.889) -- (0.050,2.886) -- (0.000,2.882) -- (-0.050,2.878) -- (-0.100,2.872) -- (-0.150,2.866) -- (-0.200,2.860) -- (-0.250,2.852) -- (-0.299,2.845) -- (-0.348,2.836) -- (-0.397,2.828) -- (-0.446,2.818) -- (-0.495,2.809) -- (-0.544,2.798) -- (-0.593,2.788) -- (-0.641,2.777) -- (-0.690,2.766) -- (-0.738,2.754) -- (-0.786,2.742) -- (-0.835,2.730) -- (-0.883,2.717) -- (-0.931,2.704) -- (-0.979,2.691) -- (-1.028,2.677) -- (-1.076,2.663) -- (-1.124,2.649) -- (-1.173,2.634) -- (-1.221,2.619) -- (-1.270,2.603) -- (-1.318,2.587) -- (-1.366,2.570) -- (-1.415,2.553) -- (-1.463,2.535) -- (-1.512,2.516) -- (-1.560,2.497) -- (-1.608,2.477) -- (-1.656,2.456) -- (-1.704,2.434) -- (-1.752,2.412) -- (-1.800,2.388) -- (-1.847,2.364) -- (-1.894,2.339) -- (-1.940,2.312) -- (-1.986,2.285) -- (-2.032,2.257) -- (-2.077,2.227) -- (-2.122,2.197) -- (-2.165,2.165) -- (-2.208,2.132) -- (-2.250,2.099) -- (-2.292,2.064) -- (-2.332,2.027) -- (-2.372,1.990) -- (-2.410,1.952) -- (-2.447,1.912) -- (-2.483,1.871) -- (-2.518,1.830) -- (-2.552,1.787) -- (-2.584,1.743) -- (-2.615,1.698) -- (-2.645,1.653) -- (-2.673,1.606) -- (-2.700,1.559) -- (-2.725,1.510) -- (-2.749,1.462) -- (-2.771,1.412) -- (-2.792,1.362) -- (-2.811,1.311) -- (-2.829,1.260) -- (-2.846,1.208) -- (-2.861,1.156) -- (-2.874,1.103) -- (-2.886,1.050) -- (-2.897,0.998) -- (-2.906,0.944) -- (-2.915,0.891) -- (-2.922,0.838) -- (-2.927,0.784) -- (-2.932,0.731) -- (-2.935,0.678) -- (-2.938,0.624) -- (-2.939,0.571) -- (-2.940,0.518) -- (-2.940,0.466) -- (-2.939,0.413) -- (-2.937,0.361) -- (-2.934,0.308) -- (-2.931,0.256) -- (-2.928,0.205) -- (-2.923,0.153) -- (-2.919,0.102) -- (-2.914,0.051) -- (-2.908,0.000) -- (-2.903,-0.051) -- (-2.897,-0.101) -- (-2.890,-0.151) -- (-2.884,-0.202) -- (-2.877,-0.252) -- (-2.870,-0.302) -- (-2.863,-0.352) -- (-2.855,-0.401) -- (-2.848,-0.451) -- (-2.840,-0.501) -- (-2.832,-0.551) -- (-2.824,-0.600) -- (-2.815,-0.650) -- (-2.807,-0.700) -- (-2.798,-0.750) -- (-2.788,-0.800) -- (-2.779,-0.850) -- (-2.768,-0.900) -- (-2.758,-0.950) -- (-2.747,-1.000) -- (-2.735,-1.050) -- (-2.723,-1.100) -- (-2.710,-1.150) -- (-2.697,-1.201) -- (-2.683,-1.251) -- (-2.668,-1.301) -- (-2.652,-1.351) -- (-2.636,-1.401) -- (-2.618,-1.451) -- (-2.600,-1.501) -- (-2.580,-1.550) -- (-2.560,-1.600) -- (-2.539,-1.649) -- (-2.516,-1.697) -- (-2.493,-1.745) -- (-2.468,-1.793) -- (-2.442,-1.840) -- (-2.416,-1.887) -- (-2.388,-1.933) -- (-2.358,-1.979) -- (-2.328,-2.024) -- (-2.297,-2.068) -- (-2.264,-2.111) -- (-2.230,-2.154) -- (-2.196,-2.196) -- (-2.160,-2.236) -- (-2.123,-2.276) -- (-2.085,-2.315) -- (-2.046,-2.354) -- (-2.006,-2.391) -- (-1.965,-2.427) -- (-1.923,-2.462) -- (-1.881,-2.496) -- (-1.837,-2.529) -- (-1.793,-2.561) -- (-1.748,-2.592) -- (-1.702,-2.621) -- (-1.656,-2.650) -- (-1.609,-2.678) -- (-1.561,-2.704) -- (-1.513,-2.730) -- (-1.465,-2.754) -- (-1.415,-2.778) -- (-1.366,-2.801) -- (-1.316,-2.822) -- (-1.266,-2.843) -- (-1.215,-2.862) -- (-1.164,-2.881) -- (-1.113,-2.898) -- (-1.061,-2.915) -- (-1.009,-2.931) -- (-0.957,-2.946) -- (-0.905,-2.960) -- (-0.853,-2.973) -- (-0.800,-2.985) -- (-0.747,-2.997) -- (-0.694,-3.008) -- (-0.641,-3.017) -- (-0.588,-3.026) -- (-0.535,-3.034) -- (-0.482,-3.041) -- (-0.428,-3.048) -- (-0.375,-3.053) -- (-0.321,-3.057) -- (-0.268,-3.061) -- (-0.214,-3.064) -- (-0.161,-3.065) -- (-0.107,-3.066) -- (-0.054,-3.066) -- (-0.000,-3.065) -- (0.053,-3.062) -- (0.107,-3.059) -- (0.160,-3.055) -- (0.213,-3.049) -- (0.266,-3.043) -- (0.319,-3.035) -- (0.372,-3.026) -- (0.424,-3.017) -- (0.476,-3.005) -- (0.528,-2.993) -- (0.579,-2.980) -- (0.630,-2.965) -- (0.681,-2.949) -- (0.731,-2.932) -- (0.781,-2.914) -- (0.830,-2.894) -- (0.879,-2.874) -- (0.927,-2.852) -- (0.974,-2.829) -- (1.021,-2.805) -- (1.067,-2.780) -- (1.112,-2.754) -- (1.157,-2.726) -- (1.201,-2.698) -- (1.244,-2.669) -- (1.287,-2.639) -- (1.329,-2.608) -- (1.369,-2.576) -- (1.410,-2.543) -- (1.449,-2.510) -- (1.487,-2.476) -- (1.525,-2.441) -- (1.562,-2.406) -- (1.598,-2.370) -- (1.634,-2.334) -- (1.669,-2.297) -- (1.703,-2.260) -- (1.736,-2.222) -- (1.769,-2.185) -- (1.802,-2.147) -- (1.833,-2.109) -- (1.865,-2.071) -- (1.895,-2.033) -- (1.926,-1.994) -- (1.956,-1.956) -- (1.986,-1.917) -- (2.015,-1.879) -- (2.044,-1.841) -- (2.073,-1.802) -- (2.102,-1.764) -- (2.131,-1.726) -- (2.160,-1.687) -- (2.188,-1.649) -- (2.217,-1.611) -- (2.246,-1.572) -- (2.274,-1.534) -- (2.303,-1.496) -- (2.332,-1.457) -- (2.361,-1.419) -- (2.391,-1.380) -- (2.420,-1.341) -- (2.449,-1.302) -- (2.479,-1.263) -- (2.509,-1.223) -- (2.538,-1.184) -- (2.568,-1.143) -- (2.598,-1.103) -- (2.627,-1.061) -- (2.657,-1.020) -- (2.686,-0.978) -- (2.716,-0.935) -- (2.745,-0.892) -- (2.773,-0.848) -- (2.802,-0.803) -- (2.830,-0.758) -- (2.857,-0.712) -- (2.884,-0.666) -- (2.910,-0.619) -- (2.936,-0.571) -- (2.960,-0.522) -- (2.984,-0.473) -- (3.007,-0.423) -- (3.028,-0.372) -- (3.049,-0.320) -- (3.068,-0.268) -- (3.086,-0.216) -- (3.103,-0.163) -- (3.119,-0.109) -- (3.132,-0.055) -- cycle;
\draw[dashed] (0,0) circle (2.590);
\draw[PTblue,thick] (0,0) circle (1.295);
\fill (0,0) circle (1.5pt) node[below left=-1pt] {$x_P$};
\node at (-2.75,2.75) {$P$};
\draw[->] (4.0,2.75) node[right]{$U_P$} -- (2.108,1.505);
\draw[->] (4.0,1.75) node[right]{$U(x_P,\tfrac12c_b\side P)$} -- (1.216,0.444);
\draw[->] (4.0,0.55) node[right]{central child $P^\downarrow$} -- (0.287,0.072);
\draw[->] (4.0,-0.6) node[right]{a child in $\Inner(P)$} -- (0.436,-0.792);
\draw[->] (4.0,-2.1) node[right,align=left]{a child $R$ with $\gap(R,\XX\setminus P)\leq 9\side R$,\\ which lies in $\Outer(P)$} -- (1.653,-1.751);
\end{tikzpicture}
\caption{Inner and outer children in Lemma \ref{l:vAR-count}. The children in $\Inner(P)$ are shown in blue, and the central child $P^\downarrow$ in a darker shade.}\label{fig:inner-outer}
\end{figure}

Figure \ref{fig:inner-outer} illustrates the sets in Lemma \ref{l:vAR-count}.

\subsection{Sparse jumps}\label{ss:variable}

For the rest of the paper, when variable scaling is needed, we use one particular choice of side lengths.

\begin{definition}[sparse jumps]\label{sparse-cubes}
Let $b>3$. Put $\xi_1:=b$ and $\xi_j:=b\lceil\log_2j\rceil$ for all $j\geq2$, and put
\begin{equation}\label{sparse-scales}
 r_k:=\begin{cases}
 b^{-k},&k\leq 0,\\
 (\xi_1\cdots\xi_k)^{-1},&k\geq 1.
 \end{cases}
\end{equation}
A system of cubes $(\Delta_k)_{k\in\ZZ}$ with side lengths $(r_k)_{k\in\ZZ}$, in the sense of Theorem \ref{vKRS-cubes} and Remark \ref{r:KRS}, is called a system of \emph{$b$-adic cubes with sparse jumps}.
\end{definition}

To distinguish the two systems, we refer to a system of $b$-adic cubes in the sense of Remark \ref{r:KRS}, for which $\xi_k=b$ for all $k$, as a system of \emph{$b$-adic cubes without jumps}.

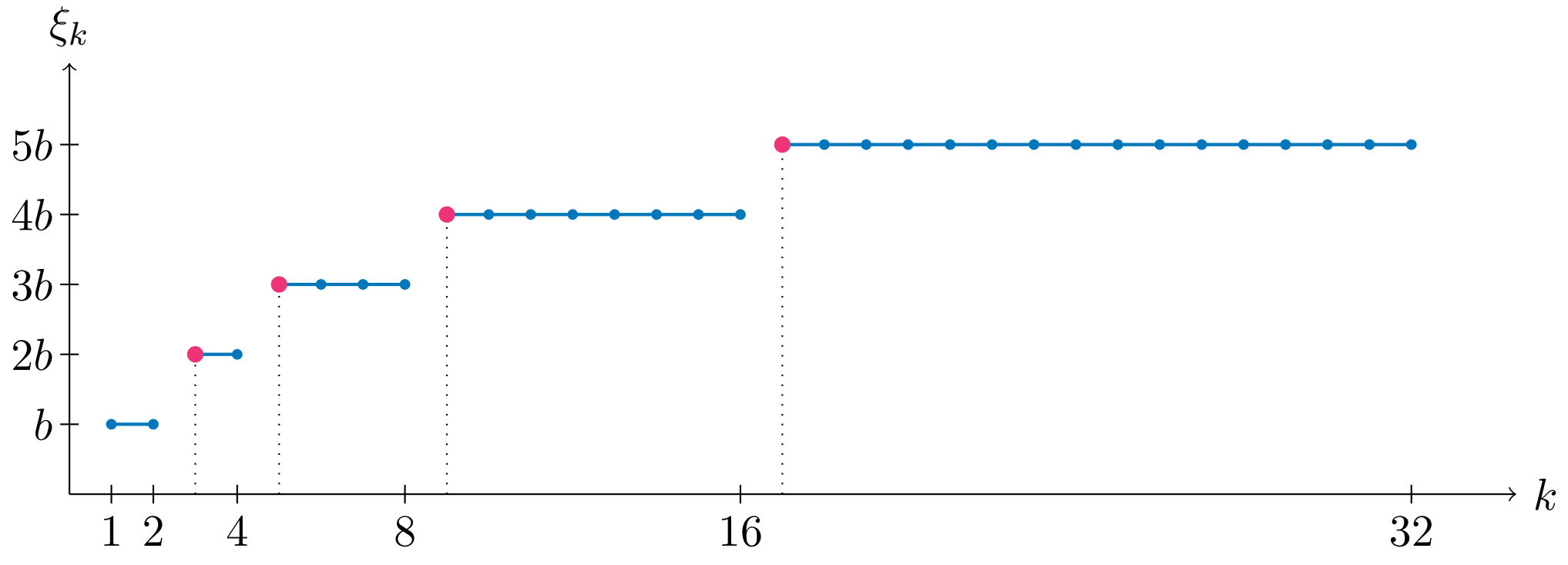
\begin{figure}[htb]
\centering
\begin{tikzpicture}[font=\small,x=0.36cm,y=1cm]
\draw[->] (0,0)--(34.5,0) node[right]{$k$};
\draw[->] (0,0)--(0,3.7) node[above]{$\xi_k$};
\foreach \y/\l in {0.6/b,1.2/2b,1.8/3b,2.4/4b,3.0/5b} {\draw (-0.22,\y)--(0.22,\y); \node[left] at (0,\y) {$\l$};}
\foreach \a/\b/\y in {1/2/0.6,3/4/1.2,5/8/1.8,9/16/2.4,17/32/3.0} {\draw[thick,PTblue] (\a,\y)--(\b,\y);}
\foreach \k in {1,2} {\fill[PTblue] (\k,0.6) circle (1.3pt);}
\foreach \k in {3,4} {\fill[PTblue] (\k,1.2) circle (1.3pt);}
\foreach \k in {5,...,8} {\fill[PTblue] (\k,1.8) circle (1.3pt);}
\foreach \k in {9,...,16} {\fill[PTblue] (\k,2.4) circle (1.3pt);}
\foreach \k in {17,...,32} {\fill[PTblue] (\k,3.0) circle (1.3pt);}
\foreach \k/\y in {3/1.2,5/1.8,9/2.4,17/3.0} {\draw[dotted] (\k,0)--(\k,\y); \fill[PTmagenta] (\k,\y) circle (2pt);}
\foreach \k in {1,2,4,8,16,32} {\draw (\k,-0.08)--(\k,0.08); \node[below] at (\k,-0.05) {$\k$};}
\end{tikzpicture}
\caption{Ratios of consecutive side lengths for a system of $b$-adic cubes with sparse jumps, for the generations $1\leq k\leq 32$. The jumps at generations $2^h+1$ are marked in magenta.}\label{fig:sparse}
\end{figure}

\begin{remark}[plateaus]\label{r:plateau}
Let $\Delta$ be a system of $b$-adic cubes with sparse jumps. For every integer $h\geq1$, we call the set of generations $k\geq1$ with $\xi_k=bh$ the $h$-th \emph{plateau} of $\Delta$. The first plateau consists of the generations $1$ and $2$, and for $h\geq2$, the $h$-th plateau consists of the $2^{h-1}$ consecutive generations $2^{h-1}<k\leq 2^h$. The ratio $\xi_k$ changes value only in the transition from $k=2^h$ to $k+1$, and we call the generations $2^h+1$ with $h\geq 1$ the \emph{jumps} of the system. Thus, $\xi_k\to\infty$, while $\xi_k$ is constant on the long plateaus between consecutive jumps; see Figure \ref{fig:sparse}.
\end{remark}

In the computation of local dimensions in the proof of Theorem \ref{t:dimension}, we use the fact that the ratios $\xi_k$ are nearly constant on average.

\begin{lemma}[ratios on average]\label{l:sparse-lag}
Let $\Delta$ be a system of $b$-adic cubes with sparse jumps. If $h\geq 2$ and $2^{h-1}<n\leq 2^h$, then
\begin{equation}\label{sparse-lag}
 \sum_{k=1}^n\log\frac{\xi_n}{\xi_k}\leq\frac{2^{h+2}}{h}+2^{(h+1)/2}\log h.
\end{equation}
Moreover,
\begin{equation}\label{sparse-lag-limit}
 \lim_{n\to\infty}\frac1n\sum_{k=1}^n\log\frac{\xi_n}{\xi_k}=0.
\end{equation}
\end{lemma}

\begin{proof}
Let $h\geq 2$ and $2^{h-1}<n\leq 2^h$. Then $\xi_n=bh$, and the generations $k\leq n$ with $\xi_k=bh$ contribute nothing to the sums below. For $1\leq i\leq h-1$, the generations $k\leq n$ with $\xi_k=bi$ form the $i$-th plateau, which contains at most $2^i$ generations. Hence
\[
 \sum_{k=1}^n\log\frac{\xi_n}{\xi_k}\leq\sum_{i=1}^{h-1}2^i\log\frac hi.
\]
If $h/2\leq i\leq h-1$, then $\log(h/i)\leq(h-i)/i\leq 2(h-i)/h$. Writing $j=h-i$ and using $\sum_{j=1}^\infty j2^{-j}=2$, we obtain
\[
 \sum_{h/2\leq i\leq h-1}2^i\log\frac hi\leq\frac2h\sum_{j=1}^\infty j2^{h-j}=\frac{2^{h+2}}h.
\]
If $1\leq i<h/2$, then $i\leq(h-1)/2$ and $\log(h/i)\leq\log h$. Thus,
\[
 \sum_{1\leq i<h/2}2^i\log\frac hi\leq\log h\sum_{i=1}^{\lfloor(h-1)/2\rfloor}2^i<2^{(h+1)/2}\log h.
\]
Adding the last two estimates gives \eqref{sparse-lag}. Since $\xi_k\leq\xi_n$ for all $k\leq n$ and $n>2^{h-1}$, dividing \eqref{sparse-lag} by $n$ gives
\[
 0\leq\frac1n\sum_{k=1}^n\log\frac{\xi_n}{\xi_k}\leq\frac8h+2^{(3-h)/2}\log h.
\]
Since $h\to\infty$ as $n\to\infty$, this yields \eqref{sparse-lag-limit}.
\end{proof}

\subsection{Measures from mass distributions}\label{ss:mass-distribution}

In \S\S\ref{s:entropy}--\ref{s:zero}, we define measures by prescribing the masses of cubes. To verify that the prescribed masses are the values of a Radon measure, we use two results from \cite{Badger-MethodI}. A \emph{Radon measure} on a proper metric space is a Borel regular outer measure that is finite on bounded sets.

\begin{definition}[nested partitions]\label{d:nested-partitions}
Let $\XX$ be a proper metric space. A sequence $(\Delta_k)_{k\geq0}$ of families of subsets of $\XX$ is a \emph{sequence of nested partitions} if the following conditions hold.
\begin{enumerate}
\item Every $\Delta_k$ is a partition of $\XX$ into nonempty bounded Borel sets.
\item Every set in $\Delta_{k+1}$ is contained in a set in $\Delta_k$.
\item Every bounded set in $\XX$ intersects only finitely many sets in $\Delta_k$ for every $k\geq 0$.
\item $\sup_{Q\in\Delta_k}\diam Q\to0$ as $k\to\infty$.
\end{enumerate}
\end{definition}

Let $(\Delta_k)_{k\geq0}$ be a sequence of nested partitions. We call the sets in $\mathcal T:=\bigcup_{k\geq0}\Delta_k$ cubes. For $Q\in\Delta_k$ and $n\geq0$, let $\Child^n(Q)$ be the finite set of cubes in $\Delta_{k+n}$ contained in $Q$, put $\Child(Q):=\Child^1(Q)$, and put
\[
 \Collar^n(Q):=\{R\in\Child^n(Q):\overline R\cap\partial Q\neq\emptyset\}.
\]
The same set may belong to $\Delta_k$ for several $k$. In that case, $\Child^n(Q)$ and $\Collar^n(Q)$ are understood relative to a specified \emph{occurrence} $Q\in\Delta_k$.

\begin{definition}[mass distributions]\label{d:mass-distribution}
Let $(\Delta_k)_{k\geq0}$ be a sequence of nested partitions of a proper metric space $\XX$. A \emph{mass distribution} is a function $w:\mathcal T\to[0,\infty)$ such that
\begin{equation}\label{mass-additivity}
 w(Q)=\sum_{R\in\Child(Q)}w(R)\quad\text{for all }Q\in\Delta_k\text{ and }k\geq0.
\end{equation}
A mass distribution $w$ is \emph{tame} if
\begin{equation}\label{tame-condition}
 \lim_{n\to\infty}\sum_{R\in\Collar^n(Q)}w(R)=0\quad\text{for every occurrence }Q\in\Delta_k.
\end{equation}
Given a mass distribution $w$, put $w(\emptyset):=0$ and
\[
 \mu_w(A):=\inf\left\{\sum_{i=1}^\infty w(Q_i):Q_i\in\mathcal T\cup\{\emptyset\}\text{ and }A\subset\bigcup_{i=1}^\infty Q_i\right\}\quad\text{for all }A\subset\XX.
\]
We call $\mu_w$ the \emph{Method I outer measure} with weight $w$.
\end{definition}

\begin{theorem}[{\cite[Theorem 6.3]{Badger-MethodI}}]\label{t:existence}
Let $(\Delta_k)_{k\geq0}$ be a sequence of nested partitions of a proper metric space $\XX$, and let $w$ be a mass distribution. Then $\mu_w$ is a Radon measure on $\XX$. Moreover, $w$ is tame if and only if
\begin{equation}\label{tame-values}
 \mu_w(Q)=w(Q)\quad\text{and}\quad\mu_w(\partial Q)=0\quad\text{for all }Q\in\mathcal T.
\end{equation}
\end{theorem}

In the following corollary, $\gap$ is given by \eqref{gap-def}, with $\gap(A,B)=\infty$ when $A$ or $B$ is empty.

\begin{corollary}[{\cite[Corollary 6.5]{Badger-MethodI}}]\label{c:tame}
Let $w$ be a mass distribution for a sequence of nested partitions of a proper metric space $\XX$. Suppose that there is $\eta>0$ such that, for every occurrence $Q\in\Delta_k$, there is a family $\mathcal I(Q)\subset\Child(Q)$ with
\[
 \gap\Big(\bigcup\mathcal I(Q),\XX\setminus Q\Big)>0\quad\text{and}\quad\sum_{R\in\mathcal I(Q)}w(R)\geq\eta\,w(Q).
\]
Then
\begin{equation}\label{tame-geometric}
 \sum_{R\in\Collar^n(Q)}w(R)\leq(1-\eta)^n w(Q)\quad\text{for all }Q\in\mathcal T\text{ and }n\geq 0.
\end{equation}
In particular, $w$ is tame, and \eqref{tame-values} holds.
\end{corollary}

Corollary \ref{c:tame} extends \cite[Lemma 4.14]{BS4}, which treats systems of $b$-adic cubes with $b>5$ and takes $\mathcal I(Q)$ to be the central child of $Q$.

\section{Measures with prescribed entropy}\label{s:entropy}

In this section, we prescribe the entropy of measures at every generation of a system of cubes and compute their local dimensions. The construction in Definition \ref{d:central} is used in \S\ref{s:zero} and a refined construction that produces doubling measures is the subject of \S\ref{s:doubling}. The dimension estimate rests on the strong law of large numbers, see e.g.~\cite{dimension-of-measures}, \cite{BS4}. The random variables in \cite[Lemma 5.6]{BS4} are bounded and uncorrelated, whereas here they are not uniformly bounded once the target entropies tend to zero. Instead, to use the strong law of large numbers, we bound their conditional second moments using only the number of children of a cube.

Scaling the metric on a space $\XX$ does not change Ahlfors regularity, local dimensions, the doubling property, or the images of Lipschitz maps. Thus, throughout \S\S\ref{s:entropy}--\ref{s:zero}, we can and do adopt the following convention, which allows us to focus on generations $k\geq 0$.

\begin{convention}[standing assumptions]\label{conv:standing}
Let $\XX$ be a complete Ahlfors $q$-regular metric space with $q>0$ and $\diam\XX\geq 2.1$. Fix an Ahlfors $q$-regular measure $\nu$ on $\XX$ satisfying \eqref{AR-background}, and let $D_\nu\geq 1$ and $C_\XX=D_\nu 6^q>1$ be the constants in Lemma \ref{l:vAR-count}. Put
\begin{equation}\label{K-and-C2}
 K:=D_\nu^2\,72^q\quad\text{and}\quad C_{\mathrm{ent}}:=C_\XX K^2.
\end{equation}
Let $b\geq\max\{47,C_\XX\}$ be large enough that
\begin{equation}\label{b-entropy}
 \log_bC_{\mathrm{ent}}\leq q/4,
\end{equation}
where the value of $b$ is chosen in each application. Let $(\Delta_k)_{k\in\ZZ}$ be a system of $b$-adic cubes without jumps or with sparse jumps, with side lengths $(r_k)_{k\in\ZZ}$.
\end{convention}

The constant $K$ bounds the Type II weights in Lemma \ref{l:flat}, and $C_{\mathrm{ent}}$ bounds the loss of entropy in Lemma \ref{l:entropy-allocation} and in the ranges \eqref{direct-range} and \eqref{endpoint-entropy-range} of target entropies.

In both cases, $r_0=1$. Every cube $Q$ in $\Delta_k$ with $k\geq0$ has diameter at most $2C_br_k\leq 2C_b<2.1$, and hence $Q\neq\XX$. Since $C_\XX^{-1}\xi_k^q\geq C_\XX^{-1}b^q\geq C_\XX^3>1$ by \eqref{b-entropy}, Lemma \ref{l:vAR-count} gives
\begin{equation}\label{child-range}
 2\leq\#\Inner(P)\leq\#\Child(P)\leq C_\XX\xi_k^q\quad\text{for all }P\in\Delta_{k-1}\text{ and }k\geq1.
\end{equation}

\begin{definition}[target entropies]\label{d:target}
Let $\mathbf s=(s_k)_{k\geq1}$ be a sequence of nonnegative numbers. A Radon measure $\muse$ on $\XX$ with $0<\muse(Q)<\infty$ for every cube $Q$ has \emph{target entropies} $\mathbf s$ if
\begin{equation}\label{target-entropy}
 \sum_{R\in\Child(P)}\frac{\muse(R)}{\muse(P)}\log_{\xi_k}\frac{\muse(P)}{\muse(R)}=s_k\quad\text{for all }P\in\Delta_{k-1}\text{ and }k\geq1,
\end{equation} where $\xi_k= b$ when the system of cubes $\Delta$ has no jumps, and $\xi_1=b$ and $\xi_k=b\lceil \log_2 k\rceil$ ($k\geq 2$) when the system of cubes $\Delta$ has sparse jumps.
\end{definition}

For a sequence $\mathbf s=(s_k)_{k\geq1}$ and $n\geq1$, we write
\begin{equation}\label{average-def}
 A_n(\mathbf s):=\frac1n\sum_{k=1}^ns_k
\end{equation}
for the average of the first $n$ terms.

Let $\muse$ be as in Definition \ref{d:target}. For $x\in\XX$ and $k\in\ZZ$, let $Q_k(x)$ be the cube in $\Delta_k$ containing $x$, and put
\begin{equation}\label{endpoint-Xk}
 X_k(x):=\log_{\xi_k}\frac{\muse(Q_{k-1}(x))}{\muse(Q_k(x))}\quad\text{for all }k\geq 1.
\end{equation}
The function $X_k$ is constant on every cube in $\Delta_k$, and for every $P\in\Delta_{k-1}$, the integral $\int_PX_k\,d\muse$ equals $\muse(P)$ times the left side of \eqref{target-entropy}. Thus, \eqref{target-entropy} is equivalent to
\begin{equation}\label{target-integral}
 \int_PX_k\,d\muse=s_k\,\muse(P)\quad\text{for all }P\in\Delta_{k-1}\text{ and }k\geq1.
\end{equation}
To apply the strong law of large numbers and related estimates, we use the following probabilistic notation. Fix $Q_0\in\Delta_0$, normalize $\muse\res Q_0$ to be a probability measure $\PP$ with expectation $\EE$, and let $\mathscr F_k$ be the $\sigma$-algebra generated by the cubes in $\Delta_k$ contained in $Q_0$. Then \eqref{target-integral} says that $\EE(X_k\mid\mathscr F_{k-1})=s_k$ for all $k\geq1$.

Our first construction produces measures with any target entropies in a range that depends only on $\XX$, $b$, and the system of cubes. It divides the mass of every cube between its central child and its remaining children.

\begin{lemma}[central entropy]\label{l:central-entropy}
Let $a>1$ and let $N\geq 2$ be an integer. For all $0<v\leq(N-1)/N$, put
\[
 h_{a,N}(v):=(1-v)\log_a\frac1{1-v}+v\log_a\frac{N-1}{v}.
\]
The function $h_{a,N}$ is continuous and strictly increasing on $(0,(N-1)/N]$, $h_{a,N}(v)\to 0$ as $v\downarrow0$, and $h_{a,N}((N-1)/N)=\log_aN$.
\end{lemma}

\begin{proof}
The assertions follow from $h_{a,N}'(v)=\log_a((N-1)(1-v)/v)>0$ for all $0<v<(N-1)/N$.
\end{proof}

\begin{definition}[central weights]\label{d:central}
Let $\mathbf s=(s_k)_{k\geq1}$ satisfy
\begin{equation}\label{direct-range}
 0<s_k\leq q-\log_{\xi_k}C_{\mathrm{ent}}\quad\text{for all }k\geq1.
\end{equation}
For $k\geq1$ and $P\in\Delta_{k-1}$, put $N_P:=\#\Child(P)$, and let $0<v_P\leq(N_P-1)/N_P$ be the unique number with $h_{\xi_k,N_P}(v_P)=s_k$. The \emph{central weights} with target $s_k$ on $\Child(P)$ are
\[
 p_P(P^\downarrow):=1-v_P\quad\text{and}\quad p_P(R):=\frac{v_P}{N_P-1}\quad\text{for all }R\in\Child(P)\setminus\{P^\downarrow\},
\]
which sum to $1$.
Declare $\muse(Q):=\nu(Q)$ for all $Q\in\Delta_0$, and recursively declare $\muse(R):=p_P(R)\muse(P)$ for all $R\in\Child(P)$, $P\in\Delta_{k-1}$, and $k\geq1$.
\end{definition}

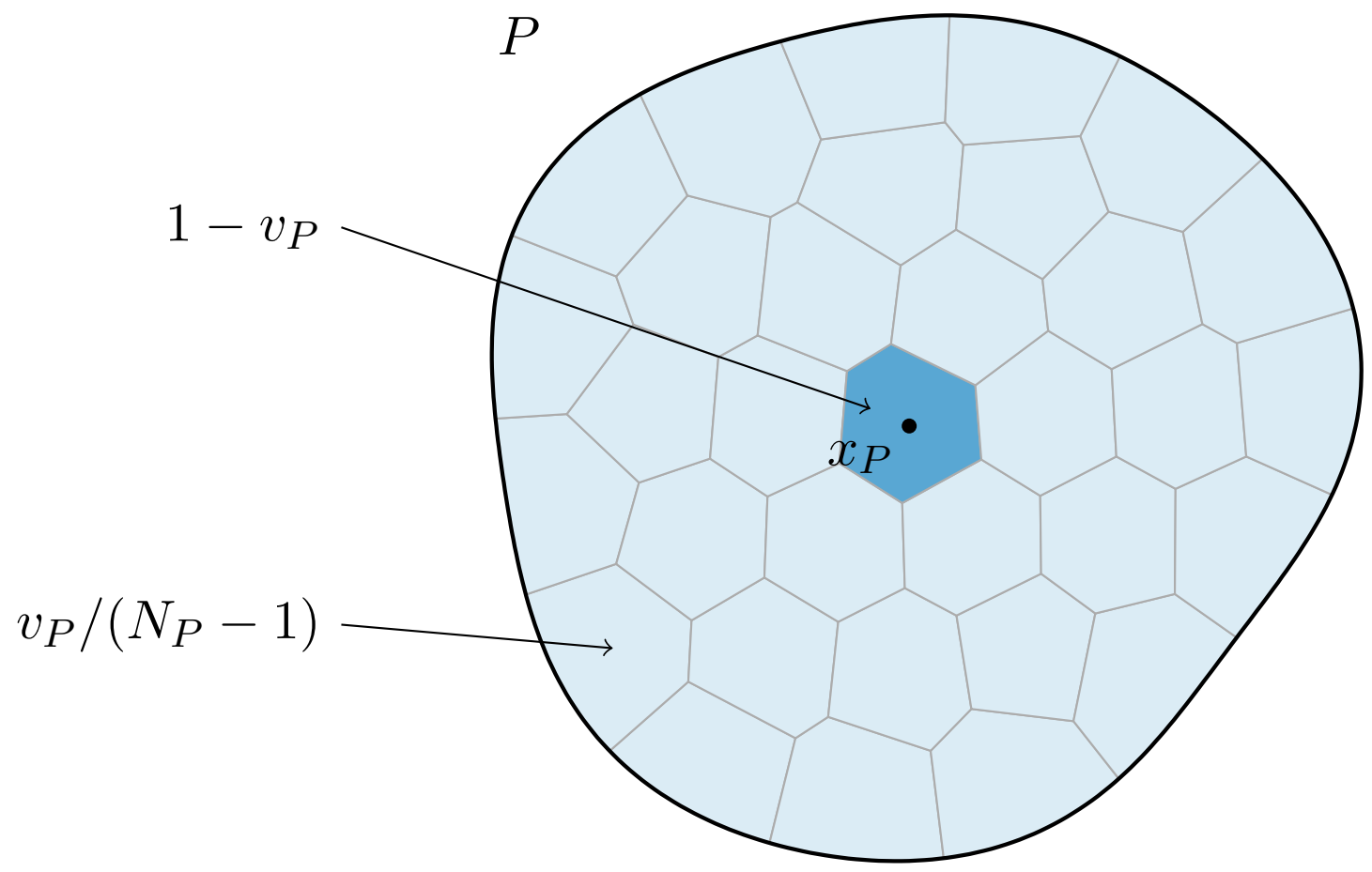
\begin{figure}[htb]
\centering
\begin{tikzpicture}[font=\small]
\path[fill=PTblue!14,draw=gray!65,line width=0.4pt] (-1.556,-1.802) -- (-0.801,-2.200) -- (-0.984,-2.938) -- (-1.009,-2.931) -- (-1.061,-2.915) -- (-1.113,-2.898) -- (-1.164,-2.881) -- (-1.215,-2.862) -- (-1.266,-2.843) -- (-1.316,-2.822) -- (-1.366,-2.801) -- (-1.415,-2.778) -- (-1.465,-2.754) -- (-1.513,-2.730) -- (-1.561,-2.704) -- (-1.609,-2.678) -- (-1.656,-2.650) -- (-1.702,-2.621) -- (-1.748,-2.592) -- (-1.793,-2.561) -- (-1.837,-2.529) -- (-1.881,-2.496) -- (-1.923,-2.462) -- (-1.965,-2.427) -- (-2.006,-2.391) -- (-2.046,-2.354) -- (-2.085,-2.315) -- (-2.110,-2.290) -- cycle;
\path[fill=PTblue!14,draw=gray!65,line width=0.4pt] (-0.801,-2.200) -- (-0.571,-2.050) -- (0.151,-2.289) -- (0.243,-3.046) -- (0.213,-3.049) -- (0.160,-3.055) -- (0.107,-3.059) -- (0.053,-3.062) -- (-0.000,-3.065) -- (-0.054,-3.066) -- (-0.107,-3.066) -- (-0.161,-3.065) -- (-0.214,-3.064) -- (-0.268,-3.061) -- (-0.321,-3.057) -- (-0.375,-3.053) -- (-0.428,-3.048) -- (-0.482,-3.041) -- (-0.535,-3.034) -- (-0.588,-3.026) -- (-0.641,-3.017) -- (-0.694,-3.008) -- (-0.747,-2.997) -- (-0.800,-2.985) -- (-0.853,-2.973) -- (-0.905,-2.960) -- (-0.957,-2.946) -- (-0.984,-2.938) -- cycle;
\path[fill=PTblue!14,draw=gray!65,line width=0.4pt] (0.151,-2.289) -- (0.438,-1.994) -- (1.157,-2.080) -- (1.476,-2.485) -- (1.449,-2.510) -- (1.410,-2.543) -- (1.369,-2.576) -- (1.329,-2.608) -- (1.287,-2.639) -- (1.244,-2.669) -- (1.201,-2.698) -- (1.157,-2.726) -- (1.112,-2.754) -- (1.067,-2.780) -- (1.021,-2.805) -- (0.974,-2.829) -- (0.927,-2.852) -- (0.879,-2.874) -- (0.830,-2.894) -- (0.781,-2.914) -- (0.731,-2.932) -- (0.681,-2.949) -- (0.630,-2.965) -- (0.579,-2.980) -- (0.528,-2.993) -- (0.476,-3.005) -- (0.424,-3.017) -- (0.372,-3.026) -- (0.319,-3.035) -- (0.266,-3.043) -- (0.243,-3.046) -- cycle;
\path[fill=PTblue!14,draw=gray!65,line width=0.4pt] (-2.064,-0.972) -- (-1.533,-1.370) -- (-1.556,-1.802) -- (-2.110,-2.290) -- (-2.123,-2.276) -- (-2.160,-2.236) -- (-2.196,-2.196) -- (-2.230,-2.154) -- (-2.264,-2.111) -- (-2.297,-2.068) -- (-2.328,-2.024) -- (-2.358,-1.979) -- (-2.388,-1.933) -- (-2.416,-1.887) -- (-2.442,-1.840) -- (-2.468,-1.793) -- (-2.493,-1.745) -- (-2.516,-1.697) -- (-2.539,-1.649) -- (-2.560,-1.600) -- (-2.580,-1.550) -- (-2.600,-1.501) -- (-2.618,-1.451) -- (-2.636,-1.401) -- (-2.652,-1.351) -- (-2.668,-1.301) -- (-2.683,-1.251) -- (-2.697,-1.201) -- (-2.700,-1.188) -- cycle;
\path[fill=PTblue!14,draw=gray!65,line width=0.4pt] (-1.533,-1.370) -- (-1.020,-1.068) -- (-0.500,-1.380) -- (-0.571,-2.050) -- (-0.801,-2.200) -- (-1.556,-1.802) -- cycle;
\path[fill=PTblue!14,draw=gray!65,line width=0.4pt] (-0.500,-1.380) -- (-0.031,-1.142) -- (0.333,-1.336) -- (0.438,-1.994) -- (0.151,-2.289) -- (-0.571,-2.050) -- cycle;
\path[fill=PTblue!14,draw=gray!65,line width=0.4pt] (0.929,-1.041) -- (1.311,-1.321) -- (1.157,-2.080) -- (0.438,-1.994) -- (0.333,-1.336) -- cycle;
\path[fill=PTblue!14,draw=gray!65,line width=0.4pt] (1.157,-2.080) -- (1.311,-1.321) -- (1.871,-1.184) -- (2.307,-1.491) -- (2.303,-1.496) -- (2.274,-1.534) -- (2.246,-1.572) -- (2.217,-1.611) -- (2.188,-1.649) -- (2.160,-1.687) -- (2.131,-1.726) -- (2.102,-1.764) -- (2.073,-1.802) -- (2.044,-1.841) -- (2.015,-1.879) -- (1.986,-1.917) -- (1.956,-1.956) -- (1.926,-1.994) -- (1.895,-2.033) -- (1.865,-2.071) -- (1.833,-2.109) -- (1.802,-2.147) -- (1.769,-2.185) -- (1.736,-2.222) -- (1.703,-2.260) -- (1.669,-2.297) -- (1.634,-2.334) -- (1.598,-2.370) -- (1.562,-2.406) -- (1.525,-2.441) -- (1.487,-2.476) -- (1.476,-2.485) -- cycle;
\path[fill=PTblue!14,draw=gray!65,line width=0.4pt] (-2.412,0.084) -- (-1.903,-0.399) -- (-2.064,-0.972) -- (-2.700,-1.188) -- (-2.710,-1.150) -- (-2.723,-1.100) -- (-2.735,-1.050) -- (-2.747,-1.000) -- (-2.758,-0.950) -- (-2.768,-0.900) -- (-2.779,-0.850) -- (-2.788,-0.800) -- (-2.798,-0.750) -- (-2.807,-0.700) -- (-2.815,-0.650) -- (-2.824,-0.600) -- (-2.832,-0.551) -- (-2.840,-0.501) -- (-2.848,-0.451) -- (-2.855,-0.401) -- (-2.863,-0.352) -- (-2.870,-0.302) -- (-2.877,-0.252) -- (-2.884,-0.202) -- (-2.890,-0.151) -- (-2.897,-0.101) -- (-2.903,-0.051) -- (-2.908,0.000) -- (-2.914,0.051) -- (-2.914,0.052) -- cycle;
\path[fill=PTblue!14,draw=gray!65,line width=0.4pt] (-1.903,-0.399) -- (-1.403,-0.230) -- (-0.998,-0.497) -- (-1.020,-1.068) -- (-1.533,-1.370) -- (-2.064,-0.972) -- cycle;
\path[fill=PTblue!14,draw=gray!65,line width=0.4pt] (-0.492,-0.263) -- (-0.049,-0.543) -- (-0.031,-1.142) -- (-0.500,-1.380) -- (-1.020,-1.068) -- (-0.998,-0.497) -- cycle;
\path[fill=PTblue!14,draw=gray!65,line width=0.4pt] (0.508,-0.236) -- (0.923,-0.492) -- (0.929,-1.041) -- (0.333,-1.336) -- (-0.031,-1.142) -- (-0.049,-0.543) -- cycle;
\path[fill=PTblue!14,draw=gray!65,line width=0.4pt] (0.923,-0.492) -- (1.459,-0.215) -- (1.876,-0.443) -- (1.871,-1.184) -- (1.311,-1.321) -- (0.929,-1.041) -- cycle;
\path[fill=PTblue!14,draw=gray!65,line width=0.4pt] (1.871,-1.184) -- (1.876,-0.443) -- (2.374,-0.215) -- (2.978,-0.485) -- (2.960,-0.522) -- (2.936,-0.571) -- (2.910,-0.619) -- (2.884,-0.666) -- (2.857,-0.712) -- (2.830,-0.758) -- (2.802,-0.803) -- (2.773,-0.848) -- (2.745,-0.892) -- (2.716,-0.935) -- (2.686,-0.978) -- (2.657,-1.020) -- (2.627,-1.061) -- (2.598,-1.103) -- (2.568,-1.143) -- (2.538,-1.184) -- (2.509,-1.223) -- (2.479,-1.263) -- (2.449,-1.302) -- (2.420,-1.341) -- (2.391,-1.380) -- (2.361,-1.419) -- (2.332,-1.457) -- (2.307,-1.491) -- cycle;
\path[fill=PTblue!14,draw=gray!65,line width=0.4pt] (-1.943,0.715) -- (-1.344,0.486) -- (-1.403,-0.230) -- (-1.903,-0.399) -- (-2.412,0.084) -- cycle;
\path[fill=PTblue!14,draw=gray!65,line width=0.4pt] (-1.344,0.486) -- (-1.068,0.638) -- (-0.438,0.387) -- (-0.492,-0.263) -- (-0.998,-0.497) -- (-1.403,-0.230) -- cycle;
\path[fill=PTblue!65,draw=gray!65,line width=0.4pt] (-0.128,0.578) -- (0.466,0.287) -- (0.508,-0.236) -- (-0.049,-0.543) -- (-0.492,-0.263) -- (-0.438,0.387) -- cycle;
\path[fill=PTblue!14,draw=gray!65,line width=0.4pt] (0.979,0.671) -- (1.426,0.401) -- (1.459,-0.215) -- (0.923,-0.492) -- (0.508,-0.236) -- (0.466,0.287) -- cycle;
\path[fill=PTblue!14,draw=gray!65,line width=0.4pt] (1.459,-0.215) -- (1.426,0.401) -- (2.065,0.718) -- (2.308,0.584) -- (2.374,-0.215) -- (1.876,-0.443) -- cycle;
\path[fill=PTblue!14,draw=gray!65,line width=0.4pt] (2.374,-0.215) -- (2.308,0.584) -- (3.128,0.829) -- (3.139,0.783) -- (3.151,0.727) -- (3.161,0.672) -- (3.169,0.616) -- (3.176,0.560) -- (3.181,0.504) -- (3.184,0.447) -- (3.185,0.391) -- (3.184,0.335) -- (3.182,0.278) -- (3.178,0.222) -- (3.172,0.166) -- (3.165,0.111) -- (3.156,0.055) -- (3.145,0.000) -- (3.132,-0.055) -- (3.119,-0.109) -- (3.103,-0.163) -- (3.086,-0.216) -- (3.068,-0.268) -- (3.049,-0.320) -- (3.028,-0.372) -- (3.007,-0.423) -- (2.984,-0.473) -- (2.978,-0.485) -- cycle;
\path[fill=PTblue!14,draw=gray!65,line width=0.4pt] (-2.064,1.054) -- (-1.943,0.715) -- (-2.412,0.084) -- (-2.914,0.052) -- (-2.919,0.102) -- (-2.923,0.153) -- (-2.928,0.205) -- (-2.931,0.256) -- (-2.934,0.308) -- (-2.937,0.361) -- (-2.939,0.413) -- (-2.940,0.466) -- (-2.940,0.518) -- (-2.939,0.571) -- (-2.938,0.624) -- (-2.935,0.678) -- (-2.932,0.731) -- (-2.927,0.784) -- (-2.922,0.838) -- (-2.915,0.891) -- (-2.906,0.944) -- (-2.897,0.998) -- (-2.886,1.050) -- (-2.874,1.103) -- (-2.861,1.156) -- (-2.846,1.208) -- (-2.829,1.260) -- (-2.811,1.311) -- (-2.799,1.343) -- cycle;
\path[fill=PTblue!14,draw=gray!65,line width=0.4pt] (-1.943,0.715) -- (-2.064,1.054) -- (-1.562,1.624) -- (-0.976,1.473) -- (-1.068,0.638) -- (-1.344,0.486) -- cycle;
\path[fill=PTblue!14,draw=gray!65,line width=0.4pt] (-0.976,1.473) -- (-0.788,1.576) -- (-0.060,1.132) -- (-0.128,0.578) -- (-0.438,0.387) -- (-1.068,0.638) -- cycle;
\path[fill=PTblue!14,draw=gray!65,line width=0.4pt] (-0.060,1.132) -- (0.330,1.384) -- (0.939,1.035) -- (0.979,0.671) -- (0.466,0.287) -- (-0.128,0.578) -- cycle;
\path[fill=PTblue!14,draw=gray!65,line width=0.4pt] (0.979,0.671) -- (0.939,1.035) -- (1.404,1.510) -- (1.928,1.368) -- (2.065,0.718) -- (1.426,0.401) -- cycle;
\path[fill=PTblue!14,draw=gray!65,line width=0.4pt] (2.308,0.584) -- (2.065,0.718) -- (1.928,1.368) -- (2.489,1.883) -- (2.493,1.879) -- (2.532,1.839) -- (2.569,1.799) -- (2.607,1.758) -- (2.643,1.717) -- (2.679,1.674) -- (2.714,1.631) -- (2.748,1.587) -- (2.782,1.542) -- (2.814,1.496) -- (2.846,1.450) -- (2.876,1.403) -- (2.905,1.355) -- (2.934,1.306) -- (2.961,1.257) -- (2.986,1.206) -- (3.011,1.156) -- (3.033,1.104) -- (3.055,1.052) -- (3.075,0.999) -- (3.093,0.946) -- (3.110,0.892) -- (3.125,0.837) -- (3.128,0.829) -- cycle;
\path[fill=PTblue!14,draw=gray!65,line width=0.4pt] (-1.562,1.624) -- (-2.064,1.054) -- (-2.799,1.343) -- (-2.792,1.362) -- (-2.771,1.412) -- (-2.749,1.462) -- (-2.725,1.510) -- (-2.700,1.559) -- (-2.673,1.606) -- (-2.645,1.653) -- (-2.615,1.698) -- (-2.584,1.743) -- (-2.552,1.787) -- (-2.518,1.830) -- (-2.483,1.871) -- (-2.447,1.912) -- (-2.410,1.952) -- (-2.372,1.990) -- (-2.332,2.027) -- (-2.292,2.064) -- (-2.250,2.099) -- (-2.208,2.132) -- (-2.165,2.165) -- (-2.122,2.197) -- (-2.077,2.227) -- (-2.032,2.257) -- (-1.986,2.285) -- (-1.940,2.312) -- (-1.896,2.338) -- cycle;
\path[fill=PTblue!14,draw=gray!65,line width=0.4pt] (-0.621,2.019) -- (-0.788,1.576) -- (-0.976,1.473) -- (-1.562,1.624) -- (-1.896,2.338) -- (-1.894,2.339) -- (-1.847,2.364) -- (-1.800,2.388) -- (-1.752,2.412) -- (-1.704,2.434) -- (-1.656,2.456) -- (-1.608,2.477) -- (-1.560,2.497) -- (-1.512,2.516) -- (-1.463,2.535) -- (-1.415,2.553) -- (-1.366,2.570) -- (-1.318,2.587) -- (-1.270,2.603) -- (-1.221,2.619) -- (-1.173,2.634) -- (-1.124,2.649) -- (-1.076,2.663) -- (-1.028,2.677) -- (-0.979,2.691) -- (-0.931,2.704) -- (-0.906,2.711) -- cycle;
\path[fill=PTblue!14,draw=gray!65,line width=0.4pt] (-0.788,1.576) -- (-0.621,2.019) -- (0.253,2.140) -- (0.382,1.981) -- (0.330,1.384) -- (-0.060,1.132) -- cycle;
\path[fill=PTblue!14,draw=gray!65,line width=0.4pt] (0.382,1.981) -- (1.205,2.043) -- (1.404,1.510) -- (0.939,1.035) -- (0.330,1.384) -- cycle;
\path[fill=PTblue!14,draw=gray!65,line width=0.4pt] (1.928,1.368) -- (1.404,1.510) -- (1.205,2.043) -- (1.485,2.605) -- (1.500,2.598) -- (1.547,2.574) -- (1.593,2.550) -- (1.639,2.524) -- (1.685,2.498) -- (1.731,2.472) -- (1.776,2.444) -- (1.821,2.416) -- (1.865,2.388) -- (1.910,2.358) -- (1.954,2.328) -- (1.997,2.298) -- (2.041,2.266) -- (2.084,2.234) -- (2.126,2.202) -- (2.169,2.169) -- (2.211,2.135) -- (2.252,2.100) -- (2.293,2.065) -- (2.334,2.029) -- (2.375,1.993) -- (2.415,1.955) -- (2.454,1.917) -- (2.489,1.883) -- cycle;
\path[fill=PTblue!14,draw=gray!65,line width=0.4pt] (0.253,2.140) -- (-0.621,2.019) -- (-0.906,2.711) -- (-0.883,2.717) -- (-0.835,2.730) -- (-0.786,2.742) -- (-0.738,2.754) -- (-0.690,2.766) -- (-0.641,2.777) -- (-0.593,2.788) -- (-0.544,2.798) -- (-0.495,2.809) -- (-0.446,2.818) -- (-0.397,2.828) -- (-0.348,2.836) -- (-0.299,2.845) -- (-0.250,2.852) -- (-0.200,2.860) -- (-0.150,2.866) -- (-0.100,2.872) -- (-0.050,2.878) -- (0.000,2.882) -- (0.050,2.886) -- (0.101,2.889) -- (0.152,2.892) -- (0.202,2.893) -- (0.253,2.894) -- (0.285,2.894) -- cycle;
\path[fill=PTblue!14,draw=gray!65,line width=0.4pt] (1.205,2.043) -- (0.382,1.981) -- (0.253,2.140) -- (0.285,2.894) -- (0.304,2.894) -- (0.355,2.892) -- (0.406,2.890) -- (0.457,2.887) -- (0.508,2.883) -- (0.560,2.878) -- (0.611,2.873) -- (0.662,2.866) -- (0.712,2.858) -- (0.763,2.849) -- (0.814,2.839) -- (0.864,2.828) -- (0.915,2.816) -- (0.965,2.803) -- (1.015,2.789) -- (1.065,2.774) -- (1.114,2.758) -- (1.163,2.741) -- (1.212,2.723) -- (1.261,2.704) -- (1.309,2.685) -- (1.357,2.664) -- (1.405,2.643) -- (1.453,2.621) -- (1.485,2.605) -- cycle;
\path[draw,line width=0.8pt] (3.145,0.000) -- (3.156,0.055) -- (3.165,0.111) -- (3.172,0.166) -- (3.178,0.222) -- (3.182,0.278) -- (3.184,0.335) -- (3.185,0.391) -- (3.184,0.447) -- (3.181,0.504) -- (3.176,0.560) -- (3.169,0.616) -- (3.161,0.672) -- (3.151,0.727) -- (3.139,0.783) -- (3.125,0.837) -- (3.110,0.892) -- (3.093,0.946) -- (3.075,0.999) -- (3.055,1.052) -- (3.033,1.104) -- (3.011,1.156) -- (2.986,1.206) -- (2.961,1.257) -- (2.934,1.306) -- (2.905,1.355) -- (2.876,1.403) -- (2.846,1.450) -- (2.814,1.496) -- (2.782,1.542) -- (2.748,1.587) -- (2.714,1.631) -- (2.679,1.674) -- (2.643,1.717) -- (2.607,1.758) -- (2.569,1.799) -- (2.532,1.839) -- (2.493,1.879) -- (2.454,1.917) -- (2.415,1.955) -- (2.375,1.993) -- (2.334,2.029) -- (2.293,2.065) -- (2.252,2.100) -- (2.211,2.135) -- (2.169,2.169) -- (2.126,2.202) -- (2.084,2.234) -- (2.041,2.266) -- (1.997,2.298) -- (1.954,2.328) -- (1.910,2.358) -- (1.865,2.388) -- (1.821,2.416) -- (1.776,2.444) -- (1.731,2.472) -- (1.685,2.498) -- (1.639,2.524) -- (1.593,2.550) -- (1.547,2.574) -- (1.500,2.598) -- (1.453,2.621) -- (1.405,2.643) -- (1.357,2.664) -- (1.309,2.685) -- (1.261,2.704) -- (1.212,2.723) -- (1.163,2.741) -- (1.114,2.758) -- (1.065,2.774) -- (1.015,2.789) -- (0.965,2.803) -- (0.915,2.816) -- (0.864,2.828) -- (0.814,2.839) -- (0.763,2.849) -- (0.712,2.858) -- (0.662,2.866) -- (0.611,2.873) -- (0.560,2.878) -- (0.508,2.883) -- (0.457,2.887) -- (0.406,2.890) -- (0.355,2.892) -- (0.304,2.894) -- (0.253,2.894) -- (0.202,2.893) -- (0.152,2.892) -- (0.101,2.889) -- (0.050,2.886) -- (0.000,2.882) -- (-0.050,2.878) -- (-0.100,2.872) -- (-0.150,2.866) -- (-0.200,2.860) -- (-0.250,2.852) -- (-0.299,2.845) -- (-0.348,2.836) -- (-0.397,2.828) -- (-0.446,2.818) -- (-0.495,2.809) -- (-0.544,2.798) -- (-0.593,2.788) -- (-0.641,2.777) -- (-0.690,2.766) -- (-0.738,2.754) -- (-0.786,2.742) -- (-0.835,2.730) -- (-0.883,2.717) -- (-0.931,2.704) -- (-0.979,2.691) -- (-1.028,2.677) -- (-1.076,2.663) -- (-1.124,2.649) -- (-1.173,2.634) -- (-1.221,2.619) -- (-1.270,2.603) -- (-1.318,2.587) -- (-1.366,2.570) -- (-1.415,2.553) -- (-1.463,2.535) -- (-1.512,2.516) -- (-1.560,2.497) -- (-1.608,2.477) -- (-1.656,2.456) -- (-1.704,2.434) -- (-1.752,2.412) -- (-1.800,2.388) -- (-1.847,2.364) -- (-1.894,2.339) -- (-1.940,2.312) -- (-1.986,2.285) -- (-2.032,2.257) -- (-2.077,2.227) -- (-2.122,2.197) -- (-2.165,2.165) -- (-2.208,2.132) -- (-2.250,2.099) -- (-2.292,2.064) -- (-2.332,2.027) -- (-2.372,1.990) -- (-2.410,1.952) -- (-2.447,1.912) -- (-2.483,1.871) -- (-2.518,1.830) -- (-2.552,1.787) -- (-2.584,1.743) -- (-2.615,1.698) -- (-2.645,1.653) -- (-2.673,1.606) -- (-2.700,1.559) -- (-2.725,1.510) -- (-2.749,1.462) -- (-2.771,1.412) -- (-2.792,1.362) -- (-2.811,1.311) -- (-2.829,1.260) -- (-2.846,1.208) -- (-2.861,1.156) -- (-2.874,1.103) -- (-2.886,1.050) -- (-2.897,0.998) -- (-2.906,0.944) -- (-2.915,0.891) -- (-2.922,0.838) -- (-2.927,0.784) -- (-2.932,0.731) -- (-2.935,0.678) -- (-2.938,0.624) -- (-2.939,0.571) -- (-2.940,0.518) -- (-2.940,0.466) -- (-2.939,0.413) -- (-2.937,0.361) -- (-2.934,0.308) -- (-2.931,0.256) -- (-2.928,0.205) -- (-2.923,0.153) -- (-2.919,0.102) -- (-2.914,0.051) -- (-2.908,0.000) -- (-2.903,-0.051) -- (-2.897,-0.101) -- (-2.890,-0.151) -- (-2.884,-0.202) -- (-2.877,-0.252) -- (-2.870,-0.302) -- (-2.863,-0.352) -- (-2.855,-0.401) -- (-2.848,-0.451) -- (-2.840,-0.501) -- (-2.832,-0.551) -- (-2.824,-0.600) -- (-2.815,-0.650) -- (-2.807,-0.700) -- (-2.798,-0.750) -- (-2.788,-0.800) -- (-2.779,-0.850) -- (-2.768,-0.900) -- (-2.758,-0.950) -- (-2.747,-1.000) -- (-2.735,-1.050) -- (-2.723,-1.100) -- (-2.710,-1.150) -- (-2.697,-1.201) -- (-2.683,-1.251) -- (-2.668,-1.301) -- (-2.652,-1.351) -- (-2.636,-1.401) -- (-2.618,-1.451) -- (-2.600,-1.501) -- (-2.580,-1.550) -- (-2.560,-1.600) -- (-2.539,-1.649) -- (-2.516,-1.697) -- (-2.493,-1.745) -- (-2.468,-1.793) -- (-2.442,-1.840) -- (-2.416,-1.887) -- (-2.388,-1.933) -- (-2.358,-1.979) -- (-2.328,-2.024) -- (-2.297,-2.068) -- (-2.264,-2.111) -- (-2.230,-2.154) -- (-2.196,-2.196) -- (-2.160,-2.236) -- (-2.123,-2.276) -- (-2.085,-2.315) -- (-2.046,-2.354) -- (-2.006,-2.391) -- (-1.965,-2.427) -- (-1.923,-2.462) -- (-1.881,-2.496) -- (-1.837,-2.529) -- (-1.793,-2.561) -- (-1.748,-2.592) -- (-1.702,-2.621) -- (-1.656,-2.650) -- (-1.609,-2.678) -- (-1.561,-2.704) -- (-1.513,-2.730) -- (-1.465,-2.754) -- (-1.415,-2.778) -- (-1.366,-2.801) -- (-1.316,-2.822) -- (-1.266,-2.843) -- (-1.215,-2.862) -- (-1.164,-2.881) -- (-1.113,-2.898) -- (-1.061,-2.915) -- (-1.009,-2.931) -- (-0.957,-2.946) -- (-0.905,-2.960) -- (-0.853,-2.973) -- (-0.800,-2.985) -- (-0.747,-2.997) -- (-0.694,-3.008) -- (-0.641,-3.017) -- (-0.588,-3.026) -- (-0.535,-3.034) -- (-0.482,-3.041) -- (-0.428,-3.048) -- (-0.375,-3.053) -- (-0.321,-3.057) -- (-0.268,-3.061) -- (-0.214,-3.064) -- (-0.161,-3.065) -- (-0.107,-3.066) -- (-0.054,-3.066) -- (-0.000,-3.065) -- (0.053,-3.062) -- (0.107,-3.059) -- (0.160,-3.055) -- (0.213,-3.049) -- (0.266,-3.043) -- (0.319,-3.035) -- (0.372,-3.026) -- (0.424,-3.017) -- (0.476,-3.005) -- (0.528,-2.993) -- (0.579,-2.980) -- (0.630,-2.965) -- (0.681,-2.949) -- (0.731,-2.932) -- (0.781,-2.914) -- (0.830,-2.894) -- (0.879,-2.874) -- (0.927,-2.852) -- (0.974,-2.829) -- (1.021,-2.805) -- (1.067,-2.780) -- (1.112,-2.754) -- (1.157,-2.726) -- (1.201,-2.698) -- (1.244,-2.669) -- (1.287,-2.639) -- (1.329,-2.608) -- (1.369,-2.576) -- (1.410,-2.543) -- (1.449,-2.510) -- (1.487,-2.476) -- (1.525,-2.441) -- (1.562,-2.406) -- (1.598,-2.370) -- (1.634,-2.334) -- (1.669,-2.297) -- (1.703,-2.260) -- (1.736,-2.222) -- (1.769,-2.185) -- (1.802,-2.147) -- (1.833,-2.109) -- (1.865,-2.071) -- (1.895,-2.033) -- (1.926,-1.994) -- (1.956,-1.956) -- (1.986,-1.917) -- (2.015,-1.879) -- (2.044,-1.841) -- (2.073,-1.802) -- (2.102,-1.764) -- (2.131,-1.726) -- (2.160,-1.687) -- (2.188,-1.649) -- (2.217,-1.611) -- (2.246,-1.572) -- (2.274,-1.534) -- (2.303,-1.496) -- (2.332,-1.457) -- (2.361,-1.419) -- (2.391,-1.380) -- (2.420,-1.341) -- (2.449,-1.302) -- (2.479,-1.263) -- (2.509,-1.223) -- (2.538,-1.184) -- (2.568,-1.143) -- (2.598,-1.103) -- (2.627,-1.061) -- (2.657,-1.020) -- (2.686,-0.978) -- (2.716,-0.935) -- (2.745,-0.892) -- (2.773,-0.848) -- (2.802,-0.803) -- (2.830,-0.758) -- (2.857,-0.712) -- (2.884,-0.666) -- (2.910,-0.619) -- (2.936,-0.571) -- (2.960,-0.522) -- (2.984,-0.473) -- (3.007,-0.423) -- (3.028,-0.372) -- (3.049,-0.320) -- (3.068,-0.268) -- (3.086,-0.216) -- (3.103,-0.163) -- (3.119,-0.109) -- (3.132,-0.055) -- cycle;
\fill (0,0) circle (1.5pt) node[below left=-1pt] {$x_P$};
\node at (-2.75,2.75) {$P$};
\draw[->] (-4.0,1.4) node[left]{$1-v_P$} -- (-0.272,0.122);
\draw[->] (-4.0,-1.4) node[left]{$v_P/(N_P-1)$} -- (-2.088,-1.566);
\end{tikzpicture}
\caption{The central weights of Definition \ref{d:central}. The central child $P^\downarrow$, which contains $x_P$, receives weight $1-v_P$, and every other child receives weight $v_P/(N_P-1)$.} 
\label{fig:central}
\end{figure}

By \eqref{child-range}, $N_P\geq 2$ and $\log_{\xi_k}N_P\geq\log_{\xi_k}\#\Inner(P)\geq q-\log_{\xi_k}C_\XX\geq q-\log_{\xi_k}C_{\mathrm{ent}}\geq s_k$. Hence $v_P$ exists and is unique by Lemma \ref{l:central-entropy}. By \eqref{b-entropy}, the range \eqref{direct-range} contains $(0,3q/4]$ for every $k\geq1$. Figure \ref{fig:central} illustrates the central weights.

\begin{lemma}[existence]\label{l:central-existence}
The masses in Definition \ref{d:central} are the values of a Radon measure $\muse$ on $\XX$ with full support, which satisfies $\muse(\partial Q)=0$ for every cube $Q\in\Delta_k$ with $k\geq0$ and has target entropies $\mathbf s$.
\end{lemma}

\begin{proof}
Since $\XX$ is complete and Ahlfors regular, $\XX$ is proper. The inner balls of the cubes in $\Delta_k$ are pairwise disjoint and have radius $c_br_k$, and every cube in $\Delta_k$ has diameter at most $2C_br_k$. By \eqref{AR-background}, $(\Delta_k)_{k\geq0}$ is a sequence of nested partitions, and by construction the masses form a mass distribution for this sequence. Let $P\in\Delta_{k-1}$ with $k\geq 1$, and put $N:=N_P$ and $M:=\#\Inner(P)$. By Lemma \ref{l:vAR-count}, every $R\in\Inner(P)$ satisfies $\gap(R,\XX\setminus P)>9\side R$. Since $\Inner(P)$ is finite, $\gap(\bigcup\Inner(P),\XX\setminus P)>0$. Since $v_P\leq(N-1)/N$, every child of $P$ other than $P^\downarrow$ receives weight $v_P/(N-1)\leq 1/N$. The $N-M$ children in $\Outer(P)$ are of this kind, because $P^\downarrow\in\Inner(P)$. The weights $p_P(R)$, $R\in\Child(P)$, sum to $1$, and therefore
\[
 \sum_{R\in\Inner(P)}p_P(R)=1-\sum_{R\in\Outer(P)}p_P(R)\geq1-\frac{N-M}N=\frac MN\geq C_\XX^{-2},
\]
where the last inequality follows from \eqref{vin-out-count}. Corollary \ref{c:tame} with $\mathcal I(P)=\Inner(P)$ and $\eta=C_\XX^{-2}$ and Theorem \ref{t:existence} give a Radon measure $\muse$ that takes the prescribed values on the cubes in $\Delta_k$ with $k\geq0$ and satisfies $\muse(\partial Q)=0$ for these cubes. A cube $Q\in\Delta_k$ with $k<0$ is the union of the finitely many cubes in $\Delta_0$ contained in $Q$, and therefore $\muse(Q)=\nu(Q)$. Every weight $p_P(R)$ is positive. Hence every cube has positive and finite mass, and since every ball contains a cube, $\muse$ has full support. Finally, the left side of \eqref{target-entropy} at $P\in\Delta_{k-1}$ is $h_{\xi_k,N_P}(v_P)=s_k$.
\end{proof}

Next, we compute the local dimensions of measures with target entropies.

\begin{lemma}[entropy estimates]\label{l:information}
Let $N\geq1$ be an integer and let $p_1,\dots,p_N$ be positive numbers with $p_1+\dots+p_N=1$. Then
\[
 \sum_{i=1}^Np_i\log\frac1{p_i}\leq\log N\quad\text{and}\quad\sum_{i=1}^Np_i\left(\log\frac1{p_i}\right)^2\leq(\log N)^2+2\log N+2.
\]
\end{lemma}

\begin{proof}
By the concavity of the logarithm, $\sum_{i=1}^Np_i\log(1/p_i)\leq\log\sum_{i=1}^Np_i(1/p_i)=\log N$. For the second inequality, put $z_i:=\log(1/p_i)\geq0$, and for $u\geq0$, put $I(u):=\{i:z_i>u\}$. Every $i\in I(u)$ satisfies $p_i<e^{-u}$, and hence
\[
 \sum_{i\in I(u)}p_i\leq\min\{1,Ne^{-u}\}\quad\text{for all }u\geq0.
\]
Since $z_i^2=\int_0^{z_i}2u\,du$ for every $i$, interchanging the sum and the integral gives
\[
 \sum_{i=1}^Np_iz_i^2=\int_0^\infty2u\sum_{i\in I(u)}p_i\,du\leq\int_0^{\log N}2u\,du+N\int_{\log N}^\infty2ue^{-u}\,du=(\log N)^2+2\log N+2.\qedhere
\]
\end{proof}

\begin{lemma}[second moments]\label{l:second-moment}
Put
\[
 C_{\mathrm{var}}:=(q+1)^2+\frac{2(q+1)}{\log b}+\frac2{(\log b)^2}.
\]
If $\muse$ has target entropies $\mathbf s=(s_k)_{k\geq1}$, then $0\leq s_k\leq q+1$ for all $k\geq1$ and
\[
 \int_P(X_k-s_k)^2\,d\muse\leq C_{\mathrm{var}}\,\muse(P)\quad\text{for all }P\in\Delta_{k-1}\text{ and }k\geq1.
\]
\end{lemma}

\begin{proof}
Fix $k\geq1$ and $P\in\Delta_{k-1}$, put $N:=\#\Child(P)$, and put $p_R:=\muse(R)/\muse(P)$ for $R\in\Child(P)$. The numbers $p_R$ are positive and sum to $1$, and $X_k\log\xi_k=\log(1/p_R)$ on every $R\in\Child(P)$. By \eqref{target-entropy}, Lemma \ref{l:information}, \eqref{child-range}, and $C_\XX\leq b\leq\xi_k$,
\[
 s_k\log\xi_k=\sum_{R\in\Child(P)}p_R\log\frac1{p_R}\leq\log N\leq\log C_\XX+q\log\xi_k\leq (q+1)\log\xi_k.
\]
Hence $0\leq s_k\leq q+1$. By \eqref{target-integral}, $\int_P(X_k-s_k)^2\,d\muse=\int_PX_k^2\,d\muse-s_k^2\muse(P)$. Lemma \ref{l:information} and $\log N\leq(q+1)\log\xi_k$ give
\begin{align*}
 \int_P(X_k-s_k)^2\,d\muse\leq\int_PX_k^2\,d\muse&=\frac{\muse(P)}{(\log\xi_k)^2}\sum_{R\in\Child(P)}p_R\left(\log\frac1{p_R}\right)^2\\
 &\leq\left((q+1)^2+\frac{2(q+1)}{\log\xi_k}+\frac2{(\log\xi_k)^2}\right)\muse(P),
\end{align*}
which is at most $C_{\mathrm{var}}\,\muse(P)$.
\end{proof}

We use the strong law of large numbers in the following form.

\begin{lemma}[strong law of large numbers {\cite[Theorem~3.3.1]{StoutASC}}]\label{l:slln}
Let $(\mathscr F_k)_{k\geq0}$ be an increasing sequence of $\sigma$-algebras on a probability space with probability $\PP$ and expectation $\EE$. Let $(Y_k)_{k\geq1}$ be square-integrable random variables such that $Y_k$ is $\mathscr F_k$-measurable and $\EE(Y_k\mid\mathscr F_{k-1})=0$ for all $k\geq1$. If
\[
 \sum_{k=1}^\infty\frac{\EE(Y_k^2\mid\mathscr F_{k-1})}{k^2}<\infty\quad\text{almost surely},
\]
then $\lim_{n\to\infty}\frac1n\sum_{k=1}^nY_k=0$ almost surely.
\end{lemma}

\begin{theorem}[local dimensions]\label{t:dimension}
Let $\muse$ be a Radon measure on $\XX$ with target entropies $\mathbf s=(s_k)_{k\geq1}$. Then, at $\muse$-a.e.~$x\in\XX$,
\begin{equation}\label{endpoint-target-Cesaro}
 \liminf_{r\downarrow0}\frac{\log\muse(B(x,r))}{\log r}=\liminf_{n\to\infty}A_n(\mathbf s),\quad
 \limsup_{r\downarrow0}\frac{\log\muse(B(x,r))}{\log r}=\limsup_{n\to\infty}A_n(\mathbf s).
\end{equation}
\end{theorem}

\begin{proof}
Fix $Q_0\in\Delta_0$, normalize $\muse(Q_0)=1$, and use the notation introduced after Definition \ref{d:target}. Put $Y_k:=X_k-s_k$. Then $\EE(Y_k\mid\mathscr F_{k-1})=0$ by \eqref{target-integral}, and $\EE(Y_k^2\mid\mathscr F_{k-1})\leq C_{\mathrm{var}}$ by Lemma \ref{l:second-moment}, for all $k\geq1$. Lemma \ref{l:slln} therefore yields
\begin{equation}\label{endpoint-SLLN}
 \lim_{n\to\infty}\frac1n\sum_{k=1}^n(X_k-s_k)=0
 \quad\text{at }\muse\text{-a.e.~}x\in Q_0.
\end{equation}
By \eqref{endpoint-Xk} and \eqref{xi-def},
\begin{equation}\label{endpoint-weighted}
 \log\frac1{\muse(Q_n(x))}=\sum_{k=1}^nX_k(x)\log\xi_k
 \quad\text{and}\quad
 \log\frac1{r_n}=\sum_{k=1}^n\log\xi_k\quad\text{for all }n\geq 1.
\end{equation}
The sequence $(\log\xi_k)_{k\geq1}$ is nondecreasing. Put $L_n:=\sum_{k=1}^n\log(\xi_n/\xi_k)$. Then $L_n=0$ when $\xi_k=b$ for all $k$, and $L_n=o(n)$ by Lemma \ref{l:sparse-lag} for a system with sparse jumps. Moreover, $\log(1/r_n)=n\log\xi_n-L_n$. By Lemma \ref{l:second-moment}, $0\leq s_k\leq q+1$, and hence
\[
 0\leq\log\xi_n\sum_{k=1}^ns_k-\sum_{k=1}^ns_k\log\xi_k=\sum_{k=1}^ns_k\log\frac{\xi_n}{\xi_k}\leq (q+1)L_n.
\]
Fix $x\in Q_0$ at which \eqref{endpoint-SLLN} holds, and put $S_n:=\sum_{k=1}^nY_k(x)$. Then $S_n=o(n)$, and summation by parts gives
\[
 \sum_{k=1}^nY_k(x)\log\xi_k=S_n\log\xi_n-\sum_{k=1}^{n-1}S_k\log\frac{\xi_{k+1}}{\xi_k}.
\]
Let $\varepsilon>0$, and choose $k_\varepsilon$ such that $|S_k|\leq\varepsilon k$ for all $k\geq k_\varepsilon$. Since $(\log\xi_k)_{k\geq1}$ is positive and nondecreasing, for all $n>k_\varepsilon$,
\[
 \left|\sum_{k=1}^nY_k(x)\log\xi_k\right|\leq C_\varepsilon+\varepsilon n\log\xi_n+\varepsilon\sum_{k=1}^{n-1}k\log\frac{\xi_{k+1}}{\xi_k}\leq C_\varepsilon+2\varepsilon n\log\xi_n,
\]
where $C_\varepsilon$ is the sum of $|S_k|\log(\xi_{k+1}/\xi_k)$ over $k<k_\varepsilon$, which does not depend on $n$.
Dividing by $n\log\xi_n$, letting $n\to\infty$, and then letting $\varepsilon\downarrow0$, we obtain $\sum_{k=1}^n(X_k(x)-s_k)\log\xi_k=o(n\log\xi_n)$. By \eqref{endpoint-weighted},
\[
 \frac{\log\muse(Q_n(x))}{\log r_n}=\frac{\log\xi_n\sum_{k=1}^ns_k-\sum_{k=1}^ns_k\log(\xi_n/\xi_k)+\sum_{k=1}^n(X_k(x)-s_k)\log\xi_k}{n\log\xi_n-L_n}.
\]
Divide the numerator and the denominator by $n\log\xi_n$. Since $0\leq s_k\leq q+1$, $L_n=o(n)$, and $\log\xi_n\geq\log b$, we obtain
\begin{equation}\label{endpoint-cubical}
 \frac{\log\muse(Q_n(x))}{\log r_n}=A_n(\mathbf s)+o(1)
 \quad\text{at }\muse\text{-a.e.~}x\in Q_0.
\end{equation}

It remains to replace the cubes $Q_n(x)$ by balls. Put $\delta_k:=\tfrac12c_br_k$ and $\Lambda:=2C_b/c_b$. For every $k\geq0$ and $Q\in\Delta_k$, the closed ball $B(x_Q,\delta_k)$ lies in the inner ball $U_Q$, and $Q\subset B(x_Q,C_br_k)=B(x_Q,\Lambda\delta_k)$. Since the inner balls of the cubes in $\Delta_k$ are pairwise disjoint, $\Delta_k$ is a $(\delta_k,\Lambda)$-partition in the sense of Definition \ref{KRS-partition}. Moreover, $\delta_{k+1}\leq\delta_k/b$ and $\log(1/\delta_k)\geq\log(1/r_k)\geq k\log b$. Since $\xi_{k+1}\leq b(1+\log_2(k+1))$,
\[
 1\leq\frac{\log\delta_{k+1}}{\log\delta_k}=1+\frac{\log\xi_{k+1}}{\log(1/\delta_k)}\leq1+\frac{\log b+\log(1+\log_2(k+1))}{k\log b}\quad\text{for all }k\geq1,
\]
and the right side tends to $1$ as $k\to\infty$. Hence $(\Delta_k)_{k\geq0}$ is an admissible family of $\delta_k$-partitions in the sense of Definition \ref{KRS-partition-2}.

The space $\XX$ is doubling, because it is Ahlfors regular. By Lemma \ref{t:mf}, at $\muse$-a.e.~$x\in\XX$, the lower and upper local dimensions of $\muse$ at $x$ are the lower and upper limits of $\log\muse(Q_n(x))/\log\delta_n$. Since $\log\delta_n/\log r_n\to1$, these are also the lower and upper limits of $\log\muse(Q_n(x))/\log r_n$, which are the right sides of \eqref{endpoint-target-Cesaro} at $\muse$-a.e.~$x\in Q_0$ by \eqref{endpoint-cubical}. Thus, \eqref{endpoint-target-Cesaro} holds at $\muse$-a.e.~$x\in Q_0$.

Finally, the cubes in $\Delta_0$ have pairwise disjoint inner balls, and $\XX$ is separable, because it is proper. Hence $\Delta_0$ is countable. Since $Q_0\in\Delta_0$ was arbitrary and \eqref{endpoint-target-Cesaro} does not change when $\muse$ is multiplied by a positive constant, \eqref{endpoint-target-Cesaro} holds at $\muse$-a.e.~$x\in\XX$.
\end{proof}

Lemma \ref{l:central-existence} and Theorem \ref{t:dimension} produce, for every sequence of target entropies satisfying \eqref{direct-range}, a measure with full support whose lower and upper local dimensions are the lower and upper limits of $A_n(\mathbf s)$. We use this construction in \S\ref{s:zero}.

\section{Doubling measures with prescribed entropy}\label{s:doubling}

When $s_H>0$, Theorem \ref{t:main} asks for a doubling measure. As in \cite[Lemma 5.5]{BS4}, whose mechanism goes back to \cite[\S3]{KRS-cubes}, we prove doubling by comparing the masses of nearby cubes through their chains of ancestors. The comparison works when, at each level, siblings have comparable mass, and children near the boundaries of their parents receive a common weight, which may depend on the level, but not on the parent. The central weights of Definition \ref{d:central} lack these features in general. Weights with both features and with a prescribed entropy are available when every cube has a bounded number of children; see \cite[\S5]{BS4}. For a system with sparse jumps, however, the number of children of a cube is unbounded. To handle this, we insert intermediate levels between consecutive generations. Between two generations, we first spread mass over several intermediate levels using weights of quasi-Bernoulli type, called Type I weights, and then, if some intermediate levels remain, reach the next generation in a single step using weights comparable to the background measure, called Type II weights. Type I weights realize entropies up to a fixed distance below $q$, while Type II weights have entropy close to $q$ when they span a large change of scale. Mixing the two yields every entropy that we need.

We continue to follow Convention \ref{conv:standing}, with the system of cubes chosen in \S\ref{ss:intermediate}.

\subsection{Intermediate levels}\label{ss:intermediate}

Our next goal is to divide the passage from each generation of a system of $b$-adic cubes with sparse jumps to the next generation into steps with bounded ratios of side lengths. The hypotheses of Theorem \ref{vKRS-cubes} constrain only the ratios of consecutive side lengths, and the properties of a system of cubes persist when we keep only a subsequence of its levels. Thus, we may first choose a system of cubes with finer side lengths and then recover the original generations as a subsequence of its levels.

\begin{lemma}[subsequences of levels]\label{l:coarsen}
Let $(\Delta^*_\ell)_{\ell\in\ZZ}$ be a system of cubes with side lengths $(\tau_\ell)_{\ell\in\ZZ}$, and let $(\ell_k)_{k\in\ZZ}$ be a strictly increasing sequence of integers with $\ell_k=k$ for all $k\leq0$. Then $(\Delta^*_{\ell_k})_{k\in\ZZ}$, with the same centers, is a system of cubes with side lengths $(\tau_{\ell_k})_{k\in\ZZ}$. Moreover, if $k\in\ZZ$ and $Q\in\Delta^*_{\ell_{k-1}}$, then for every $\ell$ with $\ell_{k-1}\leq\ell\leq\ell_k$, the cubes in $\Delta^*_\ell$ contained in $Q$ form a partition of $Q$, which refines the partition at level $\ell-1$ when $\ell>\ell_{k-1}$. At $\ell=\ell_k$, this partition is the set of children of $Q$ in $(\Delta^*_{\ell_k})_{k\in\ZZ}$.
\end{lemma}

\begin{proof}
Put $r_k:=\tau_{\ell_k}$. Then $r_k=b^{-k}$ for all $k\leq0$. For $k\geq1$, we have $\ell_k-\ell_{k-1}\geq1$, and \eqref{geometric-sides} for $(\tau_\ell)_{\ell\in\ZZ}$ gives $r_k\leq b^{-(\ell_k-\ell_{k-1})}r_{k-1}\leq r_{k-1}/b$. Thus, $(r_k)_{k\in\ZZ}$ satisfies the hypotheses of Theorem \ref{vKRS-cubes} with the same $b$, and hence with the same constants $c_b$ and $C_b$. Nesting, roundness, and the origin property for $(\Delta^*_{\ell_k})_{k\in\ZZ}$ are special cases of the same properties for $(\Delta^*_\ell)_{\ell\in\ZZ}$. For inheritance, let $Q\in\Delta^*_{\ell_k}$, and apply inheritance in $(\Delta^*_\ell)_{\ell\in\ZZ}$ a total of $\ell_{k+1}-\ell_k$ times to obtain $R\in\Delta^*_{\ell_{k+1}}$ with $R\subset Q$ and $x_R=x_Q$.

Now let $Q\in\Delta^*_{\ell_{k-1}}$ and $\ell\geq\ell_{k-1}$. By nesting, every cube in $\Delta^*_\ell$ is either contained in $Q$ or disjoint from $Q$. Since $\Delta^*_\ell$ is a partition of $\XX$, the cubes in $\Delta^*_\ell$ contained in $Q$ form a partition of $Q$. By nesting again, every cube in $\Delta^*_\ell$ is contained in a cube in $\Delta^*_{\ell-1}$, which gives the refinement. The last assertion is the definition of children.
\end{proof}

\begin{definition}[intermediate levels for $b$-adic cubes with sparse jumps]\label{d:intermediate}
Let $(r_k)_{k\in\ZZ}$ and $(\xi_k)_{k\in\ZZ}$ be the side lengths and ratios of Definition \ref{sparse-cubes}. For every $k\in\ZZ$, put
\[
 V_k:=\lfloor\log_b\xi_k\rfloor\quad\text{and}\quad a_k:=\xi_k^{1/V_k}.
\]
Thus, $V_k\geq 1$ and $b\leq a_k<b^2$. Put $\ell_k:=k$ for all $k\leq 0$ and $\ell_k:=\sum_{j=1}^kV_j$ for all $k\geq 1$. For every $\ell\in\ZZ$, let $k$ be the unique integer with $\ell_{k-1}<\ell\leq\ell_k$, and put
\[
 \tau_\ell:=r_{k-1}\,a_k^{-(\ell-\ell_{k-1})}.
\]
Then $\tau_\ell=b^{-\ell}$ for all $\ell\leq 0$, $\tau_{\ell-1}=a_k\tau_\ell\geq b\,\tau_\ell$ for all $\ell\in\ZZ$, and $\tau_{\ell_k}=r_k$ for all $k\in\ZZ$. Let $(\Delta^*_\ell)_{\ell\in\ZZ}$ be a system of cubes with side lengths $(\tau_\ell)_{\ell\in\ZZ}$ given by Theorem \ref{vKRS-cubes}, and put $\Delta_k:=\Delta^*_{\ell_k}$ for all $k\in\ZZ$. By Lemma \ref{l:coarsen}, $(\Delta_k)_{k\in\ZZ}$ is a system of $b$-adic cubes with sparse jumps. We call $\Delta^*_\ell$ the $\ell$-th \emph{intermediate level}, and we call $(\Delta_k)_{k\in\ZZ}$, together with $(\Delta^*_\ell)_{\ell\in\ZZ}$, a system of $b$-adic cubes with sparse jumps and \emph{intermediate levels}.
\end{definition}

\begin{figure}[htb]
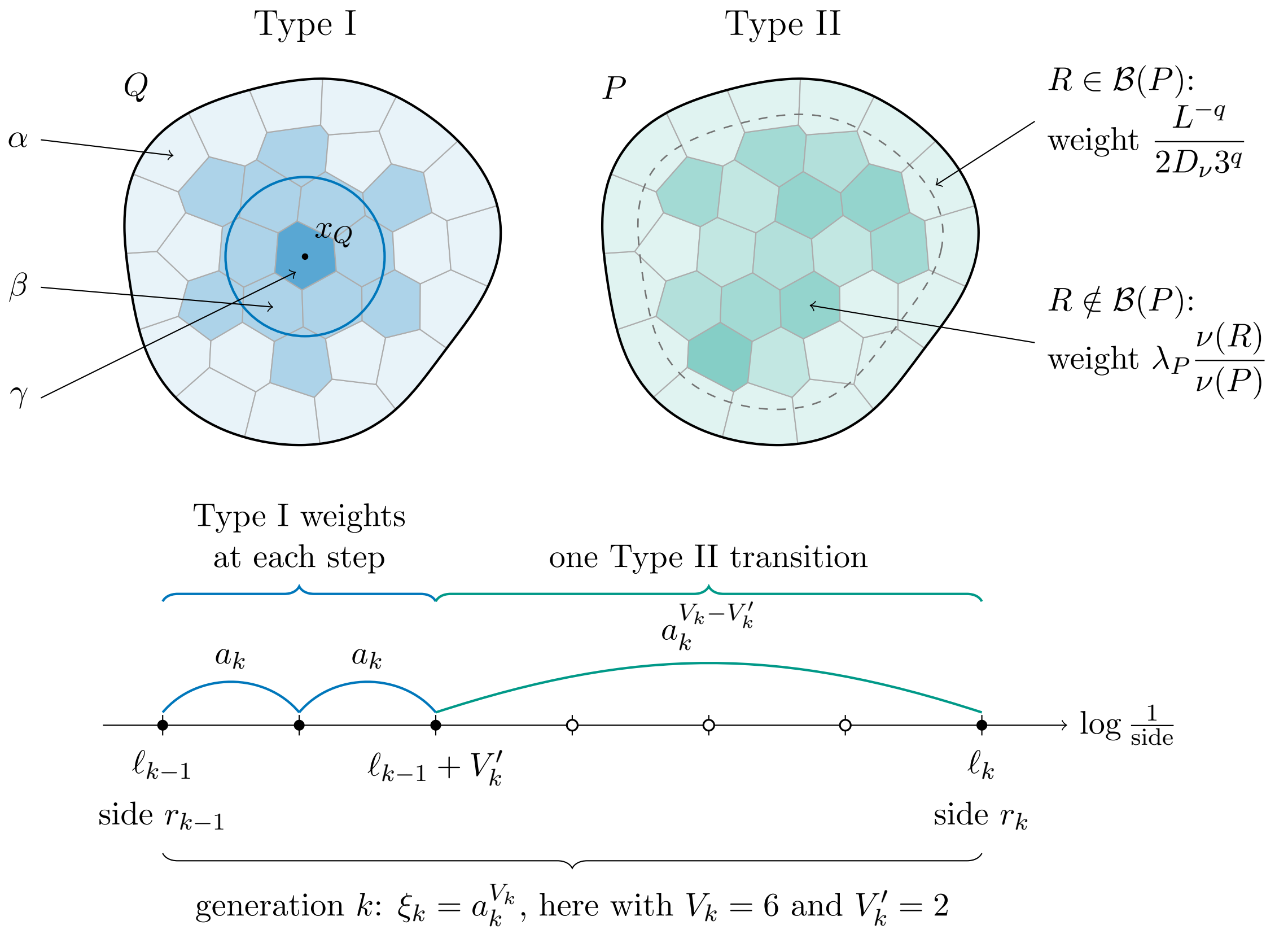

\centering

\caption{Top left: the Type I weights of Definition \ref{d:nested}; the circle is $U(x_Q,\tfrac12c_b\side Q)$. Top right: the Type II weights of Lemma \ref{l:flat}; $\mathcal B(P)$ consists of the cubes meeting the region outside the dashed contour. Shading and areas are schematic. Bottom: intermediate levels between the generations $k-1$ and $k$. The first $V_k'$ steps, each of ratio $a_k$, use Type I weights, and one Type II transition covers the rest. Filled dots mark the effective levels of \S\ref{ss:construction}, and hollow dots mark the levels skipped by the Type II transition.}\label{fig:weights}
\end{figure}

Throughout this section, the word \emph{generation} refers only to the partitions $\Delta_k=\Delta^*_{\ell_k}$, and the index $k$ always counts generations, while the index $\ell$ counts intermediate levels. We say that the intermediate levels $\ell$ with $\ell_{k-1}<\ell\leq\ell_k$ \emph{belong to} the $k$-th generation. Target entropies (Definition \ref{d:target}) are always computed between consecutive generations, with the ratios $\xi_k$.

For $Q\in\Delta_{k-1}$, we continue to write $\Child(Q)$ for the set of children of $Q$ in $\Delta_k$. For $Q\in\Delta^*_{\ell-1}$, we write $\Child^*(Q)$ for the set of cubes in $\Delta^*_\ell$ contained in $Q$, and we define $\Inner^*(Q)$ and $\Outer^*(Q)$ by \eqref{inner-def} and \eqref{outer-def} with $\Child^*(Q)$ in place of $\Child(Q)$. In this section, $Q^\downarrow$ always denotes the cube in $\Child^*(Q)$ containing $x_Q$. By Lemma \ref{l:coarsen}, for every $k\geq1$ and $Q\in\Delta_{k-1}$, the intermediate levels belonging to the $k$-th generation divide the passage from $Q$ to $\Child(Q)$ into $V_k$ steps, each with ratio of side lengths $a_k$. Lemma \ref{l:vAR-count} applies both to $(\Delta_k)_{k\in\ZZ}$, with ratios $\xi_k$, and to $(\Delta^*_\ell)_{\ell\in\ZZ}$, with ratio $a_k$ at every intermediate level belonging to the $k$-th generation; in the latter case, its conclusions hold with $\Child^*$, $\Inner^*$, and $\Outer^*$ in place of $\Child$, $\Inner$, and $\Outer$. If $\xi_k<b^2$, then $V_k=1$, and $\ell_k$ is the only intermediate level belonging to the $k$-th generation. On the $h$-th plateau (Remark \ref{r:plateau}), $\xi_k=bh$, so $V_k=1+\lfloor\log_bh\rfloor$. Thus, additional levels strictly between consecutive generations occur precisely on the plateaus with $h\geq b$.

When $\xi_k=b$ for all $k$, we use the same notation with
\[
 V_k=1,\qquad a_k=b,\qquad \ell_k=k,\qquad \tau_\ell=b^{-\ell},\qquad \Delta^*_\ell=\Delta_\ell,
\]
where $(\Delta_k)_{k\in\ZZ}$ is a system of $b$-adic cubes without jumps. In this case, every generation is an intermediate level, $\Child^*=\Child$, and only Type I weights occur in the construction of \S\ref{ss:construction}. Throughout the rest of this section, $(\Delta_k)_{k\in\ZZ}$ is a system of $b$-adic cubes without jumps, with this convention, or with sparse jumps and intermediate levels as in Definition \ref{d:intermediate}, and every result below applies in both cases. Every cube in an intermediate level $\ell\geq 0$ has diameter at most $2C_b<2.1$ and is not equal to $\XX$. The bottom panel of Figure \ref{fig:weights} shows the intermediate levels between two generations.

\subsection{Two types of weights}\label{ss:distributions}

We use two rules to assign weights to the descendants of a cube. Type I weights select an entropy at each bounded step between consecutive intermediate levels. Type II weights handle the remaining change of scale in a single step, with entropy close to the ambient value $q$ when that change is large. Both rules assign a common weight, which depends only on the levels involved and not on the parent, to all descendants sufficiently close to the boundary of the parent.

Put $N_*:=\lceil C_\XX b^{2q}\rceil$. By Lemma \ref{l:vAR-count}, if $\ell-1\geq0$ and $\ell$ belongs to the $k$-th generation, then $\#\Child^*(Q)\leq C_\XX a_k^q\leq N_*$ for every $Q\in\Delta^*_{\ell-1}$. The next lemma records the properties of the boundary weights that we use.

\begin{lemma}[{cf.~\cite[Lemma 5.3]{BS4}}]\label{l:boundary}
Let $N\geq 2$ be an integer and let $t_{\min}>0$. There are $\alpha>0$ and $0<\kappa_0\leq\alpha$, depending only on $b$, $N$, and $t_{\min}$, with the following property. Let $b\leq a\leq b^2$, let $n_{\mathrm{out}}\geq0$ and $n_{\mathrm{in}}\geq2$ be integers with $n_{\mathrm{out}}+n_{\mathrm{in}}\leq N$, and let $t_{\min}\leq t\leq\log_a n_{\mathrm{in}}$. There is a unique pair $(\beta,\gamma)$ with $0<\beta\leq\gamma$ such that
\begin{equation}\label{boundary-normalization}
 n_{\mathrm{out}}\alpha+(n_{\mathrm{in}}-1)\beta+\gamma=1,
\end{equation}
\begin{equation}\label{boundary-entropy}
 n_{\mathrm{out}}\alpha\log_a\frac1\alpha+(n_{\mathrm{in}}-1)\beta\log_a\frac1\beta+\gamma\log_a\frac1\gamma=t.
\end{equation}
For fixed $n_{\mathrm{out}}$ and $n_{\mathrm{in}}$, the pair $(\beta,\gamma)$ depends continuously on $(a,t)$. Moreover, $\beta\geq\kappa_0$ and $\gamma\geq 1/(2N)\geq\kappa_0$.
\end{lemma}

\begin{proof}
Choose $0<\alpha<1/(2N)$ with
\[
 N\alpha\big(\log(1/\alpha)+1\big)\leq\tfrac12t_{\min}\log b,
\]
which is possible because the left side tends to zero with $\alpha$. Fix $n_{\mathrm{out}}$ and $n_{\mathrm{in}}$, put $\omega:=1-n_{\mathrm{out}}\alpha$, and note that $\omega>1/2$. Writing $\gamma=\omega-(n_{\mathrm{in}}-1)\beta$, the condition \eqref{boundary-normalization} holds, and $0<\beta\leq\gamma$ is equivalent to $0<\beta\leq \omega/n_{\mathrm{in}}$. For $0\leq\beta\leq \omega/n_{\mathrm{in}}$, let $F_a(\beta)$ be the left side of \eqref{boundary-entropy}, with $0\log(1/0):=0$. Since $-(1-x)\log(1-x)\leq x$ for $0\leq x<1$, we have $\omega\log(1/\omega)\leq n_{\mathrm{out}}\alpha$, and therefore
\[
 F_a(0)=\frac{n_{\mathrm{out}}\alpha\log(1/\alpha)+\omega\log(1/\omega)}{\log a}\leq\frac{N\alpha(\log(1/\alpha)+1)}{\log b}\leq\frac{t_{\min}}2.
\]
At $\beta=\omega/n_{\mathrm{in}}$, also $\gamma=\omega/n_{\mathrm{in}}$, and
\[
 F_a(\omega/n_{\mathrm{in}})-\log_a n_{\mathrm{in}}=n_{\mathrm{out}}\alpha\log_a\frac1{n_{\mathrm{in}}\alpha}+\omega\log_a\frac1\omega \geq0,
\]
because $n_{\mathrm{in}}\alpha<1/2$ and $\omega\leq 1$. Moreover, $F_a'(\beta)=(n_{\mathrm{in}}-1)\log_a(\gamma/\beta)>0$ for $0<\beta<\omega/n_{\mathrm{in}}$. Hence $F_a$ is continuous and strictly increasing on $[0,\omega/n_{\mathrm{in}}]$ with $F_a(0)<t\leq F_a(\omega/n_{\mathrm{in}})$, and there is a unique $\beta\in(0,\omega/n_{\mathrm{in}}]$ with $F_a(\beta)=t$. By \eqref{boundary-normalization} and $\beta\leq\gamma$, we have $\gamma\geq \omega/n_{\mathrm{in}}>1/(2N)$.

Put $\mathcal K_{n_{\mathrm{out}},n_{\mathrm{in}}}:=\{(a,t):b\leq a\leq b^2,\ t_{\min}\leq t\leq\log_a n_{\mathrm{in}}\}$. Let $(a_i,t_i)\to(a,t)$ in $\mathcal K_{n_{\mathrm{out}},n_{\mathrm{in}}}$, and let $\beta_i$ and $\beta$ be the corresponding solutions. Every subsequential limit $\beta'$ of $(\beta_i)$ lies in $[0,\omega/n_{\mathrm{in}}]$ and satisfies $F_a(\beta')=t$, because $(a,\beta)\mapsto F_a(\beta)$ is continuous. As $F_a(0)<t$, uniqueness gives $\beta'=\beta$. Thus, $\beta$ and $\gamma$ depend continuously on $(a,t)$. Since $\mathcal K_{n_{\mathrm{out}},n_{\mathrm{in}}}$ is compact, $\beta$ has a positive minimum on $\mathcal K_{n_{\mathrm{out}},n_{\mathrm{in}}}$ whenever $\mathcal K_{n_{\mathrm{out}},n_{\mathrm{in}}}\neq\emptyset$. Let $\kappa_0$ be the least of $\alpha$, $1/(2N)$, and these minima over the finitely many pairs $(n_{\mathrm{out}},n_{\mathrm{in}})$ with $n_{\mathrm{out}}+n_{\mathrm{in}}\leq N$.
\end{proof}

\begin{definition}[Type I weights]\label{d:nested}
Let $0<t_{\min}\leq q/2$, and let $\alpha$ and $\kappa_0$ be given by Lemma \ref{l:boundary} with $N=N_*$. Let $k\geq 1$, let $\ell$ be an intermediate level belonging to the $k$-th generation, put $a:=a_k$, and put
\[
 t_0(a):=q-\log_aC_\XX.
\]
Let $Q\in\Delta^*_{\ell-1}$ and let $t_{\min}\leq t\leq t_0(a)$. The \emph{Type I weights} with target $t$ on $\Child^*(Q)$ assign the weight $\alpha$ to every $R\in\Outer^*(Q)$, the weight $\gamma$ to $Q^\downarrow$, and the weight $\beta$ to every other $R\in\Inner^*(Q)$, where $(\beta,\gamma)$ is given by Lemma \ref{l:boundary} with $n_{\mathrm{out}}=\#\Outer^*(Q)$ and $n_{\mathrm{in}}=\#\Inner^*(Q)$. By \eqref{boundary-normalization}, the Type I weights sum to $1$.
\end{definition}

Lemma \ref{l:boundary} applies in Definition \ref{d:nested}. Indeed, $n_{\mathrm{out}}+n_{\mathrm{in}}=\#\Child^*(Q)\leq C_\XX a^q\leq N_*$, and $n_{\mathrm{in}}\geq C_\XX^{-1}a^q$ by \eqref{vin-out-count}. Hence $t\leq t_0(a)\leq\log_a n_{\mathrm{in}}$, and $n_{\mathrm{in}}\geq2$ since $C_\XX^{-1}a^q\geq C_\XX^{-1}b^q>1$ by \eqref{b-entropy}. Thus, $\alpha$ is the common outer weight, $\beta$ is the weight of a noncentral inner child, and $\gamma$ is the central weight; see the top left panel of Figure \ref{fig:weights}. The outer weight $\alpha$ is fixed across parents, intermediate levels, and targets, whereas $\beta$ and $\gamma$ depend on the numbers of inner and outer children, on the ratio $a$, and on the target. Type I weights have the form of the weights in \cite[Definition 5.4]{BS4}, with $a$ in place of $b$.

\begin{remark}[properties of Type I weights]\label{r:nested}
Let $p$ be the Type I weights with target $t$ on $\Child^*(Q)$ as in Definition \ref{d:nested}. Then $\sum_{R\in\Child^*(Q)}p(R)\log_a(1/p(R))=t$, and $p$ depends continuously on $t$ by Lemma \ref{l:boundary}. Since the number of children is at most $N_*$, every weight of $p$ is at least $\kappa_0$, uniformly over parents, intermediate levels, and targets, where $\kappa_0$ depends only on $\XX$, $b$, and $t_{\min}$. The central child $Q^\downarrow$ receives weight at least $1/(2N_*)$. By Lemma \ref{l:vAR-count}, every $R\in\Child^*(Q)$ with $\gap(R,\XX\setminus Q)\leq 9\side R$ belongs to $\Outer^*(Q)$, and hence receives the weight $\alpha$, which depends neither on $Q$ nor on $\ell$ nor on $t$.
\end{remark}

\begin{lemma}[Type II weights]\label{l:flat}
Let $\ell<\ell'$ be integers with $\ell\geq 0$, let $P\in\Delta^*_\ell$, let $\mathcal D(P)$ be the set of cubes in $\Delta^*_{\ell'}$ contained in $P$, and put $L:=\tau_\ell/\tau_{\ell'}$. Put
\[
 \mathcal B(P):=\{R\in\mathcal D(P):\gap(R,\XX\setminus P)<\tau_\ell/5\}.
\]
There is a unique function $p_P:\mathcal D(P)\to(0,1]$ with $\sum_{R\in\mathcal D(P)}p_P(R)=1$ and a unique number $\lambda_P>0$ such that
\begin{equation}\label{flat-weights}
 p_P(R)=\frac{L^{-q}}{2D_\nu3^q}\text{ for all }R\in\mathcal B(P),\quad  p_P(R)=\lambda_P\,\frac{\nu(R)}{\nu(P)}\text{ for all }R\in\mathcal D(P)\setminus\mathcal B(P).
\end{equation}
Moreover,
\begin{equation}\label{flat-comparable}
 K^{-1}L^{-q}\leq p_P(R)\leq KL^{-q}\quad\text{for all }R\in\mathcal D(P),
\end{equation}
\begin{equation}\label{flat-interior}
 \sum_{R\in\mathcal D(P)\setminus\mathcal B(P)}p_P(R)\geq\frac12.
\end{equation}
We call $p_P$ the \emph{Type II weights} on $\mathcal D(P)$.
\end{lemma}

\begin{proof}
Since $L\geq b\geq47$, every $R\in\mathcal D(P)$ that meets $B(x_P,\tau_\ell/20)$ has $\diam R\leq 2.4\tau_\ell/L\leq 3\tau_\ell/20$ and lies in $B(x_P,\tau_\ell/5)$. Because $U(x_P,2\tau_\ell/5)\subset P$, such a cube satisfies $\gap(R,\XX\setminus P)\geq\tau_\ell/5$ and does not belong to $\mathcal B(P)$. By \eqref{AR-background},
\[
 1\geq\sum_{R\in\mathcal D(P)\setminus\mathcal B(P)}\frac{\nu(R)}{\nu(P)}\geq\frac{\nu(B(x_P,\tau_\ell/20))}{\nu(B(x_P,C_b\tau_\ell))}\geq D_\nu^{-1}24^{-q}.
\]
By \eqref{vAR-one-generation}, $\#\mathcal B(P)\leq D_\nu3^qL^q$, and the cubes in $\mathcal B(P)$ receive total weight at most $1/2$ in \eqref{flat-weights}. Hence the values of $p_P$ sum to $1$ for a unique $\lambda_P$, which satisfies $1/2\leq\lambda_P\leq D_\nu24^q$, and \eqref{flat-interior} holds. By \eqref{vAR-nu-ratio},
\[
 \frac{L^{-q}}{2D_\nu3^q}\leq p_P(R)\leq D_\nu^2\,72^qL^{-q}=KL^{-q}\quad\text{for all }R\in\mathcal D(P).
\]
It remains to check that $2D_\nu3^q\leq K$, which follows from $D_\nu4^q\geq2$. Fix distinct points $x,y\in\XX$ and put $d:=|x-y|$. The balls $B(x,d/3)$ and $B(y,d/3)$ are disjoint and contained in $B(x,4d/3)$, and $\nu(B(x,4d/3))\leq C_\nu D_\nu(4d/3)^q$ even when $4d/3>\diam\XX$, because then $B(x,4d/3)=B(x,\diam\XX)$. Hence $2C_\nu(d/3)^q\leq C_\nu D_\nu(4d/3)^q$.
\end{proof}

Type II weights are not in general equal, nor are they the normalized background masses $\nu(R)/\nu(P)$: in the boundary collar $\mathcal B(P)$ they are given by a common formula, and on the remaining descendants they are proportional to $\nu(R)/\nu(P)$. The property that we use is \eqref{flat-comparable}, which says that every weight is comparable to the background scale $L^{-q}$, with constants that do not depend on $L$. In particular, the entropy $\sum_{R\in\mathcal D(P)}p_P(R)\log(1/p_P(R))$ lies within $\log K$ of $q\log L$. Hence the entropy normalized by $\log L$ is close to $q$ when $L$ is large. Unlike for Type I weights, the lower bound on an individual Type II weight is relative to $L^{-q}$, and $L$ may be arbitrarily large. See Figure \ref{fig:weights}.

\subsection{Construction}\label{ss:construction}

We first describe the recursive mass assignment in terms of the parameters $V_k'$ and $t_P$. Lemma \ref{l:entropy-allocation} will choose these parameters so that the resulting measure $\muse$ has target entropies $\mathbf s=(s_k)_{k\geq1}$ in the original generations.

\begin{definition}[doubling measures with prescribed entropy]\label{d:endpoint-measure}
Let $0<t_{\min}\leq q/2$. Given integers $1\leq V_k'\leq V_k$ for all $k\geq 1$ and targets $t_P$ with $t_{\min}\leq t_P\leq t_0(a_k)$ for all $P\in\Delta_{k-1}$ and $k\geq 1$, we distribute mass as follows.
\begin{enumerate}
\item Declare $\muse(Q):=\nu(Q)$ for all $Q\in\Delta_k$ and $k\leq 0$.
\item Let $k\geq 1$ and $P\in\Delta_{k-1}$. For every intermediate level $\ell$ with $\ell_{k-1}<\ell\leq\ell_{k-1}+V_k'$ and every $Q\in\Delta^*_{\ell-1}$ contained in $P$, declare $\muse(R):=p(R)\muse(Q)$ for all $R\in\Child^*(Q)$, where $p$ denotes the Type I weights with target $t_P$.
\item If $V_k'<V_k$, then for every $Q\in\Delta^*_{\ell_{k-1}+V_k'}$ contained in a cube $P\in\Delta_{k-1}$, declare $\muse(R):=p_Q(R)\muse(Q)$ for all $R\in\Delta_k$ contained in $Q$, where $p_Q$ denotes the Type II weights of Lemma \ref{l:flat} with $\ell=\ell_{k-1}+V_k'$ and $\ell'=\ell_k$.
\end{enumerate}
By Lemma \ref{l:existence} below, these masses are the values of a Radon measure $\muse$ on $\XX$.
\end{definition}

Definition \ref{d:endpoint-measure} prescribes masses only at some of the intermediate levels, which we call effective. The set $\mathcal E$ of \emph{effective levels} of a weight assignment consists of all integers $\ell\leq0$ and, for every $k\geq1$, the levels $\ell$ with $\ell_{k-1}<\ell\leq\ell_{k-1}+V_k'$, together with $\ell_k$. When $V_k'<V_k$, the levels strictly between $\ell_{k-1}+V_k'$ and $\ell_k$ are not effective; these are the levels skipped by the Type II transition, shown as hollow dots in Figure \ref{fig:weights}. For $e\in\mathcal E$, let $e^-$ be the largest element of $\mathcal E$ less than $e$, i.e.~the index of the previous effective level. For $R\in\Delta^*_e$, we call the cube $R^\Uparrow\in\Delta^*_{e^-}$ containing $R$ the \emph{effective parent} of $R$. In this notation, Definition \ref{d:endpoint-measure} says that
\begin{equation}\label{effective-mass}
 \muse(R)=p(R)\,\muse(R^\Uparrow)\quad\text{for all }R\in\Delta^*_e\text{ and }e\in\mathcal E\text{ with }e\geq1.
\end{equation}
Here, if $k\geq1$ is the integer with $\ell_{k-1}<e\leq\ell_k$, then $p$ denotes the Type I weights on $\Child^*(R^\Uparrow)$ when $e\leq\ell_{k-1}+V_k'$, and the Type II weights on the cubes in $\Delta_k$ contained in $R^\Uparrow$ when $e=\ell_k$ and $V_k'<V_k$. Thus, the ratio $\side R^\Uparrow/\side R$ is $a_k$ at a Type I transition and $a_k^{V_k-V_k'}$ at a Type II transition, and the latter may be large.

\begin{lemma}[existence]\label{l:existence}
The masses in Definition \ref{d:endpoint-measure} are the values of a Radon measure $\muse$ on $\XX$, and $\muse(\partial Q)=0$ for every cube $Q\in\Delta^*_e$ with $e\in\mathcal E$.
\end{lemma}

\begin{proof}
Since $\XX$ is complete and Ahlfors regular, $\XX$ is proper. Let $0=e_0<e_1<e_2<\cdots$ enumerate the nonnegative effective levels. The inner balls of the cubes in $\Delta^*_{e_i}$ are pairwise disjoint, have radius $c_b\tau_{e_i}$, and every cube has diameter at most $2C_b\tau_{e_i}$. By \eqref{AR-background}, $(\Delta^*_{e_i})_{i\geq0}$ is a sequence of nested partitions. By construction, the masses form a mass distribution for this sequence. We verify the hypothesis of Corollary \ref{c:tame} with $\eta=1/(2N_*)$.
Let $P\in\Delta^*_{e_i}$. Suppose first that the masses of the cubes in $\Delta^*_{e_{i+1}}$ contained in $P$ are given by Type I weights. Put $\mathcal I(P):=\{P^\downarrow\}$. By Lemma \ref{l:vAR-count}, $\gap(P^\downarrow,\XX\setminus P)>9\side P^\downarrow$, and by Remark \ref{r:nested}, $\muse(P^\downarrow)\geq\muse(P)/(2N_*)\geq\eta\,\muse(P)$. Suppose instead that the masses are given by Type II weights. Put $\mathcal I(P):=\mathcal D(P)\setminus\mathcal B(P)$. By the definition of $\mathcal B(P)$ and \eqref{flat-interior}, $\gap(\bigcup\mathcal I(P),\XX\setminus P)\geq\tau_{e_i}/5$ and $\sum_{R\in\mathcal I(P)}\muse(R)\geq\muse(P)/2$. Corollary \ref{c:tame} and Theorem \ref{t:existence} complete the proof for the cubes in $\Delta^*_e$ with $e\geq0$. For $e<0$, the mass $\nu(Q)$ of a cube $Q\in\Delta^*_e$ is the sum of masses of the cubes in $\Delta_0$ contained in $Q$, and $\partial Q$ is contained in a union of the boundaries of the finitely many cubes in $\Delta_0$ contained in $Q$.
\end{proof}

It remains to choose the parameters so that every parent in $\Delta_{k-1}$ has the same prescribed entropy between generations $k-1$ and $k$. We use the number of Type I steps to select an interval of entropies attainable for every parent, and then adjust the target $t_P$ within each parent to obtain the desired value.

\begin{lemma}[entropy allocation]\label{l:entropy-allocation}
Let $k\geq1$ and let $0<t_{\min}\leq q/2$. For every $s$ with
\begin{equation}\label{endpoint-entropy-range}
 t_{\min}\leq s\leq q-\log_{\xi_k}C_{\mathrm{ent}},
\end{equation}
there is an integer $1\leq V_k'\leq V_k$ and, for every $P\in\Delta_{k-1}$, a target $t_P\in[t_{\min},t_0(a_k)]$ such that the masses of the cubes in $\Delta_k$ given by Definition \ref{d:endpoint-measure} satisfy
\begin{equation}\label{endpoint-entropy-exact}
 \sum_{R\in\Child(P)}\frac{\muse(R)}{\muse(P)}\log_{\xi_k}\frac{\muse(P)}{\muse(R)}=s\quad\text{for all }P\in\Delta_{k-1}.
\end{equation}
The interval \eqref{endpoint-entropy-range} is nonempty.
\end{lemma}

\begin{proof}
Put $a:=a_k$, $V:=V_k$, and $t_0:=t_0(a)$. For $P\in\Delta_{k-1}$, let $H(P)$ be the left side of \eqref{endpoint-entropy-exact} multiplied by $\log\xi_k$. Suppose that $1\leq V'\leq V$ and that all Type I weights below $P$ have target $t$. Grouping the cubes in $\Child(P)$ according to their ancestors in the intermediate levels gives
\[
 H(P)=V't\log a+\sum_{Q}\frac{\muse(Q)}{\muse(P)}H_{\mathrm{II}}(Q),
\]
where $Q$ ranges over the cubes in $\Delta^*_{\ell_{k-1}+V'}$ contained in $P$, and the sum is absent when $V'=V$. Here $p_Q$ denotes the Type II weights of Lemma \ref{l:flat} on the set $\mathcal D(Q)$ of cubes in $\Delta_k$ contained in $Q$, and
\[
 H_{\mathrm{II}}(Q):=\sum_{R\in\mathcal D(Q)}p_Q(R)\log\frac1{p_Q(R)}.
\] By \eqref{flat-comparable} with $L=a^{V-V'}$, each $H_{\mathrm{II}}(Q)$ lies in \[[q(V-V')\log a-\log K,\,q(V-V')\log a+\log K].\] By Remark \ref{r:nested}, $H(P)$ depends continuously on $t$. Put $I(V):=[Vt_{\min}\log a,\,Vt_0\log a]$, and for $1\leq V'<V$, put
\[
 I(V'):=\big[(V't_{\min}+(V-V')q)\log a+\log K,\,(V't_0+(V-V')q)\log a-\log K\big].
\]
For $V'<V$, the margins of $\log K$ allow for the variation in the Type II contribution. At $t=t_{\min}$, the value of $H(P)$ is at most the left endpoint of $I(V')$, and at $t=t_0$, it is at least the right endpoint. Thus, by continuity, every $\eta\in I(V')$ is attained for every parent $P$, using a target $t\in[t_{\min},t_0]$ that may depend on $P$. The same conclusion holds for $I(V)$, where no Type II contribution is present.

Since $\log_aC_\XX\leq\log_bC_{\mathrm{ent}}\leq q/4$ by \eqref{b-entropy} and $t_{\min}\leq q/2$, we have $t_0-t_{\min}\geq q/4$. Hence
\[
 V'(t_0-t_{\min})\log a\geq\frac q4\log b\geq\log C_{\mathrm{ent}}=\log C_\XX+2\log K\quad\text{for all }V'\geq 1.
\]
This inequality implies that each $I(V')$ is nonempty and $\min I(V')\leq\max I(V'+1)$ for all $1\leq V'<V$. Moreover, $\min I(V)=t_{\min}\log\xi_k$ and $\max I(1)\geq q\log\xi_k-\log C_\XX-\log K$. Thus, the union of intervals $I(V')$ with $1\leq V'\leq V$ contains $[t_{\min}\log\xi_k,\,q\log\xi_k-\log C_{\mathrm{ent}}]$. Given $s$ in \eqref{endpoint-entropy-range}, let $V_k'$ be the largest $V'$ with $s\log\xi_k\in I(V')$, and for every $P\in\Delta_{k-1}$ choose $t_P$ with $H(P)=s\log\xi_k$. This gives \eqref{endpoint-entropy-exact}.

The interval \eqref{endpoint-entropy-range} is nonempty, because $\log_{\xi_k}C_{\mathrm{ent}}\leq\log_bC_{\mathrm{ent}}\leq q/4<q-t_{\min}$.
\end{proof}

The integer $V_k'$ in Lemma \ref{l:entropy-allocation} is the same for all parents in the $k$-th generation, although the target $t_P$ may depend on $P$. Hence the effective levels, and with them the common boundary factors in the proof of Lemma \ref{l:endpoint-doubling} below, are the same for neighboring parents, while the targets $t_P$ adjust the entropy at every parent to exactly $s_k$. Choosing $s_k$ in \eqref{endpoint-entropy-range} for every $k\geq1$ and applying Lemma \ref{l:entropy-allocation} at every generation, we obtain a measure $\muse$ given by Definition \ref{d:endpoint-measure} with target entropies $\mathbf s=(s_k)_{k\geq1}$ in the sense of Definition \ref{d:target}. Every cube has positive and finite mass, and $\muse$ has full support. Its local dimensions are given by Theorem \ref{t:dimension}.

\subsection{Doubling and dimensions}\label{ss:doubling-dimensions}

\begin{figure}[htb]
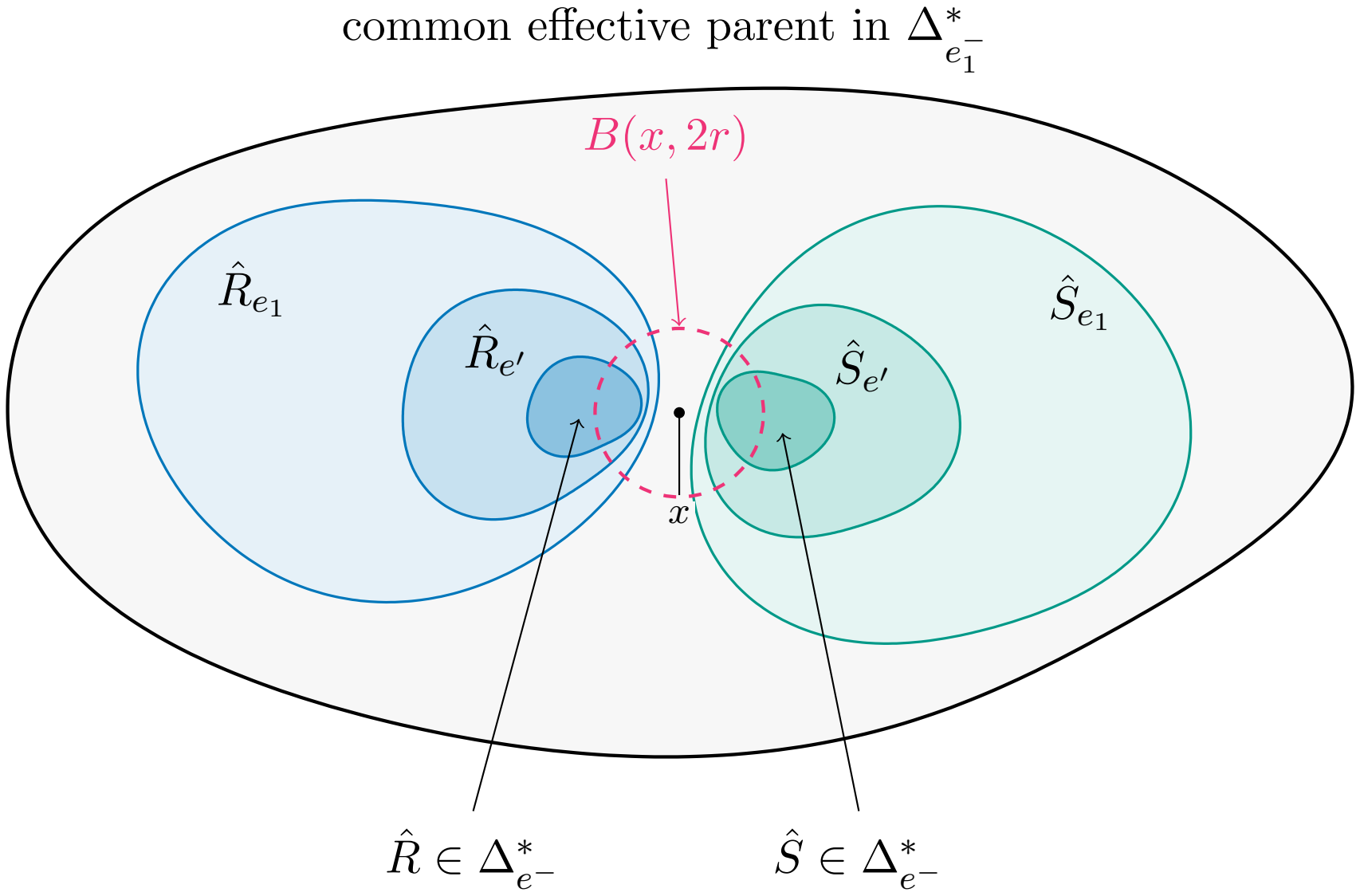

\centering

\caption{The comparison step in the proof of Lemma \ref{l:endpoint-doubling} when $e_1\geq1$. The ancestors of $\hat R$ (blue) and $\hat S$ (teal) first differ at the level $e_1$, inside a common effective parent. At every effective level $e'$ with $e_1<e'\leq e^-$, both ancestors lie near the boundaries of their parents and receive the common boundary factor $w_{e'}$. Other children and the exact ratios of side lengths are omitted.}\label{fig:doubling}
\end{figure}

\begin{lemma}[doubling]\label{l:endpoint-doubling}
Every measure given by Definition \ref{d:endpoint-measure} satisfies \eqref{doubling} with $D_{\mathsf{db}}:=D_\nu^2\,30^q\,K_*^2$, where $K_*:=\max\{\kappa_0^{-1},K^2\}$.
\end{lemma}

\begin{proof}
Let $\mathcal E$, $e^-$, and $R^\Uparrow$ be as defined before Lemma \ref{l:existence}. By \eqref{effective-mass}, the ratio $\muse(R)/\muse(R^\Uparrow)$ is the weight assigned to $R$ by Type I or Type II weights on the descendants of $R^\Uparrow$ in $\Delta^*_e$. We record two properties, valid for all $e\in\mathcal E$ with $e\geq 1$.
\begin{enumerate}[label=(\alph*)]
\item There is $m_e>0$ such that $m_e\leq\muse(R)/\muse(R^\Uparrow)\leq K_*m_e$ for all $R\in\Delta^*_e$.
\item There is $w_e>0$, the \emph{common boundary factor}, such that $\muse(R)=w_e\,\muse(R^\Uparrow)$ for all $R\in\Delta^*_e$ with $\gap(R,\XX\setminus R^\Uparrow)\leq 9\side R$.
\end{enumerate}
For Type I weights, (a) holds with $m_e=\kappa_0$, and (b) holds with $w_e=\alpha$, by Remark \ref{r:nested}. For Type II weights, put $L:=\side R^\Uparrow/\side R\geq 47$. Then (a) holds with $m_e=K^{-1}L^{-q}$ by \eqref{flat-comparable}. Moreover, (b) holds with $w_e=L^{-q}/(2D_\nu 3^q)$, because $9\side R<\side R^\Uparrow/5$ places $R$ in $\mathcal B(R^\Uparrow)$.

Let $x\in\XX$ and $r>0$. If $B(x,r)=\XX$, there is nothing to prove. Assume that $B(x,r)\neq\XX$. In particular, $r<\diam\XX$. If a cube $S\in\Delta^*_\ell$ is contained in $B(x,r)$, then $U_S\neq\XX$, and \eqref{AR-background} gives $C_\nu(2\tau_\ell/5)^q\leq C_\nu D_\nu r^q$. Hence there is a least $e\in\mathcal E$ such that some cube in $\Delta^*_e$ is contained in $B(x,r)$. Put $\sigma:=\tau_e$,
\[
 \mathcal I:=\{S\in\Delta^*_e:S\subset B(x,r)\},\quad\text{and}\quad\mathcal O:=\{R\in\Delta^*_e:R\cap B(x,2r)\neq\emptyset\}.
\]
No cube in $\mathcal O$ equals $\XX$, because $\mathcal I$ is nonempty and $B(x,r)\neq\XX$.

The inner balls of the cubes in $\mathcal O$ are pairwise disjoint and lie in $B(x,2r+2.4\sigma)$. Hence $\#\mathcal O\leq D_\nu((2r+2.4\sigma)/(0.4\sigma))^q$. If $r\geq 4.8\sigma$, then every cube in $\Delta^*_e$ meeting $B(x,r/2)$ belongs to $\mathcal I$, and $\#\mathcal I\geq D_\nu^{-1}(r/(2.4\sigma))^q$. If $r<4.8\sigma$, then $\#\mathcal I\geq 1$. In both cases,
\begin{equation}\label{IO-count}
 \#\mathcal O\leq D_\nu^2\,30^q\,\#\mathcal I.
\end{equation}

We claim that $\muse(R)\leq K_*^2\,\muse(S)$ for all $R\in\mathcal O$ and $S\in\mathcal I$. If $e\leq 0$, then $\muse=\nu$ on $\Delta^*_e$, and \eqref{AR-background} gives $\muse(R)\leq D_\nu3^q\muse(S)\leq K_*\,\muse(S)$. Suppose that $e\geq 1$, and put $\hat R:=R^\Uparrow$ and $\hat S:=S^\Uparrow$. By the minimality of $e$, the cube in $\Delta^*_{e^-}$ containing $x$ is not contained in $B(x,r)$, and hence $r<2C_b\tau_{e^-}$. Since $\hat R$ and $\hat S$ both meet $B(x,2r)$,
\[
 \gap(\hat R,\hat S)\leq 4r<8C_b\tau_{e^-}<9\tau_{e^-}.
\]
For every $e'\in\mathcal E$ with $e'\leq e^-$, let $\hat R_{e'}$ and $\hat S_{e'}$ be the cubes in $\Delta^*_{e'}$ containing $\hat R$ and $\hat S$, respectively. Suppose that $\hat R\neq\hat S$ (see Figure \ref{fig:doubling}). By Remark \ref{r:exhaustive}, there is a least $e_1\in\mathcal E$ with $\hat R_{e_1}\neq\hat S_{e_1}$. Put $e_2:=\max\{e_1,0\}$. If $e_1\geq 1$, then $\hat R_{e_1}$ and $\hat S_{e_1}$ are contained in the same cube in $\Delta^*_{e_1^-}$, and (a) gives $\muse(\hat R_{e_1})\leq K_*\,\muse(\hat S_{e_1})$. If $e_1\leq 0$, then \eqref{AR-background} gives $\muse(\hat R_0)\leq D_\nu 3^q\muse(\hat S_0)$. In both cases, $\muse(\hat R_{e_2})\leq K_*\,\muse(\hat S_{e_2})$. Now let $e'\in\mathcal E$ with $e_2<e'\leq e^-$. Then $\hat S$ is disjoint from $\hat R_{(e')^-}$, and
\[
 \gap(\hat R_{e'},\XX\setminus\hat R_{(e')^-})\leq\gap(\hat R,\hat S)<9\tau_{e^-}\leq 9\side\hat R_{e'}.
\]
By (b), $\muse(\hat R_{e'})=w_{e'}\,\muse(\hat R_{(e')^-})$, and the same identity holds with $\hat S$ in place of $\hat R$. It follows that $\muse(\hat R)/\muse(\hat R_{e_2})=\muse(\hat S)/\muse(\hat S_{e_2})$, and hence $\muse(\hat R)\leq K_*\,\muse(\hat S)$. This inequality also holds when $\hat R=\hat S$. By (a),
\[
 \muse(R)\leq K_*m_e\,\muse(\hat R)\leq K_*^2m_e\,\muse(\hat S)\leq K_*^2\,\muse(S),
\] as claimed.

Putting it all together, by \eqref{IO-count} and the previous displayed equation,
\begin{align*}
 \muse(B(x,2r))&\leq\sum_{R\in\mathcal O}\muse(R)\leq\#\mathcal O\,K_*^2\min_{S\in\mathcal I}\muse(S)\\
 &\leq D_\nu^2\,30^q\,K_*^2\sum_{S\in\mathcal I}\muse(S)\leq D_\nu^2\,30^q\,K_*^2\,\muse(B(x,r)).\qedhere
\end{align*}
\end{proof}

We can now prove an existence theorem for doubling measures with prescribed lower and upper Hausdorff and packing dimensions. Theorem \ref{t:prescribe} extends \cite[Theorem 5.7]{BS4} and establishes \cite[Conjecture 5.8]{BS4} in the affirmative.

\begin{theorem}\label{t:prescribe}
Let $\XX$ be a complete Ahlfors $q$-regular metric space.  Suppose
\[
 0<d_*\leq d^*\leq q,\qquad
 0<D_*\leq D^*\leq q,\qquad
 d_*\leq D_*,\qquad d^*\leq D^*.
\]
Then there is a doubling measure $\mu$ on $\XX$ with full support such that
\begin{equation}\label{haus-pack-dims}
 \underline{\dim}_H\mu=d_*,
 \qquad \overline{\dim}_H\mu=d^*,
 \qquad \underline{\dim}_P\mu=D_*,
 \qquad \overline{\dim}_P\mu=D^*.
\end{equation}
The constant $D_{\mathsf{db}}$ in \eqref{doubling} for $\mu$ depends only on $\XX$, $d_*$, and $d^*$. In particular, if $d_*=d^*=s_H$ and $D_*=D^*=s_P$, then $\mu$ has lower local dimension $s_H$ and upper local dimension $s_P$ at $\mu$-a.e.~point, and $D_{\mathsf{db}}$ depends only on $\XX$ and $s_H$.
\end{theorem}

\begin{proof}
We first construct, for given $0<s_H\leq s_P\leq q$, a doubling measure with full support whose lower and upper local dimensions are $s_H$ and $s_P$ at almost every point. Put $t_{\min}:=\min\{s_H,q/2\}$, choose $b\geq\max\{47,C_\XX\}$, depending only on $\XX$ and $s_H$, such that $\log_bC_{\mathrm{ent}}\leq(q-s_H)/4$ if $s_H<q$, and let $(\Delta_k)_{k\in\ZZ}$ be a system of $b$-adic cubes with sparse jumps and intermediate levels. If $s_H<q$, alternate increasingly long blocks of generations with values $s_k=s_H$ and $s_k=\min\{s_P,\,q-2\log_{\xi_k}C_{\mathrm{ent}}\}$. If $s_H=q$, put $s_k=q-2\log_{\xi_k}C_{\mathrm{ent}}$ for all $k\geq1$. Every value of $s_k$ lies in \eqref{endpoint-entropy-range}. Since $\log_{\xi_k}C_{\mathrm{ent}}\to0$, the blocks may be chosen so that $A_n(\mathbf s)$ has lower limit $s_H$ and upper limit $s_P$. By Lemma \ref{l:entropy-allocation}, there is a measure $\muse$ given by Definition \ref{d:endpoint-measure} with target entropies $\mathbf{s}=(s_k)_{k\geq1}$. Every cube has positive mass, and $\muse$ has full support. By Lemma \ref{l:endpoint-doubling} and Theorem \ref{t:dimension}, $\muse$ satisfies \eqref{doubling} with $D_{\mathsf{db}}$ depending only on $\XX$ and $s_H$, and $\muse$ has lower local dimension $s_H$ and upper local dimension $s_P$ at $\muse$-a.e.~point.

If $d_*=d^*=s_H$ and $D_*=D^*=s_P$, then the measure $\mu:=\muse$ of the first paragraph satisfies \eqref{haus-pack-dims} by the theorem of Tamashiro recalled in \S\ref{ss:dim}.

In general, let $\mu_1$ and $\mu_2$ be the measures of the first paragraph with $(s_H,s_P)=(d_*,D_*)$ and $(s_H,s_P)=(d^*,D^*)$, respectively. Then the measure $\mu:=\mu_1+\mu_2$ has full support and satisfies \eqref{doubling} with $D_{\mathsf{db}}$ equal to the larger of the constants for $\mu_1$ and $\mu_2$. The lower dimensions of a finite sum of measures are the minima of the corresponding lower dimensions of the summands, while the upper dimensions are the maxima. This gives \eqref{haus-pack-dims}.
\end{proof}

\section{\texorpdfstring{Upper endpoints: $s_P=m$ or $s_P=q$}{Upper endpoints: s\textunderscore P = m or s\textunderscore P = q}}\label{s:upper}

In this section, we treat the upper endpoints under the standing assumption $s_H>0$. We first dispose of the cases covered by Corollary \ref{c:bk}. Suppose
\[
 m-1<s_P<m.
\]
Theorem \ref{t:prescribe} with $d_*=d^*=s_H$ and $D_*=D^*=s_P$ supplies a full-support doubling measure with lower local dimension $s_H$ and upper local dimension $s_P$.  Corollary \ref{c:bk} makes this measure $m$-rectifiable, and Remark \ref{r:pure-unrectifiable} makes it purely $(m-1)$-unrectifiable, because $s_P>m-1$. This includes the endpoint $s_P=q<m$.
Consequently, the only upper-endpoint assertion requiring an additional rectifiability argument is $s_P=m\leq q$. Both cases $m<q$ and $m=q$ rely on the rectifiability criterion in \S\ref{ss:square-criterion}. When $m<q$, a fixed ratio $b$ already allows target entropies up to $m$, and we use $b$-adic cubes without jumps. When $m=q$, the target entropies must approach $q$, and we use $b$-adic cubes with sparse jumps; see Challenge 1 in the introduction.

\subsection{A criterion for rectifiability}\label{ss:square-criterion}

On a block of generations on which the ratio $\xi_k$ is constant, we retain the descendants whose masses are not much smaller than the target entropies predict. The second-moment estimate shows that these descendants carry most of the mass, while the lower bound on each retained mass controls their number.

\begin{lemma}[good descendants]\label{l:entropy-good-descendants}
Let $\muse$ be a Radon measure on $\XX$ with target entropies $(s_k)_{k\geq1}$, and let $C_{\mathrm{var}}$ be the constant in Lemma \ref{l:second-moment}. Let $k\geq0$, let $Q\in\Delta_k$, let $n\geq1$ be an integer, and suppose that $\xi_{k+1}=\xi_{k+2}=\dots=\xi_{k+n}$. For $\delta>0$, put
\[
 \mathcal G_Q(n,\delta):=
 \left\{R\in\Child^n(Q):
 \log_{\xi_{k+1}}\frac{\muse(Q)}{\muse(R)}
 \leq\sum_{i=1}^ns_{k+i}+\delta n\right\}.
\]
Then
\begin{equation}\label{good-descendant-mass}
 \muse\left(\bigcup\mathcal G_Q(n,\delta)\right)
 \geq\left(1-\frac{C_{\mathrm{var}}}{\delta^2n}\right)\muse(Q),
\end{equation}
\begin{equation}\label{good-descendant-count}
 \#\mathcal G_Q(n,\delta)
 \leq \xi_{k+1}^{\sum_{i=1}^ns_{k+i}+\delta n}.
\end{equation}
\end{lemma}

\begin{proof}
Let $X_k$ be given by \eqref{endpoint-Xk}. For $x\in Q$, we have $\log_{\xi_{k+1}}(\muse(Q)/\muse(Q_{k+n}(x)))=\sum_{i=1}^nX_{k+i}(x)$. Put $Y_i:=X_{k+i}-s_{k+i}$ for $1\leq i\leq n$. Let $1\leq i<j\leq n$. The function $Y_i$ is constant on every cube in $\Delta_{k+i}$, and hence on every $P\in\Child^{j-1}(Q)$. By \eqref{target-integral}, $\int_PY_j\,d\muse=0$ for every such $P$. Summing over $P\in\Child^{j-1}(Q)$ gives $\int_QY_iY_j\,d\muse=0$. Similarly, Lemma \ref{l:second-moment} and a sum over $P\in\Child^{i-1}(Q)$ give $\int_QY_i^2\,d\muse\leq C_{\mathrm{var}}\muse(Q)$. Expanding the square, we obtain
\[
 \int_Q\left(\sum_{i=1}^nY_i\right)^2d\muse
 =\sum_{i=1}^n\int_QY_i^2\,d\muse\leq C_{\mathrm{var}}n\,\muse(Q).
\]
Chebyshev's inequality gives
\[
 \muse\left\{x\in Q:\sum_{i=1}^nY_i(x)>\delta n\right\}
 \leq\frac{C_{\mathrm{var}}}{\delta^2n}\muse(Q).
\]
Every other point of $Q$ lies in a cube in $\mathcal G_Q(n,\delta)$, which proves \eqref{good-descendant-mass}. If $R\in\mathcal G_Q(n,\delta)$, then
\[
 \muse(R)\geq \muse(Q)\,\xi_{k+1}^{-\sum_{i=1}^ns_{k+i}-\delta n}.
\]
Summing this inequality over the pairwise disjoint good descendants gives \eqref{good-descendant-count}.
\end{proof}

\begin{proposition}[rectifiability criterion]\label{p:square}
Let $m\leq q$ be a positive integer, and let $\muse$ be a Radon measure on $\XX$ with target entropies $(s_k)_{k\geq1}$. If there is $c_0>0$ such that
\begin{equation}\label{square-hypothesis}
 m-s_k\geq\frac{c_0}{\log_2k}\quad\text{for all sufficiently large }k,
\end{equation}
then $\muse$ is $m$-rectifiable.
\end{proposition}

Proposition \ref{p:square} applies to both systems of cubes in Convention \ref{conv:standing}. In both systems, $\xi_k\leq k$ for all sufficiently large $k$. Hence \eqref{square-hypothesis} holds, in particular, when $(m-s_k)\log_2\xi_k$ is bounded below by a positive constant. The proof uses \eqref{square-hypothesis} only through the deficit accumulated over long blocks of consecutive generations; we have not tried to optimize the rate $1/\log_2k$.

\begin{figure}[htb]
\centering
\begin{tikzpicture}[font=\small]
\draw[->] (-0.2,0)--(13.4,0) node[right]{$k$};
\foreach \a/\b/\c in {0.4/1.4/15,1.4/2.4/30,2.4/3.4/15,3.4/4.4/30,4.4/6.4/15,6.4/8.4/30,8.4/10.4/15,10.4/12.4/30} {\fill[PTblue!\c] (\a,-0.25) rectangle (\b,0.25); \draw (\a,-0.25) rectangle (\b,0.25);}
\draw[PTmagenta,very thick] (4.4,-0.45)--(4.4,0.45);
\node[PTmagenta,above] at (4.4,0.45) {jump};
\node at (0.9,0) {$J_j$}; \node at (5.4,0) {$J_{j+4}$};
\draw[decorate,decoration={brace,amplitude=4pt}] (1.45,0.35)--(2.35,0.35) node[midway,above=5pt]{$2^{\lceil 2h/3\rceil}$};
\draw[decorate,decoration={brace,amplitude=4pt}] (6.45,0.35)--(8.35,0.35) node[midway,above=5pt]{$2^{\lceil 2(h+1)/3\rceil}$};
\draw[decorate,decoration={brace,amplitude=5pt,mirror}] (0.4,-0.6)--(4.35,-0.6) node[midway,below=6pt]{$H_h=\{2^{h-1}<k\leq2^h\}$};
\draw[decorate,decoration={brace,amplitude=5pt,mirror}] (4.45,-0.6)--(12.4,-0.6) node[midway,below=6pt]{$H_{h+1}=\{2^{h}<k\leq2^{h+1}\}$};
\end{tikzpicture}
\caption{The blocks of generations in the proof of Proposition \ref{p:square}, drawn to scale for $h=9$, so that $H_h$ and $H_{h+1}$ contain four blocks each. For a system with sparse jumps, $H_h$ is the $h$-th plateau, and the vertical line marks the jump at the generation $2^h+1$.}\label{fig:blocks}
\end{figure}
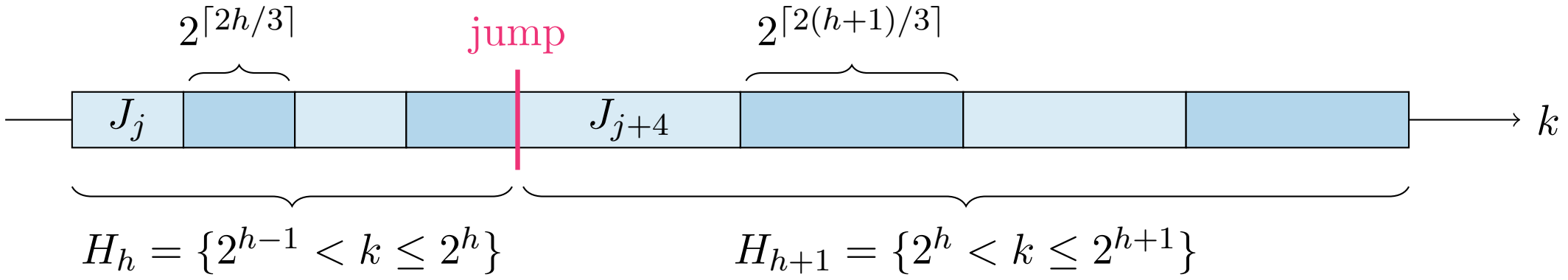

\begin{proof}
Fix $0<\theta<1$. We find an integer $h_0\geq3$ such that every cube $Q\in\Delta_{2^{h_0-1}}$ contains a set of measure at least $\theta\muse(Q)$ that is contained in a Lipschitz image of a compact subset of $[0,1]^m$. This set will be the set of leaves of a tree that retains good descendants at the end of every block in a partition of the generations after $2^{h_0-1}$. We choose the blocks long enough that the total mass discarded over all blocks is small. At the same time, each block must be short relative to the preceding generations, so that the entropy deficit accumulated earlier dominates the number of descendants created in the current block. The following choices satisfy both requirements.

For every integer $h\geq1$, put $H_h:=\{k\in\ZZ:2^{h-1}<k\leq2^h\}$. If $h\geq3$, then $\lceil2h/3\rceil\leq h-1$, and we divide $H_h$ into $2^{h-1-\lceil2h/3\rceil}\leq2^{h/3-1}$ blocks of $2^{\lceil2h/3\rceil}$ consecutive generations. Let $J_0,J_1,J_2,\dots$ be the blocks contained in the sets $H_h$ with $h\geq h_0$, listed in increasing order; see Figure \ref{fig:blocks}. For every block $J$, let $h_J$ be the integer with $J\subset H_{h_J}$, and put $k_J:=\min J-1$. Thus, $k_{J_0}=2^{h_0-1}$ and $k_{J_{j+1}}=\max J_j$ for all $j\geq0$. For a system with sparse jumps, $H_h$ is the $h$-th plateau when $h\geq2$ (Remark \ref{r:plateau}). Hence $\xi_k$ is constant on every block $J$ for both systems of cubes, with a value $\xi_J\leq bh_J$. Since the blocks preceding $J_j$ lie in the sets $H_{h'}$ with $h'\leq h_{J_j}$,
\begin{equation}\label{block-count}
 j<\sum_{h'=1}^{h_{J_j}}2^{h'/3-1}<3\cdot2^{h_{J_j}/3}\quad\text{for all }j\geq0.
\end{equation}

We next show that, if $h_0$ is sufficiently large, then
\begin{equation}\label{square-deficit}
 \sum_{k\in J}(m-s_k)\log_2\xi_k-2^{-h_J/12}|J|\log_2\xi_J\geq\frac{c_0|J|}{2h_J}\quad\text{for every block }J.
\end{equation}
Let $J$ be a block, put $h:=h_J$, and let $k\in J$. Then $\log_2k\leq h$. If $h_0$ is sufficiently large, then \eqref{square-hypothesis} and $\log_2\xi_k\geq\log_2b\geq1$ give $(m-s_k)\log_2\xi_k\geq c_0/h$. Since $\log_2\xi_J\leq\log_2(bh)$, it follows that \eqref{square-deficit} holds whenever $2^{-h/12}\log_2(bh)\leq c_0/(2h)$, which is true for all sufficiently large $h$.

To build the tree, let $J$ be a block and let $P\in\Delta_{k_J}$. Since $\xi_k=\xi_J$ for all $k\in J$, Lemma \ref{l:entropy-good-descendants} applies with $k=k_J$, $n=|J|$, and $\delta=2^{-h_J/12}$; we write $\mathcal G_J(P)$ for the resulting set of good descendants of $P$ in $\Delta_{\max J}$. Since $|J|\geq2^{2h_J/3}$, \eqref{good-descendant-mass} and \eqref{good-descendant-count} give
\begin{equation}\label{block-mass-count}
\begin{gathered}
 \muse\Big(\bigcup\mathcal G_J(P)\Big)\geq\big(1-C_{\mathrm{var}}2^{-h_J/2}\big)\muse(P)\quad\text{and}\\
 \log_2\#\mathcal G_J(P)\leq\sum_{k\in J}s_k\log_2\xi_k+2^{-h_J/12}|J|\log_2\xi_J.
\end{gathered}
\end{equation}
Every set $H_h$ contains at most $2^{h/3-1}$ blocks. Hence
\[
 \sum_{j=0}^\infty C_{\mathrm{var}}2^{-h_{J_j}/2}\leq\frac{C_{\mathrm{var}}}2\sum_{h=h_0}^\infty2^{-h/6},
\]
and we choose $h_0$ sufficiently large that the right side is at most $1-\theta$ and \eqref{square-deficit} holds. Fix $Q\in\Delta_{k_{J_0}}$ and build a tree $\mathcal T$ by putting $\mathcal T_0=\{Q\}$ and, recursively,
\[
 \Child_{\mathcal T}(P):=\mathcal G_{J_j}(P)\quad\text{for all }P\in\mathcal T_j\text{ and }j\geq 0.
\]
Thus, $\mathcal T_j\subset\Delta_{k_{J_j}}$ for all $j\geq0$. By the choice of $h_0$, we have $C_{\mathrm{var}}2^{-h_J/2}<1$ for every block $J$. Hence every cube in $\mathcal T$ has a child in $\mathcal T$, by \eqref{block-mass-count}. Since the cubes in $\mathcal T_j$ are pairwise disjoint and contained in $Q$, \eqref{block-mass-count} also gives
\[
 \muse\Big(\bigcup\mathcal T_j\setminus\bigcup\mathcal T_{j+1}\Big)\leq C_{\mathrm{var}}2^{-h_{J_j}/2}\muse(Q)\quad\text{for all }j\geq0.
\]
Summing over $j$ and using continuity from above, we obtain
\begin{equation}\label{square-leaf-mass}
 \muse(\leaves(\mathcal T))\geq\theta\muse(Q).
\end{equation}

It remains to verify \eqref{root-condition}. Let $N_j$ and $D_j$ be as in Theorem \ref{t:pack}, and put $R_j:=\xi_{J_j}^{|J_j|}$, which is the ratio of the side lengths of the cubes in $\mathcal T_j$ and $\mathcal T_{j+1}$. Then $D_j\leq2C_br_{k_{J_j}}=2C_br_{k_{J_0}}\prod_{i<j}R_i^{-1}$. Since $\lceil x\rceil\leq 2x$ for all $x\geq 1$,
\begin{equation}\label{square-log}
 \left(\prod_{i=0}^j\lceil N_i^{1/m}\rceil\right)D_j\leq 2C_b\,r_{k_{J_0}}\cdot 2N_j^{1/m}\prod_{i<j}\frac{2N_i^{1/m}}{R_i}\quad\text{for all }j\geq0.
\end{equation}
By \eqref{block-mass-count} and \eqref{square-deficit},
\begin{align*}
 \log_2\frac{2N_i^{1/m}}{R_i}&\leq 1-\frac1m\left(\sum_{k\in J_i}(m-s_k)\log_2\xi_k-2^{-h_{J_i}/12}|J_i|\log_2\xi_{J_i}\right)\\
 &\leq 1-\frac{c_0|J_i|}{2m\,h_{J_i}}\quad\text{for all }i\geq0.
\end{align*}
By \eqref{vAR-one-generation}, $N_j\leq D_\nu3^qR_j^q$. Hence $\log_2(2N_j^{1/m})\leq 1+m^{-1}\log_2(D_\nu3^q)+(q/m)|J_j|\log_2(bh_{J_j})$.

Now let $j\geq0$ with $h:=h_{J_j}>h_0$. Every block in $H_{h-1}$ precedes $J_j$, and these blocks have total length $2^{h-2}$. Keeping only these blocks in the sum over $i<j$ and using \eqref{block-count} and $|J_j|\leq2^{(2h+2)/3}$, we find that the base-$2$ logarithm of the left side of \eqref{square-log} is at most
\[
 C+3\cdot2^{h/3}+\frac qm\,2^{(2h+2)/3}\log_2(bh)-\frac{c_0}{8m}\cdot\frac{2^h}h,
\]
where $C:=\log_2(2C_br_{k_{J_0}})+1+m^{-1}\log_2(D_\nu3^q)$ does not depend on $j$. This is at most $-h$ when $h$ is sufficiently large. Since every $H_h$ contains at most $2^{h/3-1}$ blocks, the series \eqref{root-condition} converges.

Since $\XX$ is proper, the closure of $Q$ is compact. By Theorem \ref{t:pack}, $\leaves(\mathcal T)$ is contained in a Lipschitz image of a compact subset of $[0,1]^m$. Now let $Q_0\in\Delta_0$. The cube $Q_0$ is the union of finitely many cubes $Q\in\Delta_{k_{J_0}}$, and by \eqref{square-leaf-mass}, the union $F_\theta$ of the corresponding Lipschitz images satisfies $\muse(Q_0\setminus F_\theta)\leq(1-\theta)\muse(Q_0)$. Applying this with $\theta=1-1/i$ for all $i\geq 2$, we cover $\muse$-almost all of $Q_0$ by countably many Lipschitz images of compact subsets of $[0,1]^m$. Since $\Delta_0$ is countable, $\muse$ is $m$-rectifiable.
\end{proof}

\subsection{\texorpdfstring{The case $s_P=m<q$}{The case s\textunderscore P = m < q}}\label{ss:upper-fixed}

Put $t_{\min}:=s_H/2$. Choose target entropies $(s_k)_{k\geq1}\subset[t_{\min},m)$ with
\begin{equation}\label{upper-fixed-targets}
 \liminf_{n\to\infty}A_n(\mathbf s)=s_H,
 \qquad
 \limsup_{n\to\infty}A_n(\mathbf s)=m,\text{ and}
\end{equation}
\begin{equation}\label{upper-fixed-deficit}
 m-s_k\geq \frac1{\log_2(k+4)}\quad\text{for all sufficiently large }k.
\end{equation}
If $s_H=m$, take $s_k=m-1/\log_2(k+4)$ for all $k\geq1$. Otherwise, alternate blocks on which $s_k=s_H$ with blocks on which $s_k=m-1/\log_2(k+4)$, choosing the $j$-th pair of blocks long enough that the values of $A_n(\mathbf s)$ at the right endpoints of its two blocks are within $1/j$ of $s_H$ and $m$, respectively.

Every target lies in $[t_{\min},m)$, because $m-1/\log_2(k+4)\geq m/2\geq s_H/2$ for all $k\geq0$. To choose $b$ independently of $m$, note that $q-m\geq q+1-\lceil q\rceil>0$, since $m$ is an integer and $m<q$. Choose $b\geq\max\{47,C_\XX\}$, depending only on $\XX$, such that $\log_bC_{\mathrm{ent}}\leq\min\{q/4,\,q+1-\lceil q\rceil\}$, where $C_{\mathrm{ent}}$ is the constant in \eqref{K-and-C2}. Then $t_{\min}\leq s_k\leq m\leq q-\log_bC_{\mathrm{ent}}$.

Let $(\Delta_k)_{k\in\ZZ}$ be a system of $b$-adic cubes without jumps, with the convention of \S\ref{ss:intermediate}, and let $\muse$ be given by Definition \ref{d:endpoint-measure} with target entropies $(s_k)_{k\geq1}$. Such a measure exists by Lemma \ref{l:entropy-allocation}. Every cube has positive mass, and $\muse$ has full support. By Lemma \ref{l:endpoint-doubling} and Theorem \ref{t:dimension}, $\muse$ satisfies \eqref{doubling} with $D_{\mathsf{db}}$ depending only on $\XX$ and $s_H$, and $\muse$ has lower local dimension $s_H$ and upper local dimension $m$ at $\muse$-a.e.~point.

Since the upper local dimension of $\muse$ is $m>m-1$ at $\muse$-a.e.~point, Remark \ref{r:pure-unrectifiable} shows that $\muse$ is purely $(m-1)$-unrectifiable. Since $\log_2(k+4)\leq2\log_2k$ for all $k\geq3$, \eqref{upper-fixed-deficit} implies \eqref{square-hypothesis} with $c_0=1/2$. By Proposition \ref{p:square}, $\muse$ is $m$-rectifiable.

\subsection{\texorpdfstring{The case $s_P=m=q$}{The case s\textunderscore P = m = q}}\label{ss:upper-q}

Let $\muse$ be the measure constructed in the first paragraph of the proof of Theorem \ref{t:prescribe} with $s_P=q$. Thus, $\muse$ has full support, satisfies \eqref{doubling} with $D_{\mathsf{db}}$ depending only on $\XX$ and $s_H$, and has lower local dimension $s_H$ and upper local dimension $q$ at $\muse$-a.e.~point. The underlying system of $b$-adic cubes has sparse jumps, $\log_bC_{\mathrm{ent}}\leq(q-s_H)/4$ if $s_H<q$, and every target entropy $s_k$ equals either $q-2\log_{\xi_k}C_{\mathrm{ent}}$ or, when $s_H<q$, $s_H$.

Since $C_{\mathrm{ent}}>1$, each kind of generation that occurs satisfies a lower bound: generations with $s_k=s_H<q$ satisfy $(q-s_k)\log_2\xi_k\geq(q-s_H)\log_2b\geq 4\log_2C_{\mathrm{ent}}$, and generations with $s_k=q-2\log_{\xi_k}C_{\mathrm{ent}}$ satisfy $(q-s_k)\log_2\xi_k=2\log_2C_{\mathrm{ent}}$. Thus,
\begin{equation}\label{upper-q-weighted-deficit}
 (q-s_k)\log_2\xi_k\geq 2\log_2C_{\mathrm{ent}}\quad\text{for all }k\geq 1.
\end{equation}
Since $\xi_k\leq k$ for all sufficiently large $k$, \eqref{upper-q-weighted-deficit} gives $q-s_k\geq2\log_2C_{\mathrm{ent}}/\log_2k$ for all sufficiently large $k$, which is \eqref{square-hypothesis} with $m=q$ and $c_0=2\log_2C_{\mathrm{ent}}$. By Proposition \ref{p:square}, $\muse$ is $q$-rectifiable. Since the upper local dimension of $\muse$ is $q=m>m-1$ at $\muse$-a.e.~point, Remark \ref{r:pure-unrectifiable} shows that $\muse$ is purely $(m-1)$-unrectifiable. This completes the proof of the upper-endpoint assertions of Theorem \ref{t:main} when $s_H>0$.

\section{\texorpdfstring{Lower endpoint: $s_P=m-1$}{Lower endpoint: s\textunderscore P = m - 1}}\label{s:lower}

We continue under the standing assumption $s_H>0$. Since $0<s_H\leq s_P=m-1$, we have $m\geq2$. We use the measures of Definition \ref{d:endpoint-measure} on a system of $b$-adic cubes without jumps, with the convention of \S\ref{ss:intermediate}. In this case, $V_k=V_k'=1$ for all $k\geq1$, every generation is an intermediate level, and the masses of the children of every cube are given by Type I weights. Moreover, once \eqref{pack-dim-s} holds with $s_P=m-1<m$, Corollary \ref{c:bk} shows that the measure is $m$-rectifiable. Thus, the remaining task is to prove that the measure is purely $(m-1)$-unrectifiable.

We choose target entropies whose averages $A_n(\mathbf s)$ have upper limit $m-1$, but exceed $m-1$ by at least $n^{-1/4}$ along a sequence of generations $n$. The cumulative excess $n(A_n(\mathbf s)-(m-1))\geq n^{3/4}$ is then large enough to outweigh the fluctuations of the measure and to make the lower $(m-1)$-density vanish.

\begin{lemma}[target entropies at the lower endpoint]\label{l:lower-targets}
Let $m\geq2$ be an integer with $m-1<q$, and let $0<a\leq m-1$. There are an integer $j_0\geq1$, a sequence $\mathbf s=(s_k)_{k\geq1}$ in $[a,(m-1+q)/2]$ with $\inf_ks_k=a$, and integers $n_{j_0}<n_{j_0+1}<\cdots$ such that
\begin{equation}\label{lower-target-limits}
 \liminf_{n\to\infty}A_n(\mathbf s)=a,
 \qquad
 \limsup_{n\to\infty}A_n(\mathbf s)=m-1,\text{ and}
\end{equation}
\begin{equation}\label{lower-positive-excursion}
 A_{n_j}(\mathbf s)\geq m-1+n_j^{-1/4}
 \quad\text{for all }j\geq j_0.
\end{equation}
\end{lemma}

\begin{proof}
Put $u:=(m-1+q)/2$. Choose $j_0\geq1$ with $1/j_0<u-(m-1)$, and construct the sequence in stages $j\geq j_0$, starting from the empty sequence. At stage $j$, first append the value $a$ until the current length $n$ is at least $j^4$ and $A_n(\mathbf s)\leq a+1/j$. Then append the value $u$ at least once, until the first index $n_j$ for which $A_{n_j}(\mathbf s)\geq m-1+1/j$. This index exists because $u>m-1+1/j$. Since $A_{n_j-1}(\mathbf s)\leq m-1+1/j$, we have $A_{n_j}(\mathbf s)\leq m-1+1/j+u/n_j$. Since $n_j\geq j^4$, we have $A_{n_j}(\mathbf s)\geq m-1+1/j\geq m-1+n_j^{-1/4}$, which is \eqref{lower-positive-excursion}. Then proceed to stage $j+1$.

The averages $A_n(\mathbf s)$ do not increase while the value $a$ is appended, and they do not decrease while the value $u$ is appended. Hence their upper limit is the limit of their values at the indices $n_j$, which is $m-1$. Their lower limit is $a$, because every term is at least $a$ and the averages at the ends of the low blocks are at most $a+1/j$.
\end{proof}

\begin{figure}[htb]
\centering
\begin{tikzpicture}
\begin{axis}[name=avgplot,width=0.9\textwidth,height=4.2cm,xmode=log,xmin=1,xmax=8e6,ymin=0.4,ymax=1.6,
  ytick={0.5,1,1.5},yticklabels={$a=0.5$,$m-1=1$,$1.5$},xticklabels={},
  title={Running average: $A_n(\mathbf s)$},title style={font=\small},tick label style={font=\small}]
\addplot[PTblue,thick] coordinates {(1,0.5000) (2,0.5000) (3,0.5000) (4,0.5000) (5,0.5000) (6,0.5000) (7,0.5000) (8,0.5000) (9,0.5000) (10,0.5000) (11,0.5000) (12,0.5000) (13,0.5000) (14,0.5000) (15,0.5000) (16,0.5000) (17,0.5000) (18,0.5000) (19,0.5000) (20,0.5000) (21,0.5000) (23,0.5000) (24,0.5000) (26,0.5000) (27,0.5000) (29,0.5000) (31,0.5000) (32,0.5000) (34,0.5000) (37,0.5000) (39,0.5000) (41,0.5000) (44,0.5000) (46,0.5000) (49,0.5000) (52,0.5000) (55,0.5000) (59,0.5000) (62,0.5000) (66,0.5000) (70,0.5000) (74,0.5000) (79,0.5000) (81,0.5000) (83,0.5241) (88,0.5795) (94,0.6383) (99,0.6818) (105,0.7286) (112,0.7768) (118,0.8136) (126,0.8571) (133,0.8910) (141,0.9255) (150,0.9600) (159,0.9906) (168,1.0179) (178,1.0449) (189,1.0714) (200,1.0950) (212,1.1179) (225,1.1400) (239,1.1611) (253,1.1798) (268,1.1978) (284,1.2148) (301,1.2309) (320,1.2469) (339,1.2611) (359,1.2744) (381,1.2874) (404,1.2995) (428,1.3107) (454,1.3216) (481,1.3316) (486,1.3333) (510,1.2941) (540,1.2500) (573,1.2068) (607,1.1672) (644,1.1289) (682,1.0938) (723,1.0602) (767,1.0280) (813,0.9982) (861,0.9704) (913,0.9436) (968,0.9184) (1026,0.8947) (1088,0.8722) (1153,0.8513) (1222,0.8314) (1296,0.8125) (1373,0.7950) (1456,0.7782) (1543,0.7625) (1620,0.7500) (1636,0.7573) (1734,0.7993) (1838,0.8390) (1948,0.8763) (2065,0.9116) (2189,0.9450) (2321,0.9765) (2460,1.0061) (2607,1.0339) (2764,1.0604) (2930,1.0853) (3106,1.1088) (3292,1.1309) (3489,1.1518) (3699,1.1715) (3921,1.1901) (4156,1.2077) (4406,1.2242) (4670,1.2398) (4860,1.2500) (4950,1.2364) (5247,1.1947) (5562,1.1553) (5896,1.1182) (6249,1.0833) (6624,1.0503) (7022,1.0191) (7443,0.9897) (7890,0.9620) (8363,0.9358) (8865,0.9112) (9397,0.8879) (9961,0.8659) (10559,0.8452) (11192,0.8257) (11864,0.8072) (12576,0.7898) (13330,0.7734) (14130,0.7580) (14978,0.7434) (15877,0.7296) (16829,0.7166) (17839,0.7043) (18225,0.7000) (18909,0.7289) (20044,0.7726) (21247,0.8138) (22522,0.8526) (23873,0.8893) (25305,0.9238) (26824,0.9565) (28433,0.9872) (30139,1.0162) (31948,1.0436) (33864,1.0695) (35896,1.0938) (38050,1.1168) (40333,1.1385) (42753,1.1590) (45318,1.1783) (48038,1.1965) (48600,1.2000) (50920,1.1681) (53975,1.1303) (57214,1.0946) (60646,1.0610) (64285,1.0292) (68142,0.9993) (72231,0.9710) (76565,0.9443) (81159,0.9192) (86028,0.8955) (91190,0.8731) (96661,0.8520) (102461,0.8320) (108609,0.8132) (115125,0.7955) (122033,0.7788) (129355,0.7630) (137116,0.7481) (145343,0.7341) (154064,0.7208) (163308,0.7083) (173106,0.6965) (183493,0.6854) (194502,0.6749) (204120,0.6667) (206172,0.6750) (218543,0.7217) (231655,0.7657) (245555,0.8073) (260288,0.8465) (275905,0.8835) (292460,0.9184) (310007,0.9513) (328608,0.9824) (348324,1.0117) (369224,1.0393) (391377,1.0654) (414860,1.0900) (439752,1.1132) (466137,1.1351) (494105,1.1557) (510300,1.1667) (523751,1.1495) (555176,1.1128) (588487,1.0781) (623796,1.0454) (661224,1.0145) (700898,0.9854) (742951,0.9579) (787529,0.9320) (834780,0.9075) (884867,0.8845) (937959,0.8627) (994237,0.8422) (1053891,0.8228) (1117125,0.8045) (1184152,0.7873) (1255201,0.7710) (1330513,0.7557) (1410344,0.7412) (1494965,0.7276) (1584663,0.7147) (1679743,0.7025) (1780527,0.6911) (1887359,0.6803) (2000600,0.6700) (2120636,0.6604) (2247875,0.6513) (2381400,0.6429) (2382747,0.6433) (2525712,0.6918) (2677255,0.7376) (2837890,0.7807) (3008163,0.8214) (3188653,0.8599) (3379973,0.8961) (3582771,0.9303) (3797737,0.9625) (4025601,0.9929) (4267138,1.0216) (4523166,1.0487) (4794556,1.0743) (5082229,1.0984) (5387163,1.1211) (5710393,1.1425) (5715360,1.1429)};
\addplot[only marks,mark=*,mark size=1.5pt,PTmagenta] coordinates {(486,1.3333) (4860,1.2500) (48600,1.2000) (510300,1.1667) (5715360,1.1429)};
\node[above,font=\scriptsize] at (axis cs:486,1.3333) {$j=3$};
\node[above,font=\scriptsize] at (axis cs:4860,1.2500) {$j=4$};
\node[above,font=\scriptsize] at (axis cs:48600,1.2000) {$j=5$};
\node[above,font=\scriptsize] at (axis cs:510300,1.1667) {$j=6$};
\node[above left,font=\scriptsize] at (axis cs:5715360,1.1429) {$j=7$};
\addplot[dashed,gray,domain=1:8e6] {1};
\addplot[dashed,gray,domain=1:8e6] {0.5};
\end{axis}
\begin{axis}[at={(avgplot.below south west)},anchor=north west,yshift=-0.85cm,width=0.9\textwidth,height=4.2cm,xmode=log,xmin=1,xmax=8e6,ymin=-15,ymax=8,
  ytick={-10,-5,0,1,5},yticklabels={$-10$,$-5$,,$1$,$5$},xlabel={$n$},
  title={Normalized excess: $n^{1/4}(A_n(\mathbf s)-(m-1))$},title style={font=\small},tick label style={font=\small},xlabel style={font=\small}]
\addplot[PTteal,thick] coordinates {(1,-0.5000) (2,-0.5946) (3,-0.6580) (4,-0.7071) (5,-0.7477) (6,-0.7825) (7,-0.8133) (8,-0.8409) (9,-0.8660) (10,-0.8891) (11,-0.9106) (12,-0.9306) (13,-0.9494) (14,-0.9672) (15,-0.9840) (16,-1.0000) (17,-1.0153) (18,-1.0299) (19,-1.0439) (20,-1.0574) (21,-1.0703) (23,-1.0950) (24,-1.1067) (26,-1.1291) (27,-1.1398) (29,-1.1603) (31,-1.1798) (32,-1.1892) (34,-1.2074) (37,-1.2332) (39,-1.2495) (41,-1.2652) (44,-1.2878) (46,-1.3021) (49,-1.3229) (52,-1.3427) (55,-1.3616) (59,-1.3857) (62,-1.4030) (66,-1.4251) (70,-1.4463) (74,-1.4665) (79,-1.4907) (81,-1.5000) (83,-1.4364) (88,-1.2878) (94,-1.1262) (99,-1.0037) (105,-0.8689) (112,-0.7262) (118,-0.6145) (126,-0.4786) (133,-0.3702) (141,-0.2566) (150,-0.1400) (159,-0.0335) (168,0.0643) (178,0.1642) (189,0.2648) (200,0.3573) (212,0.4500) (225,0.5422) (239,0.6334) (253,0.7173) (268,0.8002) (284,0.8817) (301,0.9617) (320,1.0442) (339,1.1202) (359,1.1943) (381,1.2698) (404,1.3428) (428,1.4134) (454,1.4844) (481,1.5529) (486,1.5651) (510,1.3977) (540,1.2051) (573,1.0118) (607,0.8300) (644,0.6493) (682,0.4796) (723,0.3120) (767,0.1475) (813,-0.0099) (861,-0.1604) (913,-0.3101) (968,-0.4552) (1026,-0.5957) (1088,-0.7337) (1153,-0.8667) (1222,-0.9967) (1296,-1.1250) (1373,-1.2480) (1456,-1.3703) (1543,-1.4887) (1620,-1.5861) (1636,-1.5433) (1734,-1.2951) (1838,-1.0545) (1948,-0.8219) (2065,-0.5958) (2189,-0.3765) (2321,-0.1630) (2460,0.0429) (2607,0.2426) (2764,0.4381) (2930,0.6278) (3106,0.8124) (3292,0.9917) (3489,1.1664) (3699,1.3377) (3921,1.5045) (4156,1.6673) (4406,1.8269) (4670,1.9826) (4860,2.0874) (4950,1.9826) (5247,1.6569) (5562,1.3415) (5896,1.0359) (6249,0.7406) (6624,0.4535) (7022,0.1747) (7443,-0.0955) (7890,-0.3584) (8363,-0.6135) (8865,-0.8620) (9397,-1.1038) (9961,-1.3394) (10559,-1.5692) (11192,-1.7930) (11864,-2.0118) (12576,-2.2256) (13330,-2.4344) (14130,-2.6389) (14978,-2.8392) (15877,-3.0355) (16829,-3.2280) (17839,-3.4171) (18225,-3.4857) (18909,-3.1786) (20044,-2.7057) (21247,-2.2482) (22522,-1.8053) (23873,-1.3764) (25305,-0.9607) (26824,-0.5573) (28433,-0.1660) (30139,0.2140) (31948,0.5834) (33864,0.9422) (35896,1.2915) (38050,1.6316) (40333,1.9629) (42753,2.2859) (45318,2.6011) (48038,2.9090) (48600,2.9695) (50920,2.5253) (53975,1.9859) (57214,1.4632) (60646,0.9566) (64285,0.4650) (68142,-0.0121) (72231,-0.4756) (76565,-0.9261) (81159,-1.3642) (86028,-1.7905) (91190,-2.2058) (96661,-2.6105) (102461,-3.0052) (108609,-3.3905) (115125,-3.7668) (122033,-4.1348) (129355,-4.4947) (137116,-4.8471) (145343,-5.1924) (154064,-5.5311) (163308,-5.8636) (173106,-6.1901) (183493,-6.5112) (194502,-6.8271) (204120,-7.0852) (206172,-6.9262) (218543,-6.0180) (231655,-5.1398) (245555,-4.2900) (260288,-3.4673) (275905,-2.6704) (292460,-1.8980) (310007,-1.1491) (328608,-0.4223) (348324,0.2833) (369224,0.9689) (391377,1.6353) (414860,2.2837) (439752,2.9148) (466137,3.5297) (494105,4.1291) (510300,4.4546) (523751,4.0230) (555176,3.0785) (588487,2.1629) (623796,1.2751) (661224,0.4135) (700898,-0.4231) (742951,-1.2359) (787529,-2.0262) (834780,-2.7950) (884867,-3.5435) (937959,-4.2728) (994237,-4.9838) (1053891,-5.6775) (1117125,-6.3548) (1184152,-7.0167) (1255201,-7.6639) (1330513,-8.2974) (1410344,-8.9180) (1494965,-9.5263) (1584663,-10.1231) (1679743,-10.7091) (1780527,-11.2850) (1887359,-11.8515) (2000600,-12.4091) (2120636,-12.9585) (2247875,-13.5002) (2381400,-14.0298) (2382747,-14.0127) (2525712,-12.2853) (2677255,-10.6151) (2837890,-8.9996) (3008163,-7.4361) (3188653,-5.9222) (3379973,-4.4554) (3582771,-3.0336) (3797737,-1.6545) (4025601,-0.3160) (4267138,0.9838) (4523166,2.2470) (4794556,3.4752) (5082229,4.6704) (5387163,5.8342) (5710393,6.9682) (5715360,6.9849)};
\addplot[only marks,mark=*,mark size=1.5pt,PTmagenta] coordinates {(486,1.5651) (4860,2.0874) (48600,2.9695) (510300,4.4546) (5715360,6.9849)};
\addplot[dashed,gray,domain=1:8e6] {1};
\end{axis}
\end{tikzpicture}
\caption{The construction in the proof of Lemma \ref{l:lower-targets} with $a=0.5$, $m=2$, and $q=2$, for which $u=1.5$ and $j_0=3$, shown for the stages $3\leq j\leq 7$. Top: the averages fall to at most $a+1/j$ on each low block and first reach $m-1+1/j$ at the marked times $n_j$. Bottom: the normalized excess is at least $1$ at every $n_j$.}\label{fig:lower-targets}
\end{figure}
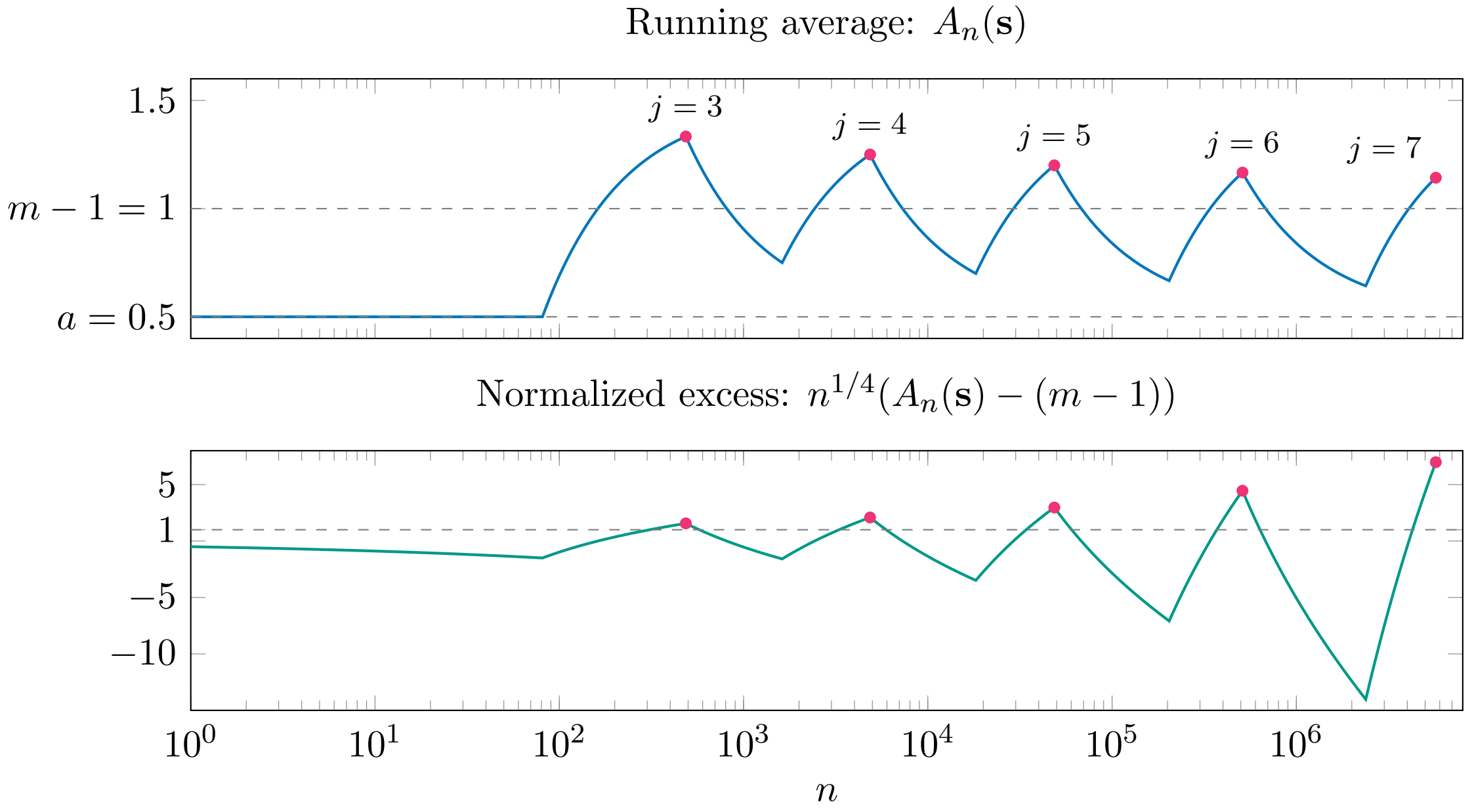

Figure \ref{fig:lower-targets} shows one run of the construction in the case $a<m-1$.

The next lemma converts \eqref{lower-positive-excursion} into vanishing lower $(m-1)$-density.

\begin{lemma}[vanishing lower density]\label{l:lower-fluctuation}
Let $m\geq2$ be an integer and let $0<t_{\min}\leq q/2$. Let $\muse$ be given by Definition \ref{d:endpoint-measure} on a system of $b$-adic cubes without jumps, with the convention of \S\ref{ss:intermediate}, and with target entropies $\mathbf s=(s_k)_{k\geq1}$. Suppose that for some integer $j_0\geq1$ and integers $n_{j_0}<n_{j_0+1}<\cdots$,
\begin{equation}\label{fluctuation-hypothesis}
 A_{n_j}(\mathbf s)\geq m-1+n_j^{-1/4}\quad\text{for all }j\geq j_0.
\end{equation}
Then
\begin{equation}\label{zero-lower-d-density}
 \liminf_{r\downarrow0}\frac{\muse(B(x,r))}{r^{m-1}}=0\quad\text{at $\muse$-a.e.~}x\in\XX.
\end{equation}
\end{lemma}

\begin{proof}
Fix $Q_0\in\Delta_0$, normalize $\muse(Q_0)=1$, and use the notation introduced after Definition \ref{d:target}. By Remark \ref{r:nested}, every child has mass at least $\kappa_0$ times the mass of its parent. Hence $0\leq X_k\leq\log_b\kappa_0^{-1}$. Since $0\leq s_k\leq q+1$ by Lemma \ref{l:second-moment}, the sums $M_n:=\sum_{k=1}^n(X_k-s_k)$ form a martingale with respect to $(\mathscr F_n)_{n\geq0}$ whose increments are bounded by $L:=q+1+\log_b\kappa_0^{-1}$. Azuma's inequality (see e.g.~\cite[top of p.~96]{concentration-survey}) gives
\begin{equation}\label{lower-azuma}
 \PP\{M_n\leq-t\}\leq\exp\left(-\frac{t^2}{2L^2n}\right)\quad\text{for all }t>0\text{ and }n\geq1.
\end{equation}
Let $j\geq j_0$. By \eqref{fluctuation-hypothesis}, $\sum_{k=1}^{n_j}s_k\geq(m-1)n_j+n_j^{3/4}$. Hence \eqref{lower-azuma} with $n=n_j$ and $t=\frac12n_j^{3/4}$ gives
\[
 \PP\left\{\sum_{k=1}^{n_j}X_k<(m-1)n_j+\tfrac12n_j^{3/4}\right\}\leq\PP\left\{M_{n_j}\leq-\tfrac12n_j^{3/4}\right\}\leq\exp\left(-\frac{n_j^{1/2}}{8L^2}\right).
\]
The right side is summable over $j\geq j_0$, because the integers $n_j$ are distinct. By the Borel--Cantelli lemma, at $\muse$-a.e.~$x\in Q_0$, we have $\sum_{k=1}^{n_j}X_k(x)\geq(m-1)n_j+\frac12n_j^{3/4}$ for all sufficiently large $j$. Since $\muse(Q_n(x))=b^{-\sum_{k=1}^nX_k(x)}$ by \eqref{endpoint-weighted}, it follows that
\begin{equation}\label{lower-cube-density}
 \lim_{j\to\infty}\frac{\muse(Q_{n_j}(x))}{b^{-(m-1)n_j}}=0\quad\text{at $\muse$-a.e.~}x\in Q_0.
\end{equation}

Let $x\in Q_0$ and $n\geq0$, and put $Q:=Q_n(x)$. Then $B(x,b^{-n})\subset B(x_Q,(C_b+1)b^{-n})$ and $B(x_Q,\tfrac12c_bb^{-n})\subset Q$. Since $c_b\geq 2/5$ and $C_b\leq 6/5$, we have $C_b+1\leq 8c_b$, and four applications of \eqref{doubling}, which holds by Lemma \ref{l:endpoint-doubling}, give $\muse(B(x,b^{-n}))\leq D_{\mathsf{db}}^4\,\muse(Q_n(x))$. Together with \eqref{lower-cube-density}, this yields
\[
 \lim_{j\to\infty}\frac{\muse(B(x,b^{-n_j}))}{(b^{-n_j})^{m-1}}=0\quad\text{at $\muse$-a.e.~}x\in Q_0,
\]
which implies \eqref{zero-lower-d-density} at $\muse$-a.e.~$x\in Q_0$. Since $\Delta_0$ is countable, \eqref{zero-lower-d-density} holds at $\muse$-a.e.~$x\in\XX$.
\end{proof}

We can now complete the proof of Theorem \ref{t:main} at the lower endpoint.

\begin{proof}[Proof of Theorem \ref{t:main} when $s_P=m-1$ and $s_H>0$]
Let $j_0$, $\mathbf s=(s_k)_{k\geq1}$, and $(n_j)_{j\geq j_0}$ be given by Lemma \ref{l:lower-targets} with $a=s_H$. Then $s_H\leq s_k\leq(m-1+q)/2$ for all $k\geq1$. Put $t_{\min}:=\min\{s_H,q/2\}$. To choose $b$ independently of $m$, note that $q-m+1\geq q+1-\lceil q\rceil>0$, since $m-1$ is an integer and $m-1<q$. Choose $b\geq\max\{47,C_\XX\}$, depending only on $\XX$, such that $\log_bC_{\mathrm{ent}}\leq\min\{q/4,\,(q+1-\lceil q\rceil)/2\}$, with $C_{\mathrm{ent}}$ as in \eqref{K-and-C2}. Then
\[
 t_{\min}\leq s_k\leq q-(q-m+1)/2\leq q-\log_bC_{\mathrm{ent}}\quad\text{for all }k\geq1.
\]

Let $(\Delta_k)_{k\in\ZZ}$ be a system of $b$-adic cubes without jumps, with the convention of \S\ref{ss:intermediate}, and let $\muse$ be given by Definition \ref{d:endpoint-measure} with target entropies $\mathbf s$, which exists by Lemma \ref{l:entropy-allocation}. Every cube has positive mass, and $\muse$ has full support. By Lemma \ref{l:endpoint-doubling}, Theorem \ref{t:dimension}, and \eqref{lower-target-limits}, $\muse$ satisfies \eqref{doubling} with $D_{\mathsf{db}}$ depending only on $\XX$ and $s_H$, and $\muse$ has lower local dimension $s_H$ and upper local dimension $m-1$ at $\muse$-a.e.~point. Since $m-1<m$, Corollary \ref{c:bk} shows that $\muse$ is $m$-rectifiable. By \eqref{lower-positive-excursion} and Lemma \ref{l:lower-fluctuation}, the lower $(m-1)$-density of $\muse$ vanishes at $\muse$-a.e.~point. Hence $\muse$ is purely $(m-1)$-unrectifiable by Remark \ref{r:pure-unrectifiable}.
\end{proof}

\section{\texorpdfstring{Zero lower local dimension: $s_H=0$}{Zero lower local dimension: s\textunderscore H = 0}}\label{s:zero}

It remains to prove Theorem \ref{t:main} when $s_H=0$. We treat this case separately, because the theorem does not ask for the measure to be doubling. This removes the only reason to keep the child weights uniformly away from zero. We therefore use the central weights of Definition \ref{d:central}, whose target entropies may tend to zero, and which require no intermediate levels. When $s_P<q$, we use $b$-adic cubes without jumps, and when $s_P=q$, we use $b$-adic cubes with sparse jumps.

\subsection{Fixed scaling ratio}\label{ss:zero-fixed}

In this subsection, $s_P<q$. Put $u_0:=(s_P+q)/2$, choose $b\geq\max\{47,C_\XX\}$ such that $\log_bC_{\mathrm{ent}}\leq\min\{q/4,(q-s_P)/2\}$, and let $(\Delta_k)_{k\in\ZZ}$ be a system of $b$-adic cubes without jumps. Then $u_0\leq q-\log_bC_{\mathrm{ent}}$, and every target in $(0,u_0]$ satisfies \eqref{direct-range}. Given target entropies $0<s_k\leq u_0$, let $\muse$ be given by Definition \ref{d:central}. By Lemma \ref{l:central-existence}, $\muse$ is a Radon measure with full support and target entropies $\mathbf s=(s_k)_{k\geq1}$. By Theorem \ref{t:dimension}, the lower and upper local dimensions of $\muse$ are the lower and upper limits of $A_n(\mathbf s)$ at $\muse$-a.e.~point.

If $m-1<s_P<m$, alternate sufficiently long blocks on which $s_k=\min\{u_0,1/k\}$ with blocks on which $s_k=s_P$. Then $\liminf_nA_n(\mathbf s)=0$ and $\limsup_nA_n(\mathbf s)=s_P$. Corollary \ref{c:bk} gives $m$-rectifiability, and Remark \ref{r:pure-unrectifiable} gives pure $(m-1)$-unrectifiability, because $s_P>m-1$.

If $s_P=m<q$, alternate sufficiently long blocks on which $s_k=1/\log_2(k+4)$ with blocks on which $s_k=m-1/\log_2(k+4)$. Then $\liminf_nA_n(\mathbf s)=0$ and $\limsup_nA_n(\mathbf s)=m$, and $m-s_k\geq1/\log_2(k+4)$ for all $k$. As in \S\ref{ss:upper-fixed}, \eqref{square-hypothesis} holds with $c_0=1/2$. By Proposition \ref{p:square}, $\muse$ is $m$-rectifiable, and by Remark \ref{r:pure-unrectifiable}, it is purely $(m-1)$-unrectifiable, because $s_P>m-1$.

It remains, among the cases with $s_P<q$, to treat $s_P=m-1$. Suppose first that $m\geq2$. Two changes are needed in the argument of Lemma \ref{l:lower-fluctuation}. First, the functions $X_k-s_k$ need no longer be uniformly bounded. In place of Azuma's inequality, we use the second-moment estimate of Lemma \ref{l:second-moment} and Chebyshev's inequality along a sequence of generations for which the resulting bounds are summable. Second, to pass from cubes to balls without doubling, after every generation $n_j$ at which $A_{n_j}(\mathbf s)$ exceeds $m-1$ by at least $n_j^{-1/4}$, we insert one generation of very low entropy. We choose this entropy so that the proportions of mass outside the central children are summable in $j$. Every point in the central child of a cube $Q\in\Delta_{n_j}$ is the center of a ball of radius $9b^{-(n_j+1)}$ contained in $Q$, which lets us transfer the small masses of cubes to small masses of balls at almost every point.

\begin{lemma}[vanishing lower density without doubling]\label{l:zero-fluctuation}
Let $m\geq2$ be an integer, let $(\Delta_k)_{k\in\ZZ}$ be a system of $b$-adic cubes without jumps, and let $\muse$ be given by Definition \ref{d:central} with target entropies $\mathbf s=(s_k)_{k\geq1}$. Suppose that for some integer $j_0\geq1$ and integers $n_{j_0}<n_{j_0+1}<\cdots$,
\begin{equation}\label{zero-fluctuation-excess}
 A_{n_j}(\mathbf s)\geq m-1+n_j^{-1/4}\quad\text{for all }j\geq j_0,
\end{equation}
\begin{equation}\label{zero-fluctuation-summability}
 \sum_{j=j_0}^\infty n_j^{-1/2}<\infty,\text{ and}
\end{equation}
\begin{equation}\label{zero-buffer-central}
 \muse(Q^\downarrow)\geq(1-2^{-j})\muse(Q)\quad\text{for all }Q\in\Delta_{n_j}\text{ and }j\geq j_0.
\end{equation}
Then
\begin{equation}\label{zero-lower-density-final}
 \liminf_{r\downarrow0}\frac{\muse(B(x,r))}{r^{m-1}}=0\quad\text{at $\muse$-a.e.~}x\in\XX.
\end{equation}
\end{lemma}

\begin{proof}
Fix $Q_0\in\Delta_0$, normalize $\muse(Q_0)=1$, and use the notation introduced after Definition \ref{d:target}. Put $Y_k:=X_k-s_k$. As in the proof of Lemma \ref{l:entropy-good-descendants}, the functions $Y_k$ are pairwise orthogonal on $Q_0$, and Lemma \ref{l:second-moment} gives
\[
 \EE\Big(\sum_{k=1}^nY_k\Big)^2=\sum_{k=1}^n\EE(Y_k^2)\leq C_{\mathrm{var}}n\quad\text{for all }n\geq1.
\]
Let $j\geq j_0$. By \eqref{zero-fluctuation-excess}, $\sum_{k=1}^{n_j}s_k\geq(m-1)n_j+n_j^{3/4}$, and Chebyshev's inequality yields
\[
 \PP\left\{\sum_{k=1}^{n_j}X_k<(m-1)n_j+\tfrac12n_j^{3/4}\right\}\leq\PP\left\{\sum_{k=1}^{n_j}Y_k<-\tfrac12n_j^{3/4}\right\}\leq4C_{\mathrm{var}}n_j^{-1/2}.
\]
By \eqref{zero-fluctuation-summability} and the Borel--Cantelli lemma, at $\muse$-a.e.~$x\in Q_0$, we have $\sum_{k=1}^{n_j}X_k(x)\geq(m-1)n_j+\frac12n_j^{3/4}$ for all sufficiently large $j$. Since $\muse(Q_n(x))=b^{-\sum_{k=1}^nX_k(x)}$ by \eqref{endpoint-weighted}, it follows that
\begin{equation}\label{zero-cube-density}
 \lim_{j\to\infty}\frac{\muse(Q_{n_j}(x))}{b^{-(m-1)n_j}}=0\quad\text{at $\muse$-a.e.~}x\in Q_0.
\end{equation}
By \eqref{zero-buffer-central}, the set of $x\in Q_0$ with $Q_{n_j+1}(x)\neq Q_{n_j}(x)^\downarrow$ has measure at most $2^{-j}$. By the Borel--Cantelli lemma again, at $\muse$-a.e.~$x\in Q_0$, we have $Q_{n_j+1}(x)=Q_{n_j}(x)^\downarrow$ for all sufficiently large $j$. For such $x$ and $j$, Lemma \ref{l:vAR-count} gives $\gap(Q_{n_j+1}(x),\XX\setminus Q_{n_j}(x))>9b^{-(n_j+1)}$, and hence $B(x,9b^{-(n_j+1)})\subset Q_{n_j}(x)$. Therefore,
\[
 \frac{\muse(B(x,9b^{-(n_j+1)}))}{(9b^{-(n_j+1)})^{m-1}}\leq\left(\frac b9\right)^{m-1}\frac{\muse(Q_{n_j}(x))}{b^{-(m-1)n_j}},
\]
and the right side tends to $0$ as $j\to\infty$ by \eqref{zero-cube-density}. This proves \eqref{zero-lower-density-final} at $\muse$-a.e.~$x\in Q_0$. Since $\Delta_0$ is countable, \eqref{zero-lower-density-final} holds at $\muse$-a.e.~$x\in\XX$.
\end{proof}

\begin{lemma}[targets at the lower packing endpoint]\label{l:zero-lower-targets}
Let $m\geq2$ be an integer with $m-1<q$, and put $u_0:=(m-1+q)/2$. Let $(\Delta_k)_{k\in\ZZ}$ be a system of $b$-adic cubes without jumps with $u_0\leq q-\log_bC_{\mathrm{ent}}$. There are an integer $j_0\geq1$, target entropies $0<s_k\leq u_0$, and integers $n_{j_0}<n_{j_0+1}<\cdots$ such that $\liminf_nA_n(\mathbf s)=0$ and $\limsup_nA_n(\mathbf s)=m-1$, \eqref{zero-fluctuation-excess} and \eqref{zero-fluctuation-summability} hold, and the measure $\muse$ given by Definition \ref{d:central} with target entropies $\mathbf s=(s_k)_{k\geq1}$ satisfies \eqref{zero-buffer-central}.
\end{lemma}

\begin{proof}
Choose $j_0\geq1$ with $1/j_0<\min\{u_0-(m-1),\,m-1\}$, and construct the sequence in stages $j\geq j_0$, starting from the empty sequence. At stage $j$, first append the value $1/j$ until the current length $n$ is at least $j^4$ and $A_n(\mathbf s)\leq2/j$. Then append the value $u_0$ until the first index $n_j$ for which $A_{n_j}(\mathbf s)\geq m-1+1/j$. This index exists because $u_0>m-1+1/j$, and it is preceded by at least one copy of $u_0$ because $2/j<m-1+1/j$. By minimality, $A_{n_j}(\mathbf s)\leq m-1+1/j+u_0/n_j$. Since $n_j\geq j^4$, we have $A_{n_j}(\mathbf s)\geq m-1+1/j\geq m-1+n_j^{-1/4}$, which is \eqref{zero-fluctuation-excess}. Also $\sum_jn_j^{-1/2}\leq\sum_jj^{-2}$, which gives \eqref{zero-fluctuation-summability}.

After generation $n_j$, append one additional target $s_{n_j+1}>0$ that is smaller than $u_0$ and than $A_{n_j}(\mathbf s)$, and small enough that the number $v$ with $h_{b,N}(v)=s_{n_j+1}$ satisfies $v\leq 2^{-j}$ for every integer $2\leq N\leq C_\XX b^q$. This is possible by Lemma \ref{l:central-entropy}, because only finitely many values of $N$ occur, and it gives \eqref{zero-buffer-central} by \eqref{child-range}. Then proceed to stage $j+1$.

The averages $A_n(\mathbf s)$ do not increase while the values $1/j$ and the additional targets are appended, and they do not decrease while the value $u_0$ is appended. Hence their upper limit is the limit of their values at the indices $n_j$, which is $m-1$. Their lower limit is $0$, because all targets are positive and the averages at the ends of the low blocks are at most $2/j$.
\end{proof}

When $m\geq2$ and $s_P=m-1$, the number $u_0$ fixed at the beginning of this subsection is $(m-1+q)/2$. Let $\mathbf s$ be the target entropies of Lemma \ref{l:zero-lower-targets}, and let $\muse$ be given by Definition \ref{d:central}. Then $\muse$ has lower local dimension $0$ and upper local dimension $m-1$ at $\muse$-a.e.~point, and by Lemma \ref{l:zero-fluctuation}, the lower $(m-1)$-density of $\muse$ vanishes at $\muse$-a.e.~point. As in \S\ref{s:lower}, Corollary \ref{c:bk} and Remark \ref{r:pure-unrectifiable} show that $\muse$ is $m$-rectifiable and purely $(m-1)$-unrectifiable.

It remains to treat $m=1$ and $s_P=0$, in which case $u_0=q/2$. Put $s_k:=u_0$ when $k$ is a perfect square and $s_k:=\min\{u_0,1/k\}$ otherwise. The perfect squares have density zero, and $s_k$ tends to $0$ along the remaining indices. Hence $\lim_{n\to\infty}A_n(\mathbf s)=0$, the measure $\muse$ has lower and upper local dimension $0$ at $\muse$-a.e.~point, and Corollary \ref{c:bk} shows that $\muse$ is $1$-rectifiable. To see that $\muse$ has no atoms, let $k$ be a perfect square and let $P\in\Delta_{k-1}$. By \eqref{child-range}, $2\leq N_P\leq C_\XX b^q$. Hence the numbers $v_P$ with $h_{b,N_P}(v_P)=u_0$ take only finitely many values, and we let $v_*>0$ be the least of them. The central weights with target $u_0$ on $\Child(P)$ are $1-v_P\leq1-v_*$ and $v_P/(N_P-1)\leq1/N_P\leq1/2$. Thus, $\muse(R)\leq\theta\,\muse(P)$ for all $R\in\Child(P)$, where $\theta:=\max\{1-v_*,1/2\}<1$. It follows that
\[
 \muse(\{x\})\leq\muse(Q_n(x))\leq\theta^{\lfloor\sqrt n\rfloor}\muse(Q_0(x))\quad\text{for all }x\in\XX\text{ and }n\geq1.
\]
Therefore, $\muse(\{x\})=0$ for every $x\in\XX$, which is \eqref{unrectifiable} when $m=1$.

\subsection{\texorpdfstring{The ambient upper endpoint $s_P=q$}{The ambient upper endpoint s\textunderscore P = q}}\label{ss:zero-q}

Choose $b\geq\max\{47,C_\XX\}$, depending only on $\XX$, such that \eqref{b-entropy} holds, and let $(\Delta_k)_{k\in\ZZ}$ be a system of $b$-adic cubes with sparse jumps. Alternate sufficiently long blocks of generations with values $s_k=2\log_{\xi_k}C_{\mathrm{ent}}$ and $s_k=q-2\log_{\xi_k}C_{\mathrm{ent}}$. These values are positive and satisfy \eqref{direct-range}, because $C_{\mathrm{ent}}>1$ and $3\log_{\xi_k}C_{\mathrm{ent}}\leq3q/4$ by \eqref{b-entropy}. Since $\log_{\xi_k}C_{\mathrm{ent}}\to0$ as $k\to\infty$, the blocks may be chosen so that $A_n(\mathbf s)$ has lower limit $0$ and upper limit $q$. Let $\muse$ be given by Definition \ref{d:central}. By Lemma \ref{l:central-existence} and Theorem \ref{t:dimension}, $\muse$ is a Radon measure with full support and target entropies $\mathbf s=(s_k)_{k\geq1}$, and it has lower local dimension $0$ and upper local dimension $q$ at $\muse$-a.e.~point.

If $q<m$, then Corollary \ref{c:bk} gives $m$-rectifiability, and Remark \ref{r:pure-unrectifiable} gives pure $(m-1)$-unrectifiability, because $q>m-1$. If $q=m$, then by \eqref{b-entropy} the two kinds of generations satisfy $(q-s_k)\log_2\xi_k=q\log_2\xi_k-2\log_2C_{\mathrm{ent}}\geq 2\log_2C_{\mathrm{ent}}$ and $(q-s_k)\log_2\xi_k=2\log_2C_{\mathrm{ent}}$, respectively. As in \S\ref{ss:upper-q}, \eqref{square-hypothesis} holds with $m=q$ and $c_0=2\log_2C_{\mathrm{ent}}$. By Proposition \ref{p:square}, $\muse$ is $q$-rectifiable. Since the upper local dimension of $\muse$ is $q=m>m-1$ at $\muse$-a.e.~point, Remark \ref{r:pure-unrectifiable} shows that $\muse$ is purely $(m-1)$-unrectifiable.

\begin{proof}[Completion of the proof of Theorem \ref{t:main} when $s_H=0$]
If $s_P<q$, the constructions in \S\ref{ss:zero-fixed} treat the cases $m-1<s_P<m$, $s_P=m$, $s_P=m-1>0$, and $m=1$ with $s_P=0$. If $s_P=q$, the construction in \S\ref{ss:zero-q} treats the cases $q<m$ and $q=m$. In every case, the measure has full support and the required lower and upper local dimensions. Theorem \ref{t:main} does not assert doubling when $s_H=0$, and the proof is complete.
\end{proof}

\providecommand{\bysame}{\leavevmode\hbox to3em{\hrulefill}\thinspace}
\providecommand{\MR}{\relax\ifhmode\unskip\space\fi MR }
\providecommand{\MRhref}[2]{%
  \href{http://www.ams.org/mathscinet-getitem?mr=#1}{#2}
}
\providecommand{\href}[2]{#2}

\end{document}